\documentclass[11pt]{article}
\usepackage{amsfonts, amsthm, amssymb, amsmath, stmaryrd}
\usepackage{mathrsfs,array}
\usepackage{xy}
\usepackage{hyperref}
\usepackage{pst-all,pst-node,pst-tree,graphicx}
\input xy
\xyoption{all}

\def\from{\colon}

\DeclareMathOperator{\nr}{nr}

\def\C{\mathbf{C}}

\def\Q{\mathbf{Q}}
\def\R{\mathcal{R}}

\def\Z{\mathbf{Z}}
\def\FF{\mathbf{F}}
\def\OO{\mathcal{O}}
\def\A{\mathcal{A}}

\def\K{\mathcal{K}}
\def\LL{\mathcal{L}}
\def\M{\mathcal{M}}
\def\MM{\mathcal{M}}

\def\gm{\mathfrak{m}}
\def\gp{\mathfrak{p}}

\def\isom{\cong}

\newcommand{\powerseries}[1]{\llbracket #1 \rrbracket}

\newcommand{\injects}{\hookrightarrow}

\DeclareMathOperator{\End}{End}
\DeclareMathOperator{\Aut}{Aut}

\DeclareMathOperator{\Spec}{Spec}

\DeclareMathOperator{\Spa}{Spa}

\DeclareMathOperator{\Hom}{Hom}

\DeclareMathOperator{\Spf}{Spf}
\DeclareMathOperator{\univ}{univ}
\DeclareMathOperator{\GL}{GL}
\DeclareMathOperator{\SL}{SL}

\DeclareMathOperator{\JL}{JL}
\DeclareMathOperator{\LLC}{LLC}
\DeclareMathOperator{\Ind}{Ind}

\DeclareMathOperator{\perf}{perf}

\DeclareMathOperator{\ad}{ad}

\DeclareMathOperator{\W}{\mathscr{W}}

\DeclareMathOperator{\St}{St}

\DeclareMathOperator{\tr}{tr}

\DeclareMathOperator{\cl}{cl}
\DeclareMathOperator{\Ends}{Ends}
\DeclareMathOperator{\nonCM}{non-CM}
\newcommand{\abs}[1]{\left\lvert #1 \right\rvert}
\newcommand{\class}[1]{\left< #1 \right>}
\newcommand{\set}[1]{\left\{ #1 \right\}}

\newcommand{\floor}[1]{\left\lfloor #1 \right\rfloor}
\newcommand{\tbt}[4]{\left(\begin{matrix}#1 & #2\\#3 & #4\end{matrix}\right)}
\newcommand{\tbo}[2]{\left(\begin{matrix}#1 \\ #2 \end{matrix}\right)}
\newcommand{\tatealgebra}{\class}

\DeclareMathOperator{\Nil}{Nil}
\DeclareMathOperator{\redu}{red}
\DeclareMathOperator{\Ig}{Ig}

\def\Vect{\mathbf{Vect}}

\def\Mod{\mathbf{Mod}}
\def\Alg{\mathbf{Alg}}

\def\Vect{\mathbf{Vect}}

\def\gp{\mathfrak{p}}
\def\gP{\mathfrak{P}}

\def\YY{\mathcal{Y}}
\def\ZZ{\mathcal{Z}}
\def\T{\mathcal{T}}
\def\X{\mathbf{X}}
\def\Y{\mathbf{Y}}
\def\x{\mathbf{x}}
\def\y{\mathbf{y}}
\def\t{\mathbf{t}}
\def\L{\mathcal{L}}

\numberwithin{equation}{subsection}
\newtheorem{Theorem}{Theorem}[subsection]
\numberwithin{Theorem}{subsection}
\newtheorem{lemma}[Theorem]{Lemma}
\newtheorem{Cor}[Theorem]{Corollary}
\newtheorem{prop}[Theorem]{Proposition}

\theoremstyle{definition}
\newtheorem{defn}[Theorem]{Definition}

\newtheorem{exmp}[Theorem]{Example}

\newtheorem{rmk}[Theorem]{Remark}

\def\bluedot{\pscircle*[linecolor=blue](0,0){.1}}
\def\blackdot{\psdots[dotstyle=o](0,0)}
\def\greendot{\pscircle*[linecolor=green](0,0){.1}}
\def\reddot{\pscircle*[linecolor=red](0,0){.1}}

\newcommand{\greekcross}{
\rput(0,0){\rnode{0}{\bluedot}}
\rput(0,2){\rnode{01p}{\blackdot}}
\rput(2,0){\rnode{02p}{\blackdot}}
\rput(-2,0){\rnode{03p}{\blackdot}}
\rput(0,-2){\rnode{04p}{\blackdot}}
\ncline {0}{01p}
\ncline {0}{02p}
\ncline {0}{03p}
\ncline {0}{04p}
}

\newcommand{\romancrosssmall}{
\rput(0,1){\rnode{0}{\bluedot}}
\rput(-.5,1){\rnode{01p}{\blackdot}}
\rput(0,1.5){\rnode{02p}{\blackdot}}
\rput(.5,1){\rnode{03p}{\blackdot}}
\rput(0,0){\rnode{04p}{\blackdot}}
\ncline {0}{01p}
\ncline {0}{02p}
\ncline {0}{03p}
\ncline {0}{04p}
}

\newcommand{\romancross}{
\rput(0,2){\rnode{A0}{\bluedot}}
\rput{90}(-1,2){\rnode{A01p}{\romancrosssmall}}
\rput{0}(0,3){\rnode{A02p}{\romancrosssmall}}
\rput{270}(1,2){\rnode{A03p}{\romancrosssmall}}
\rput(0,0){\rnode{A04p}{\blackdot}}
\ncline {A0}{A01p}
\ncline {A0}{A02p}
\ncline {A0}{A03p}
\ncline {A0}{A04p}
}

\newcommand{\depthzerotree}{
\pspicture(-3,-3)(3,3)
\psset{unit=.5, nodesep=2pt}
\rput(0,0){\greekcross}
\rput{0}(0,2){\romancross}
\rput{90}(-2,0){\romancross}
\rput{180}(0,-2){\romancross}
\rput{270}(2,0){\romancross}
\endpspicture
}

\newcommand\smallunramifiedtree{
\rput(0,0){\rnode{0g}{\greendot}}
\rput(-.10,1){\rnode{01g}{\greendot}}
\rput(-.05,1){\rnode{02g}{\greendot}}
\rput(.05,1){\rnode{03g}{\greendot}}
\rput(.10,1){\rnode{04g}{\greendot}}
\ncline {0g}{01g}
\ncline {0g}{02g}
\ncline {0g}{03g}
\ncline {0g}{04g}
}

\newcommand\unramifiedtree{
\rput(0,0){\rnode{0a}{\bluedot}}
\rput(-.45,2){\rnode{01a}{\smallunramifiedtree}}
\rput(-.15,2){\rnode{02a}{\smallunramifiedtree}}
\rput(.15,2){\rnode{03a}{\smallunramifiedtree}}
\rput(.45,2){\rnode{04a}{\smallunramifiedtree}}
\ncline {0a}{01a}
\ncline {0a}{02a}
\ncline {0a}{03a}
\ncline {0a}{04a}
}

\newcommand\tinyramifiedtree{
\rput(0,0){\rnode{0black}{\blackdot}}
\rput(-.10,1){\rnode{01red}{\reddot}}
\rput(-.05,1){\rnode{02red}{\reddot}}
\rput(.05,1){\rnode{03red}{\reddot}}
\rput(.10,1){\rnode{04red}{\reddot}}
\ncline {0black}{01red}
\ncline {0black}{02red}
\ncline {0black}{03red}
\ncline {0black}{04red}
}

\newcommand\redblackred{
\rput(0,0){\rnode{0red}{\reddot}}
\rput(-.20,1){\rnode{01black}{\tinyramifiedtree}}
\rput(.20,1){\rnode{02black}{\tinyramifiedtree}}
\ncline {0red}{01black}
\ncline {0red}{02black}
}

\newcommand\slantBTtreeWithGreenSprouts{
\rput(0,0){\rnode{0}{\unramifiedtree}}
\rput(2,0){\rnode{01p}{\blackdot}}
\rput(1.73,1){\rnode{02p}{\blackdot}}
\rput(-2,0){\rnode{03p}{\blackdot}}
\rput(-1.73,-1){\rnode{04p}{\blackdot}}

\rput(4,0){\rnode{01b}{\bluedot}}
\rput(5,0){\rnode{011p}{\blackdot}}
\rput(4.866,.5){\rnode{012p}{\blackdot}}
\rput(3.134,-.5){\rnode{013p}{\blackdot}}

\rput(3.46,2){\rnode{02b}{\bluedot}}
\rput(4.33,2.5){\rnode{021p}{\blackdot}}
\rput(2.46,2){\rnode{022p}{\blackdot}}
\rput(4.46,2){\rnode{023p}{\blackdot}}

\rput(-4,0){\rnode{03b}{\bluedot}}
\rput(-5,0){\rnode{031p}{\blackdot}}
\rput(-4.866,-.5){\rnode{032p}{\blackdot}}
\rput(-3.134,.5){\rnode{033p}{\blackdot}}

\rput(-3.46,-2){\rnode{04b}{\bluedot}}
\rput(-4.33,-2.5){\rnode{041p}{\blackdot}}
\rput(-2.46,-2){\rnode{042p}{\blackdot}}
\rput(-4.46,-2){\rnode{043p}{\blackdot}}

\ncline {0}{01p}
\ncline {0}{02p}
\ncline {0}{03p}
\ncline {0}{04p}

\ncline {01p}{01b}
\ncline {01b}{011p}
\ncline {01b}{012p}
\ncline {01b}{013p}

\ncline {02p}{02b}
\ncline {02b}{021p}
\ncline {02b}{022p}
\ncline {02b}{023p}

\ncline {03p}{03b}
\ncline {03b}{031p}
\ncline {03b}{032p}
\ncline {03b}{033p}

\ncline {04p}{04b}
\ncline {04b}{041p}
\ncline {04b}{042p}
\ncline {04b}{043p}

}

\newcommand{\BTtreewithsprouts}
{
\pspicture(-2,-3)(3,3)
\psset{unit=.5, nodesep=2.5pt}
\slantBTtreeWithGreenSprouts
\rput(2.46,2){\rnode{blacknode}{\blackdot}}
\rput(1.96,3){\rnode{redbranch1}{\redblackred}}
\rput(2.96,3){\rnode{redbranch2}{\redblackred}}
\rput(0,4){$L_0$}
\rput(1.86,6){$L_1$}
\rput(3.16,6){$L_2$}
\ncline{blacknode}{redbranch1}
\ncline{blacknode}{redbranch2}
\endpspicture
}

\newcommand\newgreekcross{
\rput(0,0){\rnode{0}{\broadunramifiedtree}}
\rput(0,2){\rnode{01p}{\blackdot}}
\rput(2,0){\rnode{02p}{\blackdot}}
\rput(-2,0){\rnode{03p}{\blackdot}}
\rput(0,-2){\rnode{04p}{\blackdot}}
\ncline {0}{01p}
\ncline {0}{02p}
\ncline {0}{03p}
\ncline {0}{04p}
}

\newcommand\newromancrosssmall{
\rput(0,1){\rnode{0}{\bluedot}}
\rput(-.5,1){\rnode{01p}{\blackdot}}
\rput(0,1.5){\rnode{02p}{\blackdot}}
\rput(.5,1){\rnode{03p}{\blackdot}}
\rput(0,0){\rnode{04p}{\blackdot}}
\rput(.35,1.35){\rnode{05g}{\greendot}}
\rput(-.35,1.35){\rnode{06g}{\greendot}}
\rput(-.35,.65){\rnode{07g}{\greendot}}
\rput(.35,.65){\rnode{08g}{\greendot}}
\rput(.35,.35){\rnode{1r}{\reddot}}
\rput(.35,-.35){\rnode{2r}{\reddot}}
\rput(-.35,.35){\rnode{3r}{\reddot}}
\rput(-.35,-.35){\rnode{4r}{\reddot}}

\ncline  {0}{01p}
\ncline  {0}{02p}
\ncline  {0}{03p}
\ncline  {0}{04p}
\ncline  {0}{05g}
\ncline  {0}{06g}
\ncline  {0}{07g}
\ncline  {0}{08g}
\ncline  {04p}{1r}
\ncline  {04p}{2r}
\ncline  {04p}{3r}
\ncline  {04p}{4r}
}

\newcommand\smallbroadramifiedtree{
\rput(0,0){\rnode{0srt}{\reddot}}
\rput(0,.25){\rnode{01srt}{\blackdot}}
\rput(.22,-.125){\rnode{02srt}{\blackdot}}
\rput(-.22,-.125){\rnode{03srt}{\blackdot}}
\ncline  {0srt}{01srt}
\ncline  {0srt}{02srt}
\ncline  {0srt}{03srt}
}

\newcommand\broadramifiedtree{
\rput(0,0){\rnode{0s}{\blackdot}}
\rput{75}(.50,.50){\rnode{01s}{\smallbroadramifiedtree}}
\rput{-15}(.50,-.50){\rnode{02s}{\smallbroadramifiedtree}}
\rput{165}(-.50,.50){\rnode{03s}{\smallbroadramifiedtree}}
\rput{-105}(-.50,-.50){\rnode{04s}{\smallbroadramifiedtree}}
\ncline  {0s}{01s}
\ncline  {0s}{02s}
\ncline  {0s}{03s}
\ncline  {0s}{04s}
}

\newcommand\smallbroadunramifiedtree{
\rput(0,0){\rnode{0ut}{\greendot}}
\rput(1.00000000000000, 0.000000000000000 ){\rnode{00ut}{\greendot}}
\rput(0.766271892468299, 0.642516448671201 ){\rnode{01ut}{\greendot}}
\rput(0.174345226373896, 0.984684590130583 ){\rnode{02ut}{\greendot}}
\rput(-0.49908019935562, 0.866555800056266 ){\rnode{03ut}{\greendot}}
\rput(-0.939207484081270, 0.343350115546407 ){\rnode{04ut}{\greendot}}
\rput(-0.940296393139069, -0.340356714418356 ){\rnode{05ut}{\greendot}}
\rput(-0.501837909222310, -0.864961682889699 ){\rnode{06ut}{\greendot}}
\rput(0.171207824314841, -0.985234936902553 ){\rnode{07ut}{\greendot}}
\rput(0.764221396308537, -0.644953996362709 ){\rnode{08ut}{\greendot}}
\ncline{0ut}{00ut}
\ncline{0ut}{01ut}
\ncline{0ut}{02ut}
\ncline{0ut}{03ut}
\ncline{0ut}{04ut}
\ncline{0ut}{05ut}
\ncline{0ut}{06ut}
\ncline{0ut}{07ut}
\ncline{0ut}{08ut}
}

\newcommand\broadunramifiedtree{
\rput(0,0){\rnode{0t}{\bluedot}}
\rput(2,2){\rnode{00t}{\smallbroadunramifiedtree}}
\rput{90}(-2,2){\rnode{01t}{\smallbroadunramifiedtree}}
\rput{180}(2,-2){\rnode{02t}{\smallbroadunramifiedtree}}
\rput{270}(-2,-2){\rnode{03t}{\smallbroadunramifiedtree}}
\ncline{0t}{00t}
\ncline{0t}{01t}
\ncline{0t}{02t}
\ncline{0t}{03t}
}

\newcommand\newromancross{
\rput(0,2){\rnode{A0}{\bluedot}}
\rput{90}(-1,2){\rnode{A01p}{\newromancrosssmall}}
\rput{0}(0,3){\rnode{A02p}{\newromancrosssmall}}
\rput{270}(1,2){\rnode{A03p}{\newromancrosssmall}}
\rput(0,0){\rnode{A04p}{\broadramifiedtree}}
\rput(.35,2.35){\rnode{05pg}{\greendot}}
\rput(-.35,2.35){\rnode{06pg}{\greendot}}
\rput(-.35,1.65){\rnode{07pg}{\greendot}}
\rput(.35,1.65){\rnode{08pg}{\greendot}}
\ncline {A0}{A01p}
\ncline {A0}{A02p}
\ncline {A0}{A03p}
\ncline {A0}{A04p}
\ncline{A0}{05pg}
\ncline{A0}{06pg}
\ncline{A0}{07pg}
\ncline{A0}{08pg}

}

\newcommand\fullpicture{
\rput(0,0){\newgreekcross}
\rput{0}(0,2){\newromancross}
\rput{90}(-2,0){\newromancross}
\rput{180}(0,-2){\newromancross}
\rput{270}(2,0){\newromancross}
}

\begin{document}
\title{Semistable models for modular curves of arbitrary level}
\author{Jared Weinstein}
\maketitle
\begin{abstract}  We produce an integral model for the modular curve $X(Np^m)$ over the ring of integers of a sufficiently ramified extension of $\Z_p$ whose special fiber is a {\em semistable curve} in the sense that its only singularities are normal crossings.  This is done by constructing a semistable covering (in the sense of Coleman) of the supersingular part of $X(Np^m)$, which is a union of copies of a Lubin-Tate curve.   In doing so we tie together non-abelian Lubin-Tate theory to the representation-theoretic point of view afforded by Bushnell-Kutzko types.

For our analysis it was essential to work with the Lubin-Tate curve not at level $p^m$ but rather at infinite level.  We show that the infinite-level Lubin-Tate space (in arbitrary dimension, over an arbitrary nonarchimedean local field) has the structure of a perfectoid space, which is in many ways simpler than the Lubin-Tate spaces of finite level.
\end{abstract}
\tableofcontents
\pagebreak

\section{Introduction:  The Lubin-Tate tower}
Let $K$ be a non-archimedean local field with uniformizer $\pi$ and residue field $k\isom\FF_q$, and let $n\geq 1$.  The {\em{Lubin-Tate tower}} is a projective system of affine formal schemes $\MM_m$ which parameterize deformations with level $\pi^m$ structure of a one-dimensional formal $\OO_K$-module of height $n$ over $\overline{\FF}_q$.  (For precise definitions, see \S\ref{ModuliOfFormalModules};  for a comprehensive historical overview of Lubin-Tate spaces, see the introduction to~\cite{Strauch}.)  After extending scalars to a separable closure of $K$, the Lubin-Tate tower admits an action of the triple product group $\GL_n(K)\times D^\times\times W_K$, where $D/K$ is the central division algebra of invariant $1/n$, and $W_K$ is the Weil group of $K$.  Significantly, the $\ell$-adic \'etale cohomology of the Lubin-Tate tower realizes both the Jacquet-Langlands correspondence (between $\GL_n(K)$ and $D^\times$) and the local Langlands correspondence (between $\GL_n(K)$ and $W_K$).  When $n=1$, this statement reduces to classical Lubin-Tate theory~\cite{LubinTate}.  For $n=2$ the result was proved by Deligne and Carayol (see \cite{Carayol:ladicreps}, \cite{Carayol:ladicreps2});  Carayol conjectured the general phenomenon under the name ``non-abelian Lubin-Tate theory".  Non-abelian Lubin-Tate theory was established for all $n$ by Boyer~\cite{Boyer} for $K$ of positive characteristic and by Harris and Taylor~\cite{HarrisTaylor} for $p$-adic $K$.  In both cases, the result is established by embedding $K$ into a global field and appealing to results from the theory of Shimura varieties or Drinfeld modular varieties.

In this paper we focus on the case that $n=2$ and $q$ is odd.  We construct a compatible family of {\em{semistable models}} $\hat{\MM}_m$ for each $\MM_m$ over the ring of integers of a sufficiently ramified extension of $K$.  For our purposes this means that the rigid generic fiber of $\hat{\MM}_m$ is the same as that of $\MM_m$, but that the special fiber of $\hat{\MM}_m$ is a locally finitely presented scheme of dimension 1 with only ordinary double points as singularities.   The weight spectral sequence would then allow for the computation of the cohomology of the Lubin-Tate tower of curves (along with the actions of the three relevant groups), and one could recover the result of Deligne-Carayol in a purely local manner, although we do not do this here.


The study of semistable models for modular curves begins with the Deligne-Rapoport model for $X_0(Np)$ in~\cite{DeligneRapoport}.  A semistable model for $X_0(Np^2)$ was constructed by Edixhoven in~\cite{Edixhoven}.   A stable model for $X(p)$ was constructed by Bouw and Wewers in~\cite{BouwWewers}.  A stable model for $X_0(Np^3)$ was constructed by Coleman and McMurdy in~\cite{ColemanMcMurdy}, using the notion of {\em{semistable coverings}} of a rigid-analytic curve by ``basic wide opens".  The special fiber of their model is a union of Igusa curves which are linked at each supersingular point of $X_0(N)\otimes\overline{\FF}_p$ by a certain configuration of projective curves, including in every case a number of copies of the curve with affine model $y^p-y=x^2$.  The same method was employed by Tsushima \cite{Tsushima} and Imai-Tsushima \cite{ImaiTsushima} for the curves $X_0(p^4)$ and $X_1(p^3)$, respectively;  the curve $y^p+y=x^{p+1}$ appears in the former.  In each of these cases the interesting part of the special fiber of the modular curve is the supersingular locus.  Inasmuch a Lubin-Tate curve (for $K=\Q_p$) appears as the rigid space attached to the $p$-adic completion of a modular curve at one of its mod $p$ supersingular points, the problem of finding a semistable model for a modular curve is essentially the same as finding one for the corresponding Lubin-Tate curve.  In this sense our result subsumes the foregoing results;  however our method cannot produce the ``intersection multiplicities" for the singular points of the special fiber.

We now summarize our main result.  Let $\hat{K}^{\nr}$ be the completion of the maximal unramified extension of $K$;  then $\MM_m$ is defined over $\Spf \OO_{\hat{K}^{\nr}}$.  

\begin{Theorem}
\label{mainthm}
Assume that $q$ is odd.  For each $m\geq 1$, there is a finite extension $L_m/\hat{K}^{\nr}$ for which $\MM_m$ admits a semistable model $\hat{\MM}_m$;  every connected component of the special fiber of $\hat{\MM}_m$ admits a purely inseparable morphism to one of the following smooth projective curves over $\overline{\FF}_q$:

\begin{enumerate}
\item The projective line $\mathbf{P}^1$,
\item The curve with affine model $xy^q-x^qy=1$,
\item The curve with affine model $y^q + y = x^{q+1}$,
\item The curve with affine model $y^q - y = x^2$.
\end{enumerate}
\end{Theorem}

\begin{rmk}
The mere existence of a semistable model of $\MM_m$ (after passing to a finite extension of the field of scalars) follows from the corresponding theorem about proper (algebraic) curves.   The formal scheme $\MM_m$ appears as the completion along a point in the special fiber of a proper curve over $\OO_K$ (e.g., the appropriate modular curve), and a semistable model of the proper curve restricts to a semistable model of $\MM_m$.  Furthermore, the theorem Drinfeld-Carayol allows one to predict in advance the field $L_m$ over which a semistable model appears.
The real content of the theorem is the assertion about the equations for the list of curves appearing therein.   A semistable model is unique up to blowing up, so the above theorem holds for all semistable models of $\MM_m$ if it holds for one of them.
\end{rmk}

\begin{rmk}  A purely inseparable morphism between nonsingular projective curves induces an equivalence on the level of \'etale sites and therefore an isomorphism on the level of  $\ell$-adic cohomology.
\end{rmk}

\begin{rmk} The equations for the curves appearing in Thm.~\ref{mainthm} were known some time ago by S. Wewers to appear in the stable reduction of $\MM_m$ (unpublished work).  Furthermore, it so happens that $xy^q-x^qy=1$ and $y^q+y=x^{q+1}$ determine isomorphic projective curves, but we have listed them separately because the nature of the group actions on these curves is different.
\end{rmk}

Let us explain some more features of our semistable models $\hat{\MM}_m$.  It is not the case that one can arrange for the semistable models $\hat{\MM}_m$ to be compatible:  there is no tower $\dots \to \hat{\MM}_2\to\hat{\MM}_1$ with finite transition maps.  Loosely speaking, the problem is that as $m\to \infty$, the singularities of the $\hat{\MM}_m$ accumulate around the {\em CM points}, that is, the points corresponding to deformations of $G_0$ with extra endomorphisms.

This problem can be remedied by removing the CM points entirely.  Let $\MM_{m,\overline{\eta}}^{\ad}$ be the geometric adic generic fiber\footnote{In much of this paper we work with adic spaces rather than rigid spaces, because presently we will be using adic spaces which do not come from rigid spaces.  There is a fully faithful functor (see \cite{HuberEtale}) from the category of rigid analytic varieties to the category of adic spaces, which identifies admissible opens with opens, and admissible open covers with open covers.  A separated adic space lies in the image of this functor if it is locally topologically of finite type.} of $\MM_{m}$.  The CM points constitute a closed subset of $\MM_{m,\overline{\eta}}^{\ad}$ whose topology is locally profinite.  Let $\MM_{m,\overline{\eta}}^{\ad,\nonCM}$ be the complement in $\MM_{m,\overline{\eta}}^{\ad}$ of the set of CM points.  This is an adic space (indeed we give a covering of it by affinoids).  Then $\MM_{m,\overline{\eta}}^{\ad,\nonCM}$ admits a semistable model which varies compatibly in $m$.  In fact much more is true.

\begin{Theorem}
\label{MainTheoremInDepth} The tower of adic spaces $\MM_{m,\overline{\eta}}^{\ad,\nonCM}$ admits a tower of semistable models $\hat{\MM}_m^{\nonCM}$ with finite transition maps.  Let $\hat{\MM}_{m,s}^{\nonCM}$ be the special fiber of $\hat{\MM}_m^{\nonCM}$.  For each $m$, let $C_m$ be an irreducible component of $\hat{\MM}_{m,s}^{\nonCM}$, such that the transition maps carry $C_{m+1}$ onto $C_m$.  Assume that $C_m$ has positive genus for some $m$.  Then for $m$ large enough, the morphism $C_{m+1}\to C_m$ is purely inseparable, and $\varprojlim C_m$ is the perfection of one of the curves listed in Thm. \ref{mainthm}.
\end{Theorem}


Thm. \ref{MainTheoremInDepth} allows us to associate a ``dual graph'' $\T$ to the tower of semistable models $\hat{\MM}_m^{\nonCM}$.  The vertices of $\T$ will correspond to towers $\cdots C_{m+1}\to C_m\to \cdots$ of irreducible components of $\hat{\MM}_{m,s}^{\nonCM}$;  two of these are adjacent when the corresponding irreducible components cross.  The graph $\T$ admits an action of the triple product group $\GL_2(K)\times D^\times \times W_K$, and the stabilizer in this group of a vertex of $\T$ acts on the corresponding scheme $\varprojlim C_m$.  

The geometric generic fiber of $\MM_m$ is highly disconnected, owing to the existence of the determinant morphism (see \S\ref{detmodules}).  In the limit, the set of connected components is a principal homogeneous space for $K^\times$.  The graph $\T$ has the same set of connected components.  One connected component of $\T$ is displayed in \S\ref{figures}, where it is called $\T^\circ$.

In theory one could draw a picture of the special fiber of any particular $\hat{\MM}_m$ by forming the quotient of the pictures described in \S\ref{figures} by the congruence subgroup $1+\pi^mM_2(\OO_K)$.  This would allow one to determine the structure of the reduction of a semistable model of the appropriate modular curve at level $m$, see \S\ref{figures}.

\subsection{Lubin-Tate space at infinite level}
\label{LTSpaceInfiniteLevel}
Before elaborating further, let us give some precise definitions.  Much of our notation has been borrowed from
\S2.1.1 of \cite{StrauchDeformationSpaces}.

Let $G_0$ be a one-dimensional formal $\OO_K$-module over $\overline{\FF}_q$ of height\footnote{Through the paper, the ``height" of a formal $\OO_K$-module will be understood to mean its height relative to $K$.  The same convention holds for quasi-isogenies between formal $\OO_K$-modules.} $n$.  Then $G_0$ is unique up to isomorphism.  Let $K_0=\hat{K}^{\nr}$ be the completion of the maximal unramified extension of $K$.  Let $\mathcal{C}$ be the category of complete local Noetherian $\OO_{K_0}$-algebras with residue field $\overline{\FF}_q$.  Let $j\in\Z$.  We consider the functor $\MM_{G_0,0}^{(j)}$ which associates to each $R\in\mathcal{C}$ the set of pairs $(G,\iota)$, where $G$ is a formal $\OO_K$-module over $R$ and $\iota\from G_0\to G\otimes_R\overline{\FF}_q$ is a quasi-isogeny of height $j$.
An isomorphism between pairs $(G,\iota)$ and $(G',\iota')$ is a quasi-isogeny of formal $\OO_K$-modules $f\from G\to G'$ which intertwines $\iota$ with $\iota'$.

Since $D=\End G_0\otimes_{\OO_K}K$ is a division algebra, a quasi-isogeny from $G_0$ to another formal $\OO_K$-module over $\overline{\FF}_q$ has height 0 if and only if it is an isomorphism.  Thus $\MM_{G_0,0}^{(0)}$ classifies formal $\OO_K$-modules $G$ together with an isomorphism $\iota\from G_0\to G\otimes_R\overline{\FF}_q$.  By~\cite{DrinfeldElliptic}, Prop. 4.2, $\MM_{G_0,0}^{(0)}=\Spf A_0$, where $A_0$ is a (noncanonically) isomorphic to the formal power series ring $\OO_{K_0}\powerseries{u_1,\dots,u_{n-1}}$.

One adds level structures to the moduli problem (see \ref{ModuliOfFormalModules}) to obtain formal schemes $\MM_{G_0,m}^{(j)}$, $m\geq 1$.  We put
\[\MM_{G_0,m}=\coprod_{j\in\Z} \MM_{G_0,m}^{(j)}. \]
We note that $\MM_{G_0,m}^{(j)}$ is isomorphic to $\MM_{G_0,m}^{(0)}$, though not canonically so.  In some of the paper we work with $\MM_{G_0,m}$ rather than $\MM_{G_0,m}^{(0)}$, so that we can take advantage of larger symmetry groups.  When $G_0$ is fixed in the discussion, we will drop it from the notation and simply write $\MM_m^{(j)}$ and $\MM_m$.

The formal schemes $\MM_m$ are rather mysterious.  Even at level zero, $\MM_0^{(0)}$ is the formal open unit ball of dimension $n-1$, but the action of the group $\OO_D^\times=\Aut G_0$ on $\MM_0^{(0)}$ is very difficult to write down.   It turns out however that the infinite-level deformation space
\[\MM_{\infty}=\varprojlim \MM_m\]
seems to be simpler than all of the spaces at finite level.  To prove Thm. \ref{mainthm} it was indispensable to work at infinite level, where a surprisingly nice description of $\MM_{\infty}$ emerges.  Results gathered about $\MM_{\infty}$ can then be translated into results about the individual $\MM_m$.

It will be helpful to first describe the case of $n=1$, so that $\MM_0^{(0)}$ parameterizes lifts of a Lubin-Tate formal $\OO_K$-module $G_0$ over $\overline{\FF}_q$.  One can find a coordinate $T$ on $G_0$ with respect to which $[\pi]_{G_0}(T)=T^q$.  $G_0$ lifts uniquely to $G/\OO_{\hat{K}^{\nr}}$, so that $\MM_0^{(0)}=\Spf\OO_{\hat{K}^{\nr}}$.  For each $m\geq 0$ we have $\MM_m^{(0)}=\Spf\OO_{K_m}$, where $K_m/\hat{K}^{\nr}$ is the field obtained by adjoining a $\pi^m$-torsion element $\lambda_m$ of $G$ (chosen compatibly).  Let $K_\infty=\cup_{m\geq 1} K_m$ and let $\hat{K}_\infty$ be its completion;  then $\hat{K}_\infty$ is the completion of the maximal abelian extension of $K$.  We have $\MM_{\infty}^{(0)}=\Spf \OO_{\hat{K}_\infty}$.  One finds in $\hat{K}_\infty$ an element $t=\lim_{m\to\infty} \lambda_m^{q^m}$ which admits arbitrary $q$th power roots.
If $K$ has positive characteristic, then in fact $\OO_{\hat{K}_\infty}\isom\overline{\FF}_q\powerseries{t^{1/q^\infty}}$
is a ring of {\em fractional power series} in $t$.
If $K$ has characteristic 0, then we can form the inverse limit $\OO_{K_\infty^{\flat}}=\varprojlim\OO_{K_\infty}/\pi$ along the Frobenius map, and then once again $\OO_{K_\infty^{\flat}}\isom\overline{\FF}_q\powerseries{t^{1/q^\infty}}$.  The fraction field $K_\infty^\flat$ of $\OO_{K_\infty^{\flat}}$ is the field of norms of $K_\infty$, as in \cite{FontaineWintenberger}.  In either case the field $\hat{K}_\infty$ is an example of a {\em perfectoid field}; see \cite{ScholzePerfectoidSpaces}.  See \S\ref{height1} for proofs of these claims.

Now return to the case of general $n$.  Let $A_m$ be the coordinate ring of $\MM_m^{(0)}$, so that $\MM_m^{(0)}=\Spf A_m$.  Each $A_m$ is complete with respect to the topology induced by the maximal ideal $I$ of $A_0$.  Let $A$ be the completion of $\varinjlim A_m$ with respect to the $I$-adic topology, so that $\MM_{\infty}^{(0)}=\Spf A$.  We show in Cor. \ref{Ainf} that if $K$ has positive characteristic, then
\[ A\isom \overline{\FF}_q\powerseries{X_1^{1/q^\infty},\dots,X_n^{1/q^\infty}}\]
is a ring of fractional power series in $n$ variables.  This is defined as the completion of $\overline{\FF}_q[X_1^{1/q^\infty},\dots,X_n^{1/q^\infty}]$ with respect to the ideal $(X_1,\dots,X_n)$.   If $K$ has characteristic 0, then $A$ contains topologically nilpotent elements $X_1,\dots,X_n$ admitting arbitrary $q$th power roots in $A$.  We define $A^{\flat}=\varprojlim A/\pi$, where the limit is taken with respect to the $q$th power Frobenius map.  Then
\[ A^\flat \isom \overline{\FF}_q\powerseries{X_1^{1/q^\infty},\dots,X_n^{1/q^\infty}}. \]
See Cor. \ref{Ainf}.

In either case, the parameters $X_1,\dots,X_n$ arise from Drinfeld's parameters on $A_m$ through a limiting process.  Furthermore, the action of the group $\GL_n(\OO_K)\times \OO_D^\times$ on these parameters can be determined directly from the formal $\OO_K$-module $G_0$ itself.

There is a continuous homomorphism morphism $\OO_{\hat{K}_\infty}\to A$, which sends $t$ to a rather complicated fractional power series in $X_1,\dots,X_n$.  This power series can be interpreted as the determinant of a {\em formal vector space}, see \S\ref{detFVS}.  Let $\MM^{(0),\ad}_{\infty,\overline{\eta}}$ be the geometric adic generic fiber of $\MM_{\infty}^{(0)}$, where $\overline{\eta}=\Spa(\C,\OO_{\C})$ and $\C/K$ is a complete algebraically closed field.  That is, $\MM^{(0),\ad}_{\infty,\overline{\eta}}$ is the set of continuous valuations $\abs{\;}$ (in the sense of Huber, \cite{Hub94}) on $A\hat{\otimes}\OO_{\C}$ for which $\abs{\pi}\neq 0$.  The above descriptions of $A$ show that $\MM^{\ad}_{\infty,\overline{\eta}}$ is a {\em perfectoid space}, see \cite{ScholzePerfectoidSpaces}.   In light of the above description of $A$ it would appear that $\MM^{\ad}_{\infty,\overline{\eta}}$ is a very simple sort of space, let alone that it encodes the Langlands correspondence!  In fact it is the complexity of the element $t\in A$ which accounts for the interesting cohomological behavior of $\MM^{\ad}_{\infty,\overline{\eta}}$.

Much of the above was probably known to the experts, although perhaps not in this precise form.  In \cite{FarguesGenestierLafforgue}, an isomorphism between the Lubin-Tate and Drinfeld towers is constructed.  For this it is necessary to work with the infinite-level versions of both towers.   Roughly speaking, the authors work with an integral model not of the whole Lubin-Tate space (as we have done), but rather with an integral model of a ``fundamental domain", whose coordinate ring carries the structure of a perfectoid affinoid algebra.  Certainly the important role of the determinant is recognized in \cite{FarguesGenestierLafforgue}.

All Rapoport-Zink spaces at infinite level are perfectoid spaces which can be described in terms of $p$-adic Hodge theory, \cite{ScholzeWeinstein}.   There we prove a general duality theorem relating basic Rapoport-Zink spaces to one another, and in particular we arrive at an isomorphism between the infinite-level Lubin-Tate and Drinfeld spaces which does not require any integral models at all.

\subsection{Outline of the paper}
In \S\ref{LTperfectoid}, we review the construction of the Lubin-Tate tower attached to a one-dimensional formal $\OO_K$-module $G_0$.  We consider $\tilde{G}=\varprojlim G$ (limit along multiplication by $\pi$), where $G$ is any lift of $G_0$.  Then $\tilde{G}$ does not depend on the choice of lift.  $\tilde{G}$ carries the structure of a $K$-vector space object in the category of formal schemes.  Following \cite{FarguesFontaine}, we call $\tilde{G}$ a {\em formal vector space}.  A choice of coordinate on $G_0$ allows us to identify $\tilde{G}$ with $\Spf\OO_{\hat{K}^{\nr}}\powerseries{T^{1/q^\infty}}$.

We also review relevant results on determinants of formal $\OO_K$-modules, as these play an important role.  Let $\wedge^n G_0$ be the top exterior power of $G_0$.  Then $\wedge^n G_0$ is a formal $\OO_K$-module of height one and dimension one;  {\em i.e.} it is a Lubin-Tate formal $\OO_K$-module.  If $\MM_{G_0,\infty}$ is the Lubin-Tate deformation space of $G_0$ at infinite level, and similarly for $\MM_{\wedge^n G_0,\infty}$, then there is a determinant morphism $\MM_{G_0,\infty}\to\MM_{\wedge^n G_0,\infty}$.  The main result of the section is that there is a Cartesian diagram
\[
\xymatrix{
\mathcal{M}_{G_0,\infty} \ar[r] \ar[d] & \mathcal{M}_{\wedge^n G_0,\infty} \ar[d] \\
\tilde{G}^n \ar[r] & \widetilde{\wedge^n G}
}
\]
where the horizontal arrows are determinant morphisms (Thm. \ref{Cartesian}).   Passing to the geometric generic fiber, we arrive at the infinite-level Lubin-Tate space $\MM_{G_0,\overline{\eta}}^{\ad}$, which is a perfectoid space.

In \S\ref{RepresentationTheoreticPreparations} we review Carayol's description of the cohomology of the Lubin-Tate tower for $\GL_2(K)$ (non-abelian Lubin-Tate theory), see Thm. \ref{JLinCohomology}.  Informally, the cohomology splits up into a sum of representations of $\GL_2(K)\times D^\times\times W_K$ of the form $\Pi\otimes\JL(\check{\Pi})\otimes\mathcal{H}(\Pi)'$, where $\Pi$ runs over discrete series representations, $\JL$ is the Jacquet-Langlands correspondence and $\Pi\mapsto\mathcal{H}(\Pi)'$ is a normalized local Langlands correspondence.  We also review the theory of Bushnell-Kutzko types for $\GL_2(K)$ and its inner form $D^\times$, as these play a major role in our work. This theory furnishes a classification of supercuspidal representations of these groups according to which ``strata" they contain.  A stratum is essentially a one-dimensional character of a compact subgroup of $\GL_2(K)$ (or $D^\times$) which belongs to a certain explicit class.

In \S\ref{CMPointsAndLinkingOrders} we specialize to the case that $G_0$ has height 2.  We consider the set of CM points in $\MM_{G_0,\infty,\overline{\eta}}^{\ad}$, where $G_0$ has height 2.  Points in $\MM_{G_0,\infty,\overline{\eta}}^{\ad}$ with CM by a quadratic extension $L/K$ correspond to pairs of embeddings $L\injects M_2(K)$ and $L\injects D$.  If $\x$ is a point with CM by $L$, then its stabilizer in $\GL_2(K)\times D^\times$ is the diagonally embedded $L^\times$.  To each such $\x$ and each integer $m\geq 0$ we associate the following data:
\begin{enumerate}
\item An $\OO_L$-order $\mathcal{L}_{\x,m}\subset M_2(K)\times D$, which is normalized by the diagonally embedded $L^\times$,
\item The group $\K_{\x,m}=L^\times\L_{\x,m}^\times$ and its subgroup $\K_{\x,m}^1$ consisting of pairs $(g_1,g_2)\in\K_{\x,m}$ with $\det(g_1)=N(g_2)$, where $N\from D\to K$ is the reduced norm,
\item A smooth affine curve $C_{\x,m}/\overline{k}$ equipped with an action of $\K_{\x,m}^1$ (only given outside the case that $L/K$ is ramified and $m$ is even).
\end{enumerate}

The orders $\mathcal{L}_{\x,m}$, which we call ``linking orders", were first constructed in \cite{WeinsteinFourier}.  Their study links together the Bushnell-Kutzko type theories for $\GL_2(K)$ and $D^\times$.  In Thm. \ref{CurvesRealizeJLC} we show that if $\Pi$ is a supercuspidal representation of $\GL_2(K)$ with coefficients in $\overline{\Q}_\ell$ ($\ell\nmid q$), then there exists a pair $(\x,m)$ (depending on the strata contained in $\Pi$) such that $\Pi\otimes\JL(\check{\Pi})$ is contained in the representation of $\GL_2(K)\times D^\times$ induced from the representation of $\K_{\x,m}^1$ on $H^1_c(C_{x,m},\overline{\Q}_\ell)$.

In \S\ref{specialaffinoids} we identify a family of special affinoid subspaces $\set{\ZZ_{\x,m}}$ of the $\MM_{\infty,\overline{\eta}}^{\ad}$, parameterized by pairs $(\x,m)$.  These have the following properties (Thm. \ref{ExistenceOfAffinoid}):
\begin{enumerate}
\item $\ZZ_{\x,m}$ is stabilized by $\K_{\x,m}^1$, and
\item There exists a nonconstant $\K_{\x,m}^1$-equivariant map $\overline{\ZZ}_{\x,m}\to C_{\x,m}$, where $\overline{\ZZ}_{\x,m}$ is the reduction of the affinoid $\ZZ_{\x,m}$ (once again the case of $L/K$ ramified and $m$ even is excluded).
\end{enumerate}
The proof of Thm. \ref{ExistenceOfAffinoid} is a long, delicate calculation, and we regret that we could not find a coordinate-free method.  The payoff of Thm. \ref{ExistenceOfAffinoid} is the observation that the special affinoids $\ZZ_{\x,m}$ exhaust the entire supercuspidal part of the cohomology of the Lubin-Tate tower.

In \S\ref{semistablecoverings}, we translate the results of \S\ref{specialaffinoids} back to finite level.  We construct a $\GL_2(K)\times D^\times$-equivariant graph $\T$ whose vertices are equivalence classes of pairs $(\x,m)$, and also a finite-level version $\T^{(m)}$.  For every vertex $v$ of $\T^{(m)}$ we get an open affinoid $\ZZ_v^{(m)}\subset \MM_{m,\overline{\eta}}^{\ad}$, equal to the image of the corresponding $\ZZ_{x,m}$.   The cohomology of the $\overline{\ZZ}_v^{(m)}$ exhausts all of $H^1_c(\MM_{m,\overline{\eta}},\overline{\Q}_\ell)$, except for the part coming from the boundary.  At this point our argument starts to resemble the method employed in \cite{ColemanMcMurdy}.  We find a covering of $\MM_{m,\overline{\eta}}^{\ad}$ by ``wide opens", whose underlying affinoids are the $\ZZ_{v}^{(m)}$.  For cohomological reasons this must be a semistable covering.  This means that pairs of wide opens intersect in annuli, and that the $\ZZ_{v}^{(m)}$ have smooth reduction.  By the general theory of \cite{Coleman:StableMapsOfCurves}, a semistable covering of $\MM_{m,\eta}^{\ad}$ corresponds to a semistable model $\hat{\MM}_m$.  The dual graph of this model is $\T^{(m)}$.  Finally, we complete the proof of Thm. \ref{MainTheoremInDepth}.

The tree $\T$ (or rather one of its connected component $\T^\circ$) is depicted in \S\ref{figures}.  We note that $\T^\circ$ contains a copy of the Bruhat-Tits tree for $\text{PGL}_2(K)$:  this reflects the structure of fundamental domains in $\MM_{\infty,\eta}^{\ad}$ already observed in \cite{FarguesGenestierLafforgue}.  The ends (infinite paths) of the Bruhat-Tits tree are in correspondence with $\mathbf{P}^1(K)$.  On the other hand $\T$ has additional ends which are in correspondence with the set of CM points.  We sketch a procedure for computing the special fiber of $\hat{\MM}_m$.

\subsection{Acknowledgments}
This project originated in conversations with Robert Coleman while the author was a graduate student at Berkeley.  We thank Peter Scholze for his work on perfectoid spaces, thus introducing a category which contains our deformation spaces at infinite level.   We thank Mitya Boyarchenko, Jay Pottharst and Peter Scholze for very helpful conversations.  Parts of this project were completed at the Institute for Advanced Study, and supported by NSF grant DMS-0803089.

\section{The Lubin-Tate deformation space at infinite level}
\label{LTperfectoid}

\subsection{Moduli of one-dimensional formal modules}
\label{ModuliOfFormalModules}

Let $K_0=\hat{K}^{\nr}$.  Recall that the functor  $\M_{G_0,0}$ is represented by the formal scheme $\Spf A_0$, where $A_0$ is (noncanonically) isomorphic to a formal power series ring
$\OO_{K_0}\powerseries{u_1,\dots,u_{n-1}}$ in one variable.  Thus there is a universal formal $\OO_K$-module $G^{\univ}$ defined over $A_0$.   We follow the construction of $G^{\univ}$ in \cite{GrossHopkins}, \S5 and \S12.  Over the polynomial ring $\OO_{K_0}[v_1,v_2,\dots]$ we can consider the universal $p$-typical formal $\OO_K$-module $F$, whose logarithm $f(T)=\log_F(T)$ satisfies Hazewinkel's ``functional equation"
\[ f(T) = T + \sum_{i\geq 1}\frac{v_i}{\pi}f^{q^i}(X^{q^i}). \]
Here $f^{q^i}$ is the series obtained from $f(X)$ by replacing each variable $v_j$ by $v_j^{q^i}$.   Then multiplication by $\pi$ in $F$ satisfies the congruences
\begin{equation}
\label{picongs}
 [\pi]_{F}(T)\equiv v_k T^{q^k}\pmod{\pi,v_1,\dots,v_{k-1},T^{q^k+1}},
\end{equation}
as in \cite{GrossHopkins}, Prop. 5.7.  Then $G^{\univ}$ is the push-forward of $F$ through the homomorphism $\OO_{K_0}[v_1,v_2,\dots]\to A_0$ which sends
\[ v_i\mapsto \begin{cases}
u_i,& i=1,\dots,n-1 \\
1,& i=n\\
0,& i>n.
\end{cases}\]
Let $[\pi]_{G^{\univ}}(T)=c_1 T+c_2T^2+\dots$. (Thus $c_1=\pi$.) It follows from Eq. \eqref{picongs} that in $A_0=\OO_{K_0}\powerseries{u_1,\dots,u_{n-1}}$ we have the congruences
\[
\begin{matrix}
c_q&\equiv& u_1 &\pmod{\pi}\\
c_{q^2}&\equiv& u_2&\pmod{\pi,u_1}\\
&\vdots&&\\
c_{q^{n-1}}&\equiv& u_{n-1}&\pmod{\pi,u_1,\dots,u_{n-2}}\\
c_{q^n}&\equiv& 1&\pmod{\pi,u_1,\dots,u_{n-1}}.
\end{matrix}
\]
These congruences have the following immediate consequence.
\begin{lemma}
\label{generate}
The coefficients $c_1,c_q,\dots,c_{q^{n-1}}$ of $[\pi]_{G^{\univ}}(T)$ generate the maximal ideal of $A_0$, and $c_{q^n}\in A_0$ is a unit.
\end{lemma}

\subsection{Level structures}
\label{levelstructures}

For an algebra $R\in\mathcal{C}$ and a deformation $(G,\iota)\in \M_{G_0,0}(R)$, a {\em Drinfeld level $\pi^m$ structure} on $G$ is an $\OO_K$-module homomorphism
\[\phi\from(\pi^{-m}\OO_K/\OO_K)^{\oplus n}\to G(R)\] for which the relation
\[\prod_{x\in (\gp^{-1}_K/\OO_K)^{\oplus n}}\left(X-\phi(x)\right)\biggm\vert [\pi]_{G}(X)\] holds in $R\powerseries{X}$.  The images under $\phi$ of the standard basis elements $(\pi^{-n},0)$ and $(0,\pi^{-n})$ of $(\pi^{-n}/\OO_K)^{\oplus n}$ form a {\em{Drinfeld basis}} of $G[\pi^m]$ over $R$.

\begin{rmk}\label{drinfeld1}Note that $x_1,\dots,x_n$ is a Drinfeld basis of $G[\pi^m](R)$ if and only if $\pi^{m-1}x_1,\dots,\pi^{m-1}x_n$ is a Drinfeld basis of $G[\pi](R)$.
\end{rmk}

Let $\MM_{G_0,m}^{(j)}$ denote the functor which assigns to each $R\in\mathcal{C}$ the set of deformations $(G,\iota)\in\MM_{G_0,0}^{(j)}(R)$ together with a Drinfeld level $\pi^m$ structure on $G$ over $R$.  Let $\MM_{G_0,m}$ be the union of the $\MM_{G_0,m}^{(j)}$, $j\in\Z$.  

By~\cite{DrinfeldElliptic}, Prop. 4.3, the functor $\MM_{G_0,m}^{(0)}$ is represented by a formal scheme $\Spf A_m$, where $A_m$ is finite, flat, and generically \'etale over $A_0\isom\OO_{K_0}\powerseries{u_1,\dots,u_{n-1}}$.  The universal Drinfeld level structure on $A_m$ corresponds to $m$ topologically nilpotent elements $X_1^{(m)},\dots,X_n^{(m)}\in A_m$.  Drinfeld shows that $A_m$ is a regular local ring with parameters $X_1^{(m)},\dots,X_n^{(m)}$.

\subsection{The case of height one}\label{height1}  In this paragraph we assume $n=1$. Then $G_0$ is a Lubin-Tate formal $\OO_K$-module over $\overline{\FF}_q$, which admits a unique deformation $G$ to $\OO_{K_0}$.  In fact after choosing a suitable coordinate on $G$, we may assume $[\pi]_{G}(T)=\pi T + T^q$.  For each $m\geq 1$, write $\Phi_m(T)=[\pi^m]_G(T)/[\pi^{m-1}]_G(T)$.  Then $\Phi_m(T)$ is an Eisenstein polynomial of degree $q^{m-1}(q-1)$ and a unit in $\OO_{K_0}\powerseries{T}$.

\begin{lemma} Let $R\in\mathcal{C}$.  An element $x\in G[\pi^m](R)$ constitutes a Drinfeld basis if and only if it is a root of $\Phi_m(T)$.
\end{lemma}


\begin{proof} \label{M1} By Remark \ref{drinfeld1}, and because $\Phi_1([\pi^{m-1}]_G(T))=\Phi_m(T)$, we may assume $m=1$.  The condition for $x$ to be is a Drinfeld basis of $G[\pi](R)$ is the condition that $T\prod_{a\in k^\times}(T-[a]_G(x))$ is divisible by $[\pi]_G(T)=T\Phi_1(T)$.  This is equivalent to the condition that $x$ is a root of $\Phi_1(T)$.
\end{proof}

Let $K_m$ be the field obtained by adjoining the $\pi^m$-torsion in $G$ to $K_0$.  Lemma \ref{M1} implies that  $\MM_{G_0,m}^{(0)}=\Spf\OO_{K_m}$.   Note that by local class field theory, the union $K_\infty=\bigcup_m K_m$ is the compositum of $K_0$ with the maximal abelian extension of $K$.  The following fact will be useful later.

\begin{lemma}
\label{OKinfty} The $q$th power Frobenius map is surjective on $\OO_{K_{\infty}}/\pi$.
\end{lemma}

\begin{proof}  Let $\lambda_1,\lambda_2,\dots\in \OO_{K_\infty}$ be a compatible sequence of roots of $[\pi^m]_G(T)$, $m\geq 1$.  Then $\lambda_m$ generates $\OO_{K_m}$ over $\OO_{K_0}$.
Since $\lambda_{m}=[\pi]_G(\lambda_{m+1})$, and $[\pi]_G(T)\pmod \pi$ is a power series in $T^q$, we have that every element of $\OO_{K_m}/\pi$ is a $q$th power in $\OO_{K_{m+1}}/\pi$.  The result follows.
\end{proof}

Let $\hat{K}_\infty$ be the $\pi$-adic completion of $K_\infty$.

\begin{prop}
\label{OKinf} If $K$ has positive characteristic, then $\OO_{\hat{K}_{\infty}}\isom \overline{\FF}_q\powerseries{t^{1/q^\infty}}$, where $\overline{\FF}_q\powerseries{t^{1/q^\infty}}$ is the $t$-adic completion of $\overline{\FF}_q[t^{1/q^\infty}]$.  If $K$ has characteristic 0, let $\OO_{K_\infty^{\flat}}=\varprojlim \OO_{K_\infty}/\pi$, where the inverse limit is taken with respect to the $q$th power Frobenius map.  Then $\OO_{K_\infty^{\flat}}\isom\overline{\FF}_q\powerseries{t^{1/q^\infty}}$.
\end{prop}

\begin{proof}  The element $\varpi=\lim_{m\to\infty} \lambda_m^{q^{m-1}}$ belongs to $\hat{K}_\infty$ and has an obvious system of $q$th power roots, which we write as $\varpi^{1/q^m}$, $m\geq 1$.
We have the congruences $\lambda_m\equiv \lambda_{m-1}^{q}\equiv \lambda_{m_2}^{q^2}\equiv \dots$ modulo $\pi\OO_{K_\infty}$, which shows that $\lambda_m\equiv \varpi^{1/q^{m-1}}\pmod{\pi\OO_{\hat{K}_\infty}}$, and therefore (since $\OO_{K_\infty}$ is generated by the $\lambda_m$) there is a surjection $\overline{\FF}_q[t^{1/q^\infty}]\to\OO_{K_\infty}/\pi$ which sends $t$ to $\varpi$.

Assume $K$ has positive characteristic.  Then there is a continuous $\overline{\FF}_q$-algebra homomorphism $\phi\from \overline{\FF}_q\powerseries{t^{1/q^\infty}}\to \OO_{\hat{K}_{\infty}}$ which sends $t$ to $\varpi$ and which is a surjection modulo $\pi$.  In particular it is a surjection modulo $\varpi$, because in $\hat{K}_\infty$ we have $\abs{\varpi}=\abs{\pi}^{1/(q-1)}$.   Thus any $b\in \OO_{\hat{K}_\infty}$ can be written $b=\phi(a_1)+\varpi b_1=\phi(a_1)+\phi(t)b_1$, with $a_1\in\overline{\FF}_q[t^{1/q^\infty}]$ and $b_1\in \OO_{\hat{K}_\infty}$.  But then we can write $b_1=\phi(a_2)+\phi(t)b_2$, and so forth, the result being that $b=\phi(a_1+ta_2+\dots)$.  Thus $\phi$ is surjective.  The injectivity of $\phi$ follows from the fact that any nonzero element of $\overline{\FF}_q\powerseries{t^{1/q^\infty}}$ equals $t^\alpha$ times a unit for some $\alpha\in\Z[1/q]$, so that if $\phi$ has a nonzero kernel, we would have $\phi(t^\alpha)=\varpi^\alpha=0$ for some $\alpha$, which is absurd.  Thus $\phi$ is an isomorphism.

Now assume $K$ has characteristic 0.  We have put $\OO_{K^{\flat}_\infty}=\varprojlim \OO_{K_\infty}/\pi$;  this makes $\OO_{K^{\flat}}$ the ring of integers in a complete nonarchimedean valuation field $K^\flat$ containing $\overline{\FF}_q$.  We have a continuous $\overline{\FF}_q$-algebra homomorphism $\phi\from \overline{\FF}_q\powerseries{t^{1/q^\infty}}\to\OO_{K^{\flat}}$ which sends $t$ to the sequence $\varpi^\flat=(\varpi,\varpi^{1/q},\dots)$.  We have an isomorphism $\OO_{K^{\flat}_\infty}/\varpi^{\flat}\to\OO_{\hat{K}_\infty}/\varpi$ given by projection onto the first coordinate.  We see that $\phi$ is once again surjective modulo $\varpi$.  The argument now continues as in the previous paragraph.
\end{proof}

\subsection{Formal vector spaces}

Suppose $A$ is a topological ring which is separated and complete for the topology induced by an ideal of definition $I$.  For such a ring we write $\Nil(A)$ for the set of topologically nilpotent elements of $A$, which is to say that $\Nil(A)$ is the radical of $I$.  Of course we allow for the trivial case in which $I=0$ and $A$ is discrete, in which case $\Nil(A)$ is the set of nilpotent elements of $A$.  Let $\Alg_A$ be the category of topological $A$-algebras $R$ which are separated and complete for the topology induced by an ideal $J$ (which may be assumed to contain the image of $I$).  Also let $\Mod_{\OO_K}$ be the category of $\OO_K$-modules, and let $\Vect_K$ be the category of $K$-vector spaces.

Recall that $G_0$ is a formal $\OO_K$-module over $\overline{\FF}_q$ of dimension 1 and height $n$.  $G_0$ induces a functor $\Alg_{\overline{\FF}_q}\to \Mod_{\OO_K}$ whose value on an object $R$ is $\Nil(R)$ with the  $\OO_K$-module structure afforded by $G_0$.  This functor is representable by a formal scheme which we will simply call $G_0$.  A choice of coordinate on $G_0$ is equivalent to a choice of isomorphism $G_0\isom \Spf \overline{\FF}_q\powerseries{X}$.  Now consider the functor $\tilde{G}_0\from\Alg_{\overline{\FF}_q}\to \Vect_K$ defined by
\[ \tilde{G}_0(R)= \varprojlim G_0(R), \]
where the transition map is multiplication by a uniformizer $\pi$.  Let us call $\tilde{G}$ the {\em formal $K$-vector space} associated to $G$.

\begin{prop}
\label{prorep} $\tilde{G}_0$ is representable by an affine formal scheme isomorphic to $\overline{\FF}_q\powerseries{X^{1/q^\infty}}$.
\end{prop}

See Prop. 3.1.2 of \cite{WeinsteinFormalVectorSpaces} for a proof (in that context $K$ has positive characteristic, but it makes no difference).

Now let $A$ be an object of $\Alg_{\OO_{K_0}}$ with ideal of definition $I$.

\begin{prop} \label{crystalline} Let $G$ be a one-dimensional formal $\OO_K$-module over $A$, and define a functor $\tilde{G}\from \Alg_A\to\Vect_K$ by $\tilde{G}(R)=\varprojlim G(R)$ (inverse limit along multiplication by $\pi$).
\begin{enumerate}
\item The natural reduction map $\tilde{G}(R)\to \tilde{G}(R/I)$ is an isomorphism.
\item If $A/I$ is a perfect field, and if $G\otimes A/I$ has finite height, then $\tilde{G}$ is representable by $A\powerseries{X^{1/q^\infty}}$.
\end{enumerate}
\end{prop}

\begin{proof} Choose a coordinate on $G$, so that the $G(R)$ may be identified with $\Nil(R)$ of $R$.  Let $I_R$ be the extension of $I$ to $R$, so that $I_R$ is nilpotent.  If $(x_1,x_2,\dots)\in\tilde{G}(R)$ lies in the kernel of $\tilde{G}(R)\to \tilde{G}_0(R/I)$, then each $x_i$ lies in $I_R$.  But the power series giving multiplication by $\pi$ in $G$ has $\pi\in I$ as its linear terms, so it carries $I_R^m$ onto $I_R^{m+1}$.  It follows that each $x_i$ lies in $\cap_{m\geq 1} I_R^m =0$.

We show that $\tilde{G}(R)\to \tilde{G}_0(R/I)$ is surjective using the standard ``Teichm\"uller lift".
Suppose $(x_1,x_2,\dots)\in\tilde{G}_0(R/I)$.
Since $I$ is nilpotent in $R$, we may lift each $x_i$ to an element $y_i\in G(R)$.
Then the sequence $y_i$,$\pi y_{i+1}$, $\pi^2 y_{i+2},\dots$ must converge to an element $z_i\in G(R)$.
 Then $(z_1,z_2,\dots)\in\tilde{G}(R)$ lifts $(x_1,x_2,\dots)\in \tilde{G}(R/I)$.  This settles part (1).

For part (2), let $G_0=G\otimes_A A/I$.  By Lemma \ref{prorep}, the functor $\tilde{G}_0\from \Alg_{R/I} \to \Vect_K$ is representable by a formal scheme isomorphic to $\Spf (A/I)\powerseries{X^{1/q^\infty}}$.  Thus if $R$ is an $A/I$-algebra, then $\tilde{G}_0(R)$ may be identified with $\varprojlim \Nil(R)$ (limit taken with respect to $x\mapsto x^q$).   Now suppose $R$ is an object of $\Alg_A$;  then by part (1) we have
\begin{eqnarray*}
\tilde{G}(R) &\isom& \tilde{G}(R/I) \\
&=& \tilde{G}_0(R/I) \\
&\isom& \varprojlim_{x\mapsto x^q} \Nil(R/I) \\
&\isom& \varprojlim_{x\mapsto x^q} \Nil(R).
\end{eqnarray*}
In the last step, we have used the standard Teichm\"uller lift procedure.  This functor is representable by $\Spf A\powerseries{X^{1/q^\infty}}$.
\end{proof}

\begin{rmk}
\label{Giscrystalline}
The first part of the proposition shows that the functor $\tilde{G}$ only depends on $G_0=G\otimes_A A/I$, in a functorial sense.  That is, there is a functor
\begin{eqnarray*}
\set{\text{Formal $\OO_K$-modules over $A/I$}} &\to& \set{\text{Formal schemes over $\Spf A$}} \\
G_0 &\mapsto & \tilde{G},
\end{eqnarray*}
where $\tilde{G}$ represents the functor $R\mapsto \tilde{G}_0(R/\pi)$ for any object $R$ of $\Alg_A$.  Then if $G'$ is a lift of $G_0$ to a formal $\OO_K$-module over $R$, then we have a canonical isomorphism of functors $\tilde{G}'\isom\tilde{G}$.
\end{rmk}

\begin{rmk}
\label{boldfaceX} In the situation of the second part of the proposition, where $A/I$ is a perfect field and $G$ is a formal $\OO_K$-module over $A$, we will often use boldface letters, such as $\mathbf{X}$, to denote elements of the $K$-vector space $\tilde{G}(R)$, where $R$ is an object of $\Alg_A$.  Such an element corresponds to a compatible sequence $(\mathbf{X}^{(1)},\mathbf{X}^{(2)},\dots)$ in the inverse limit $\varprojlim G(R)$.  Assume that a coordinate on $G$ has been chosen, so that $G(R)$ may be identified as a set with $\Nil(R)$.  Then the proposition shows that $\mathbf{X}$ corresponds to a topologically nilpotent element $X\in R$ admitting arbitrary $q$th roots, which will simply be written $X^{1/q^\infty}$.

Let us record the relationship between $\mathbf{X}$ and $X$.  The formal module $G_0$ has height $n$, so $[\pi]_{G_0}(T)$ is a power series in $T^{q^n}$.  If $A/I$ is algebraically closed, we may even perform a change of coordinate so that $[\pi]_{G_0}(T)=T^{q^n}$.  Then
\[ X = \lim_{m\to\infty} \left(\X^{(m)}\right)^{q^{mn}}. \]
This pattern (boldface for elements of $\tilde{G}(R)$, Roman for elements of $R$) will be useful later on.
\end{rmk}

\subsection{Determinants of formal modules}
\label{detmodules}
Assume for the moment that $K$ has characteristic 0.  Let $\wedge^n G_0$ be the formal group whose (covariant) Dieudonn\'e module is the top exterior power of the Dieudonn\'e module of $G_0$.  Then $\wedge^n G_0$ has height one and dimension one;  {\em i.e.} it is the Lubin-Tate formal $\OO_K$-module over $\overline{\FF}_q$.  Therefore $\wedge^n G_0$ admits a unique deformation $\wedge^n G$ to any $R\in\mathcal{C}$.

Now let $(G,\iota)$ be a deformation of $G_0$ to $R\in\mathcal{C}$.
\begin{Theorem}
\label{lambda} For every $m\geq 1$ there exists an alternating and multilinear map of $\OO_K$-module schemes
\[ \lambda_m\from G[\pi^m]^n \to \wedge^n G[\pi^m] \]
of formal $\OO_K$-module schemes over $R$, which is universal in the sense that any alternating and multilinear map from $G[\pi^m]^n$ into another $\OO_K$-module scheme must factor through $\lambda$.
\end{Theorem}

\begin{proof} This is a special case of the main theorem of Hadi Hedayatzadeh's thesis, \cite{Hedayatzadeh}, Thm. 9.2.36.  There the author constructs arbitrary exterior powers of arbitrary $\pi$-divisible $\OO_K$-modules $G$ over arbitary locally Noetherian $\OO_K$-schemes, so long as $\dim G\leq 1$.  Hedayatzadeh shows that if $\dim G=1$ then $\wedge^kG$ has dimension $\tbo{n}{k}$ and height $\tbo{n-1}{k-1}$.  Thus in our case, the $n$th exterior power of $G$ has dimension one and height one;  {\em i.e.} it is isomorphic to the unique Lubin-Tate formal module, which we have called $\wedge^n G$.
\end{proof}

\begin{prop} \label{DrinfeldtoDrinfeld}  Let $R\in\mathcal{C}$, let $(G,\iota)$ be a deformation of $G_0$ to $R$, and let $x_1,\dots,x_n\in G[\pi^m]$ be a Drinfeld level $\pi^m$ level structure.  Then $\lambda_m(x_1,\dots,x_n)\in\wedge^nG[\pi^m]$ is a Drinfeld level $\pi^m$ structure.
\end{prop}

\begin{proof} It suffices to treat the universal case, where $R=A_m$, $G$ is the universal deformation, and $X_1,\dots,X_n\in G[\pi^m](A_m)$ is the universal level structure.
Let $X=\lambda_m(X_1,\dots,X_n)$.  By Remark \ref{drinfeld1}, we can reduce to the case that $m=1$.
It suffices to show that $X$ is a primitive element of $\wedge^n G[\pi](A_1)$.
Now we appeal to the fact that $A_1$ is a domain:  if $X$ isn't a primitive element, then it must be 0.
But this would mean that $\lambda_1=0$, which would contradict the fact that $\wedge^nG[\pi]\neq 0$.
\end{proof}

Thm \ref{lambda} and Prop. \ref{DrinfeldtoDrinfeld} are true in the case that $K$ has positive characteristic;  see \cite{WeinsteinFormalVectorSpaces}, Prop. 4.4.1.

From Prop. \ref{DrinfeldtoDrinfeld} we deduce the existence of a morphism of formal schemes $\MM_{G_0,m}^{(0)}\to\MM_{\wedge^n G_0,m}^{(0)}$.  Recall from \S\ref{height1} that $\MM_{\wedge^n G_0,m}^{(0)}=\Spf\OO_{K_n}$.  After passing to the geometric generic fiber, it breaks up as the union of $q^{m-1}(q-1)$ points.  The fibers of $\MM_{G_0,m}^{(0)}$ over each of these points are connected;  this is the main result of \cite{Strauch}.

\subsection{Determinants of formal vector spaces}
\label{detFVS}
For an object $R$ of $\Alg_{\overline{\FF}_q}$, one has an injection
\[ \varprojlim_m G_0[\pi^m](R)\injects \varprojlim G_0(R)=\tilde{G}_0(R). \]
Since $\tilde{G}_0(R)$ is a $K$-vector space, and $\varprojlim_m G_0[\pi^m](R)$ is a torsion-free $\OO_K$-module, this extends to an injective map of $K$-vector spaces
\[ \varprojlim_m G_0[\pi^m](R)\otimes K\injects \tilde{G}_0(R).\]

\begin{lemma} \label{Rdiscrete} Suppose $R$ is discrete.  Then $\varprojlim_m G_0[\pi^m](R)\otimes_{\OO_K} K \to \tilde{G}_0(R)$ is an isomorphism.
\end{lemma}

\begin{proof} The only thing to check is surjectivity.  Since $R$ is discrete, every element of $G_0(R)$ is $\pi^m$-torsion for some $m$.  Let $x=(x_0,x_1,\dots)\in \tilde{G}_0(R)$.  If $\pi^mx_0=0$, then $\pi^mx$ lies in $\varprojlim_m G_0[\pi^m](R)$.
\end{proof}

Now let $R$ be an object of $\Alg_{\OO_{K_0}}$ with ideal of definition $I$ which we assume contains $\pi$.  Then $R/I$ is a discrete $\overline{\FF}_q$-algebra.  Let $G$ be a lift of $G_0\otimes_{\overline{\FF}_q} R/I$ to $R$.  Thm. \ref{lambda} applied to $G_0$ shows that there exists an alternating $\OO_K$-multilinear map
\[  \lambda_m\from G_0[\pi^m](R/I)^n\to \wedge^n G_0[\pi^m](R/I). \]
Taking inverse limits and tensoring with $K$, we get an alternating $K$-multilinear map
\[ \lambda\from \varprojlim_m G_0[\pi^m](R/I)^n\otimes_{\OO_K} K \to \varprojlim_m\wedge^n G_0[\pi^m](R/I)\otimes_{\OO_K} K. \]
By Lemmas \ref{crystalline} and \ref{Rdiscrete}, there are isomorphisms
\[ \tilde{G}(R) \isom \tilde{G}_0(R/I)\isom\varprojlim_m G_0[\pi^m](R/I)\otimes_{\OO_K}K \]
(and similarly for $\wedge^n G$).  Thus we have defined an alternating $K$-multilinear map
\[ \lambda\from \tilde{G}^n\to \widetilde{\wedge^n G} \]
of formal $K$-vector spaces over $\OO_{K_0}$.


After choosing coordinates on $G$ and $\wedge^n G$, we get isomorphisms $\tilde{G}\isom \Spf \OO_{K_0}\powerseries{X^{1/q^\infty}}$ and $\wedge^n G\isom \Spf\OO_{K_0}\powerseries{T^{1/q^\infty}}$.  The morphism $\lambda$ above amounts to having an element
\[ \delta(X_1,\dots,X_n)\in \OO_{K_0}\powerseries{X_1^{1/q^\infty},\dots,X_n^{1/q^\infty}}\]
which comes equipped with a distinguished family of $q^m$th power roots for $m=1,2,\dots$.  These will simply be written $\delta(X_1,\dots,X_n)^{1/q^m}$.

\subsection{The structure of the Lubin-Tate moduli problem at infinite level}

\begin{defn}  Write $\MM_{G_0,m}^{(j)}=\Spf A_m^{(j)}$.  Let $A^{(j)}$ be the completion of $\varinjlim A_m^{(j)}$ with respect to the topology induced by the maximal ideal of $A_0^{(j)}$ (or any $A_m^{(j)}$, it doesn't matter). Let $\MM_{G_0,\infty}^{(j)}=\Spf A^{(j)}$, and let $\MM_{G_0,\infty}=\coprod_{j\in\Z}\MM_{G_0,\infty}^{(j)}$.  $\MM_{G_0,\infty}$ is the {\em Lubin-Tate deformation space at infinite level}.   Write $A=A^{(0)}$.
\end{defn}

\begin{rmk} The completion of a non-noetherian ring at an ideal $I$ is not necessarily $I$-adically complete.  However, this is true if $I$ is finitely generated.  Thus the $A^{(j)}$ are complete, and $\Spf A^{(j)}$ makes sense as a formal scheme.  We have $\Spf A^{(j)}=\varprojlim \Spf A_m^{(j)}$ in the category of formal schemes over $\Spf \OO_{K_0}$.
\end{rmk}

Recall that if $G_0$ has height $n$, then $\wedge^nG_0$ has height 1.   By \S\ref{height1} we have $\MM^{(0)}_{\wedge^n G_0,m}=\Spf \OO_{K_m}$, where $K_m/K_0$ is the totally ramified abelian extension of degree $q^{m-1}(q-1)$.  Let $\hat{K}_\infty$ be the $\pi$-adic completion of $K_\infty=\bigcup_m K_m$.  Then $\MM^{(0)}_{\wedge^n G_0,\infty}=\Spf \OO_{\hat{K}_\infty}$.

In \S\ref{detmodules} we constructed a morphism $\MM_{G_0,m}^{(0)}\to \MM_{\wedge^n G_0,m}^{(0)}$.  Taking inverse limits with respect to $m$, we get a morphism $\MM_{G_0,\infty}^{(0)}\to\MM_{\wedge^n G_0,\infty}^{(0)}$.  

Let $G^{\univ}$ be the universal deformation of $G_0$ to $A_0$.  Then over $A_m$, we have a universal Drinfeld basis $X_1^{(m)},\dots,X_n^{(m)}\in G^{\univ}[\pi^m](A_m)$.  In the limit, we get $n$ distinguished elements
\begin{equation}
\label{X_i}
\mathbf{X}_1,\dots,\mathbf{X}_n\in \tilde{G}^{\univ}(A).
\end{equation}

Now suppose that $G$ is an arbitrary lift of $G_0$ to $\OO_{K_0}$.  Let $I_0\subset A_0$ be the maximal ideal, so that $A_0/I_0=\overline{\FF}_q$.   We have $G\otimes_{\OO_{K_0}} \overline{\FF}_q = G^{\univ}\otimes_{A_0} A_0/I_0 =G_0$.  Twice applying part (1) of Prop. \ref{crystalline}, we get isomorphisms
\[ \tilde{G}^{\univ}(A)\isom \tilde{G}^{\univ}(A/I_0) = \tilde{G}_0(A/I_0) =\tilde{G}(A/I_0)\isom \tilde{G}(A). \]
Let $\mathbf{Z}_i$ be the image of $\mathbf{X}_i$ under the above isomorphism.  By unwinding the proof of Lemma \ref{crystalline}, we can say what these are explicitly.  Choose coordinates on $G$ and $G^{\univ}$, so that $G(A)$ and $G^{\univ}(A)$ may be identified with $\Nil(A)$.
Then $\mathbf{Z}_i=(Z_i^{(1)},Z_i^{(2)},\dots)\in\varprojlim G(A)$, where
\begin{equation}
\label{Zim}
 Z_i^{(m)}=\lim_{r\to\infty} [\pi^r]_G\left(X_i^{(r+m)}\right).
 \end{equation}

The tuple $(\mathbf{Z}_1,\dots,\mathbf{Z}_n)$ represents an $A$-point of $\tilde{G}$, which is to say a morphism of formal schemes $\MM_{G_0}^{(0)}\to \tilde{G}^n$ over $\Spf \OO_{K_0}$.  Recall by Prop. \ref{crystalline}, $\tilde{G}$ is representable by a formal scheme isomorphic to $\Spf \OO_{K_0}\powerseries{X^{1/q^\infty}}$.  Thus in fact we have a continuous $\OO_{K_0}$-algebra homomorphism $\OO_{K_0}\powerseries{X^{1/q^\infty}}\to A$ which sends $X_i^{1/q^m}$ to $Z_i^{1/q^m}$.

Applying the same constructions to $\wedge^n G_0$, we have a morphism of formal schemes $\MM_{\wedge^n G_0,\infty}\to \widetilde{\wedge^n G,\infty}$.  By the naturality of the determinant morphism, the diagram
\begin{equation}
\label{Mdiagram}
\xymatrix{
\mathcal{M}_{G_0,\infty}^{(0)} \ar[r] \ar[d] & \mathcal{M}_{\wedge^n G_0,\infty}^{(0)} \ar[d] \\
\tilde{G}^n \ar[r]_{\lambda} & \widetilde{\wedge^n G}
}
\end{equation}
commutes.

\begin{Theorem}\label{Cartesian}
The above diagram is Cartesian.  That is, $\mathcal{M}_{G_0,\infty}^{(0)}$ is isomorphic to the fiber product of $\tilde{G}^n$ and $\mathcal{M}_{\wedge^n G_0,\infty}^{(0)}$ over $\widetilde{\wedge^n G}$.
\end{Theorem}

We remark there is a similar diagram for the entire space $\MM_{G_0,\infty}$.  

\subsection{Proof of Thm. \ref{Cartesian}}

The fiber product of $\tilde{G}^n$ and $\mathcal{M}^{(0)}_{\wedge^n G_0,\infty}$ over $\widetilde{\wedge^n G}$ is an affine formal scheme, say $\Spf B$, where $B$ is a complete local ring.  We have a homomorphism of local rings $\phi\from A\to B$ which we claim is an isomorphism.  We need a few lemmas.

\begin{lemma}\label{trivialbasis} Let $R$ be an object of $\mathcal{C}$ in which $\pi=0$.  Any $n$-tuple of elements in $G_0[\pi^{m-1}](R)$ constitutes a Drinfeld basis for $G_0[\pi^m](R)$.  Similarly, any element in $\wedge^n G_0[\pi^{m-1}](R)$ constitutes a Drinfeld basis for $\wedge^n G_0[\pi^m](R)$.
\end{lemma}

\begin{proof} The claim for $G_0$ is equivalent to the assertion that the $n$ elements $0,\dots,0$ constitute a Drinfeld basis for $G_0[\pi](R)$.  This in turn is equivalent to the assertion that $T^{q^n}$ be divisible by $[\pi]_{G_0}(T)$ in $R\powerseries{T}$.  But $[\pi]_{G_0}(T)$ equals $T^{q^n}$ times a unit in $\overline{\FF}_q\powerseries{T}$, because $G_0$ has height $n$.  The claim for $\wedge^n G_0$ is proved similarly.
\end{proof}

Recall the parameters $X_1^{(m)},\dots,X_n^{(m)}\in A_m$, which represent the universal Drinfeld basis for the $G^{\univ}[\pi^m](A_m)$.   Let $I\subset A_1$ be the ideal generated by $(X_1^{(1)},\dots,X_n^{(1)})$, which is to say that $I$ is the maximal ideal of $A_1$.  We will often be considering the extension of $I$ to the rings $A_m$ and $A$, and we will abuse notation in calling these ideals $I$ as well.  Note that $I\subset A_1$ is the maximal ideal of $A_1$, so that $A_1/I=\overline{\FF}_q$.  In particular $\pi\in I$.

Recall that $I_0$ is the maximal ideal of $A_0$.  Thus $I_0\subset I$.  In fact:
\begin{lemma}
\label{I0I2} $I_0\subset I^2$.
\end{lemma}

\begin{proof} $X_1^{(1)},\dots,X_n^{(1)}$ is a Drinfeld basis for $G^{\univ}[\pi](A_1)$.  Thus the polynomial
\[ \prod_{(a_1,\dots,a_n)\in k^n} \left(T - ([a_1]_{G_{\univ}}(X_1^{(1)})+_{G^{\univ}}\dots+_{G^{\univ}}[a_n]_{G^{\univ}}(X_n^{(1)})\right) \]
is divisible by $[\pi]_{G^{\univ}}(T)$ in $A_1\powerseries{T}$.  This product is congruent to $T^{q^n}$ modulo $I^2$.  Now we apply Lemma \ref{generate}.  Since the coefficient of $T^{q^n}$ in $[\pi]_{G^{\univ}}(T)$ is a unit, we find that the the coefficients of $T,T^2,\dots,T^{q^n-1}$ in $[\pi]_{G^{\univ}}(T)$ must lie in $I^2$.  But these coefficients generate $I_0$, whence the lemma.
\end{proof}

\begin{lemma}  \label{modI2} In $A_1\powerseries{T}$, The congruence $[\pi]_{G^{\univ}}(T)\equiv [\pi]_{G}(T)$ holds modulo $I^2\powerseries{T}$.
\end{lemma}

\begin{proof}  Indeed, both sides of the congruence lie in $A_0\powerseries{T}$ and are both congruent to $[\pi]_{G_0}(T)$ modulo $I_0\powerseries{T}$, so this follows from Lemma \ref{I0I2}.
\end{proof}

The next lemma describes the closed subscheme $\Spec A_m/I$ of the formal scheme $\MM_{G_0,m}=\Spf A_m$.

\begin{lemma} \label{AmI} There is an isomorphism of affine $k$-schemes $\Spec A_m/I\to G_0[\pi^{m-1}]^n$.
\end{lemma}

\begin{proof} For an object $R$ of $\mathcal{C}$ in which $\pi=0$, we have that $\Hom_{\mathcal{C}}(A_m/I,R)$ is the set of deformations $G'$ of $G_0$ to $R$ equipped with a Drinfeld basis $x_1,\dots,x_n$ for $G'[\pi^m](R)$ which satisfy $\pi^{m-1}x_i=0$, $i=1,\dots,n$.  For such a deformation we have
\[ G' = G^{\univ} \otimes_{A_0} R = (G^{\univ}\otimes_{A_0} A_0/I_0) \otimes_{A_0/I_0} R = G_0\otimes_{\overline{\FF}_q} R \]
Thus $\Hom_{\mathcal{C}}(A_m/I,R)$ is the set of Drinfeld bases $x_1,\dots,x_n$ for $G_0[\pi^m](R)$ which satisfy $\pi^{m-1}x_i=0$, $i=1,\dots,n$; that is, $x_1,\dots,x_n\in G_0[\pi^{m-1}](R)$.  But by Lemma \ref{trivialbasis}, any such $n$-tuple is automatically a Drinfeld basis.  Thus $\Hom_{\mathcal{C}}(A_m/I,R)$ is simply the set of $n$-tuples of elements of $G_0[\pi^{m-1}](R)$.  This identifies $\Spec A_m/I$ with $G_0[\pi^{m-1}]^n$.
\end{proof}

We now turn to $B$, which by definition is the coordinate ring of the affine formal scheme $\tilde{G}_0^n \times_{\widetilde{\wedge^n G}} \MM_{\wedge^n G_0}^{(0)}$.

\begin{lemma} \label{Bsemiperfect} The $q$th power Frobenius map is surjective on $B/\pi$.
\end{lemma}

\begin{proof} We have the following presentation of $B$:
\begin{equation}
\label{Bpresentation}
B\approx \OO_{K_0}\powerseries{X_1^{1/q^\infty},\dots,X_n^{1/q^\infty}} \hat{\otimes}_{\OO_{K_0}\powerseries{X^{1/q^\infty}}} \OO_{\hat{K}_{\infty}}
\end{equation}
Since the Frobenius map is surjective on $\OO_{\hat{K}_{\infty}}/\pi$ (Lemma \ref{OKinfty}) and on $\OO_{K_0}/\pi=\overline{\FF}_q$ , it is surjective on $B/\pi$.
\end{proof}

For an object $R$ of $\Alg_{\OO_{K_0}}$, $\Hom(B,R)$ is in bijection with the set of $n$-tuples $x_1,\dots,x_n\in \tilde{G}(R)$ such that $\lambda(x_1,\dots,x_n)$, a priori just an element of $\widetilde{\wedge^n G}(R)$, actually lies in $T(\wedge^n G)(R)=\varprojlim_m \wedge^n G[\pi^m](R)$, and constitutes a Drinfeld basis for each $\wedge^n G[\pi^m](R)$.  The identity homomorphism $\Hom(B,B)$ corresponds to an $n$-tuple of universal elements $Y_1,\dots,Y_n\in\tilde{G}(B)$.  For $i=1,\dots,n$, let us write $Y_i=(Y_i^{(1)},Y_i^{(2)},\dots)$.  After choosing a coordinate on $G$, we can identify $Y_i^{(m)}$ with a (topologically nilpotent) element of $B$.

Let $J\subset B$ be the ideal generated by $\pi$ and by $Y_1^{(1)},\dots,Y_n^{(1)}$.

\begin{lemma}  \label{BJAI} $\phi(J)\subset I$, and $\phi$ descends to an isomorphism $B/J\to A/I$.
\end{lemma}

\begin{proof}
For an $\OO_{K_0}$-algebra $R$ in which $\pi=0$, $\Hom(B/J,R)$ is in bijection with the set of $n$-tuples $x_1,\dots,x_n\in \tilde{G}_0(R)$ such that (a) $x_i^{(1)}=0$ for $i=1,\dots,n$ and such that (b) $\lambda(x_1,\dots,x_n)$ constitutes a compatible family of Drinfeld bases for $\wedge^n G_0[\pi^m](R)$.
However, if condition (a) is satisfied, then
\[ \lambda(x_1,\dots,x_n)^{(1)} = \lambda_1(x_1^{(1)},\dots,x_n^{(1)}) =0, \]
and $0$ is always a Drinfeld basis for $\wedge^n G[\pi](R)$ by Lemma \ref{trivialbasis}, so that condition (b) is satisfied as well.

Therefore $\Spec B/J=\varprojlim_m G[\pi^{(m-1)}]^n$.  By Lemma \ref{AmI} this is isomorphic to $\varprojlim_m \Spec A_m/I = \Spec A/I$.
\end{proof}

\begin{lemma} \label{phiJ} $\phi(J)A$ is a dense subset of $I$.
\end{lemma}

\begin{proof} By Eq. \eqref{Zim} we have
\[ \phi(Y_i^{(r)}) = \lim_{m\to \infty} [\pi^m]_{G}\left(X_i^{(m+r)}\right). \]
Since the limit converges $I$-adically, $\phi\left(Y_i^{(1)}\right)\equiv[\pi^{m-1}]_{G_0}\left(X_i^{(m)}\right)\pmod{I^2}$ for some sufficiently large $m$.  By Lemma \ref{modI2} we have
\[ X_i^{(1)} = [\pi^{m-1}]_{G^{\univ}}\left(X_i^{(m)}\right) \equiv [\pi^{m-1}]_{G}\left(X_i^{(m)}\right)\equiv \phi\left(Y_i^{(1)}\right)\pmod{I^2}, \]
and therefore $X_i^{(1)}\in \phi(J)A+I^2$.  Since the $X_i^{(1)}$ generate $I$, we have $I\subset \phi(J)A+I^2$, which when iterated yields $I\subset \phi(J)A+I^m$ for all $m\geq 1$.  Since $I$ generates the topology on $A$, the closure of $\phi(J)A$ must equal $I$.
\end{proof}

\begin{lemma} \label{Aalmostsemiperfect} The $q$th power Frobenius map on $A/\pi$ has dense image.
\end{lemma}

\begin{proof} Let $\overline{A}=A/\pi$, and let $\overline{I}$ be the image of $I$ in $\overline{A}$, so that $\overline{I}$ generates the topology on $\overline{A}$.  Similarly define $\overline{B}$ and $\overline{J}$.  By Lemma \ref{BJAI} we have $\overline{A}/\overline{I}\isom \overline{B}/\overline{J}$, so the $q$th power Frobenius map is also surjective on $\overline{A}/\overline{I}$.  Thus $\overline{A}=\overline{A}^q+\overline{I}$.

We will prove by induction that for all $m\geq 1$, $\overline{A}=\overline{A}^q+\overline{I}^m$ and $\overline{I}=(\overline{A}^q\cap\overline{I})+\overline{I}^m$.  The first claim proves the lemma, since $\overline{I}$ generates the topology on $\overline{A}$.  As for the base case $m=1$, the first claim is discussed above, and the second claim is vacuous.  For the induction step, assume both claims for $m$.  By Lemma \ref{phiJ}, $\phi(\overline{J})\overline{A}$ is dense in $\overline{I}$, so that $\overline{I}=\phi(\overline{J})\overline{A}+\overline{I}^{m+1}$.  Since Frobenius is surjective on $\overline{B}$ (Lemma \ref{Bsemiperfect}), we have $\phi(\overline{J})\subset \overline{A}^q\cap \overline{I}$.  Thus
\begin{eqnarray*}
\overline{I}&\subset& (\overline{A}^q\cap \overline{I})\overline{A}+\overline{I}^{m+1} \\
&=& (\overline{A}^q\cap \overline{I})(\overline{A}^q+\overline{I}^m)+\overline{I}^{m+1}\\
&\subset& (\overline{A}^q\cap \overline{I})+\overline{I}^{m+1}.
\end{eqnarray*}
The reverse containment is obvious, so that $\overline{I}=(\overline{A}^q\cap \overline{I})+\overline{I}^{m+1}$, thus establishing the second claim for $m+1$.
Inserting this into $\overline{A}=\overline{A}^q+\overline{I}$ gives $\overline{A}=\overline{A}^q+\overline{I}^{m+1}$, which establishes the first claim for $m+1$.
\end{proof}

\begin{lemma} \label{Asemiperfect} The $q$th power Frobenius map on $A/\pi$ is surjective.
\end{lemma}

\begin{proof}  Once again let $\overline{A}=A/\pi$ and let $\overline{I}$ be the image of $I$ in $\overline{A}$.  The ideal $\overline{I}$ is finitely generated;  let $f_1,\dots,f_n$ be a set of generators ({\em e.g.}, the images of the elements $X_i^{(1)}$, $i=1,\dots,n$).   Recall that for $m\geq 0$, $\overline{I}^{[q^m]}$ is the ideal generated by the $q^m$th powers of elements of $\overline{I}$, so that $\overline{I}^{[q^m]}$ is generated by the $f_i^{q^m}$, $i=1,\dots,n$.  Obviously we have $\overline{I}^{[q^m]}\subset \overline{I}^{q^m}$.  But also we have $\overline{I}^{q^N}\subset \overline{I}^{[q^m]}$ for $N$ large enough.  Thus the sequence of ideals $\overline{I}^{[q^m]}$ generates the topology on $A$.

Let $a\in \overline{A}$.  By Lemma \ref{Aalmostsemiperfect} there exists $b\in \overline{A}$ such that $a-b^q\in \overline{I}^{[q]}$. Let us write
\[ a = b^q + \sum_i a_if_i^q,\;a_i\in \overline{A}. \]
For each $i$ we may also find $b_i\in \overline{A}$ with $a_i-b_i^q\in \overline{I}^{[q]}$; write $a_i=b_i^q+\sum_j a_{ij}f_j^{q}$, $a_{ij}\in A$.  Thus
\[ a  =b^q + \sum_i b_i^qf_i^q + \sum_{i,j} a_{ij}f_i^qf_j^{q}.\]
Continuing this process, we find a $q$th root of $a$ in $\overline{A}$, namely
\[ b + \sum_i b_if_i + \sum_{i,j} b_{ij}f_if_j + \dots. \]
This completes the proof.
\end{proof}

\begin{lemma} \label{Surjmodpi} The induced map $B/\pi\to A/\pi$ is surjective.
\end{lemma}

\begin{proof}  Let us write $\overline{\phi}\from \overline{B}\to\overline{A}$ for the induced map.  Let $a\in \overline{A}$.  By Lemma \ref{BJAI}, there exists $b_0\in \overline{B}$ and $a_0\in \overline{I}$ with $a=\phi(b_0)+a_0$.  By Lemma \ref{Asemiperfect}, $a_0$ has a $q$th root in $\overline{A}$, call it $a_0^{1/q}$.  Apply Lemma \ref{BJAI} to write $a_0^{1/q} = \phi(b_1)+a_1$, with $b_1\in \overline{B}$, $a_1\in \overline{I}$.  Then $\phi(b_1^q)=a_0-a_1^q\in\overline{I}$, so that (by Lemma \ref{BJAI}) $b_1^q\in \overline{J}$.  Similarly write $a_1^{1/q}=\phi(b_2)+a_2$, and so on.  We have $b_m^q\in J$ for all $m\geq 1$.  Therefore the series $b_0+b_1^q+b_2^{q^2}+\dots$ converges to an element $b$ with $\phi(b)=a$.
\end{proof}

Recall that $A_m$ is a regular local ring.  Let $I_m$ be the maximal ideal of $A_m$.  Then for $m\geq 1$, $A_m$ is generated by the Drinfeld parameters $X_1^{(m)},\dots,X_n^{(m)}$.  We have $[\pi]_{G^{\univ}}\left(X_i^{(m+1)}\right)=X_i^{(m)}$ for $m\geq 1$.  Also recall that we had set $I=I_1$.

\begin{lemma} \label{GrowthOfMthIdeal}  As ideals in $A_m$ we have $IA_m\subset I_m^{[q^{n(m-1)]}}$.  Furthermore, $I_m^{nq^{n(m-1)}}\subset IA_m$.
\end{lemma}

\begin{proof} From Lemma \ref{modI2} we have $[\pi]_{G^{\univ}}(T)\equiv [\pi]_{G_0}(T)$ modulo $I^2\powerseries{T}$.  Since $G_0$ has height $n$, $[\pi]_{G^{\univ}}(T)$ is congruent to a power series in $T^{q^n}$ modulo $I_0\powerseries{T}\subset I^2\powerseries{T}$ (see Lemma \ref{I0I2}).  Thus
\[ \left(X_i^{(m)}\right)^{q^{n(m-1)}}\equiv [\pi^{m-1}]_{G_0}\left(X_i^{(m)}\right)\equiv [\pi^{m-1}]_{G^{\univ}}\left(X_i^{(m)}\right)\equiv X_i^{(1)}\pmod{I^2A_m}. \]
Since the $X_i^{(1)}$ generate $I$ we have $I\subset I_m^{[q^{n(m-1)}]}+I^2A_m$.  Iterating this containment shows that $I\subset I_m^{[q^{n(m-1)}]}+I^NA_m$ for any $N$.  $IA_m$ generates the topology on $A_m$, and $I_m^{[q^{n(m-1)}]}$ is open, so that $I^N\subset I_m^{(q^n(m-1)}$ for sufficiently large $N$.  Thus $IA_m\subset I_m^{[q^{n(m-1)}]}$.

For the second claim, note that $I_m$ has $n$ generators $X_1^{(m)},\dots,X_n^{(m)}$, so that $I_m^{nq^{n(m-1)}}$ is contained in the ideal generated by the $\left(X_i^{(m)}\right)^{q^{n(m-1)}}$.  The above congruence shows that each of these elements lies in $IA_m$.  
\end{proof}

\begin{lemma} \label{GradedRingArgument} Let $R$ be a regular local ring with maximal ideal $M$ whose residue field has characteristic $p$.  Let $Q$ be a power of $p$.  Let $\eta\in R$ be an element for which $\eta^Q$ divides $p$ in $R$.  Suppose $N\geq 1$.  If $x^Q\in \eta^QR+M^{QN+1}$, then $x\in \eta R+M^{N+1}$.
\end{lemma}

\begin{proof}  For $x\in R$, let $v(x)\in\Z_{\geq 0}\cup\set{\infty}$ denote the maximal $s\geq 0$ such that $x\in M^s$.  If $x$ satisfies the hypothesis of the lemma, then every element of $x+\eta R$ satisfies the hypothesis as well.  Indeed if $x'=x+\eta y$ for $y\in R$, then $(x')^Q\in x^Q+pR+\eta^QR\subset x^Q+\eta^QR\in \eta^QR+M^{QN+1}$.  Similarly, if $x$ satisfies the conclusion of the lemma, then so does every element of $x+\eta R$.  Thus to prove the lemma, we assume that $x$ satisfies $v(x)\geq v(x+\eta y)$ for all $y\in R$.  Under this assumption, we will show that $x\in M^{N+1}$, which suffices.

Assume otherwise, so that $v(x)\leq N$.  Consider the graded ring $\text{Gr} R=\bigoplus_{i\geq 0}M^i/M^{i+1}$.  Since $R$ is a regular local ring, $\text{Gr} R$ is isomorphic to a power series ring over the residue field of $R$. If $z\in R$, define $\overline{z}\in \text{Gr} R$ as follows.  If $z\in R$ is nonzero, let $\overline{z}\in \text{Gr}R$ denote the image of $z$ in $I_m^{v(z)}/I_m^{v(z)+1}$.  Define $\overline{0}=0$.  We have $x^Q\in\eta^QR +M^{QN+1}$ and $v(x^Q)\leq QN$.  This implies that $\overline{\eta}^Q$ divides $\overline{x}^Q$ in $\text{Gr}R$, and therefore (since this ring is a unique factorization domain) $\overline{\eta}$ divides $\overline{x}$.
Thus we can find $y\in R$ with $v(x-\eta y)>v(x)$.  But this contradicts our assumption about $x$.
\end{proof}

\begin{lemma} \label{AIsIntegrallyClosed} For every $r\geq 1$, $A$ contains an element $\eta_r$ with $\eta_r^{q^r}A=\pi A$.  If $f\in A$ satisfies $f^{q^r}\in \pi A$, then $f\in\eta_rA$.
\end{lemma}

\begin{proof}  We have seen that $A_r$ contains the ring of integers $\OO_{K_r}$ in the Lubin-Tate extension $K_r/K_0$.  A uniformizer $\pi_r$ of $\OO_{K_r}$ satisfies $\pi_r^{q^{r-1}(q-1)}\OO_{K_r}=\pi\OO_{K_r}$.  Thus $\eta_r=\pi_{r+1}^{q-1}$ satisfies $\eta_r^{q^r}A=\pi A$.

Now suppose that $f\in A$ satisfies $f^{q^r}\in \pi A$.  Let $N\geq 1$ be arbitrary.  Since $A$ is the $I$-adic completion of the direct limit of the $A_m$, there exist $m\geq 1$ and $f_m\in A_m$ such that $f-f_m\in I^{nq^rN+1}A$.  Then $f_m^{q^r}\in \pi A_m+ I^{nq^rN+1}A_m$.  By Lemma \ref{GrowthOfMthIdeal} we have $f_m^{q^r}\in \pi A_m+I_m^{[q^{n(m-1)}](nq^rN+1)}\subset \pi A_m+I_m^{nq^{r+n(m-1)}N+1}$.  After possibly enlarging $m$ we may assume that $\eta_r\in A_m$.  Applying Lemma \ref{GradedRingArgument} we find $f_m\in \eta_rA_m+I_m^{nq^{n(m-1)}N+1}$.  Applying Lemma \ref{GrowthOfMthIdeal} again, we get $f_m\in\eta_rA_m+I^NA$.  Thus $f\in \eta_rA+I^NA$.  Since $N$ was arbitrary, $f$ lies in the closure of $\eta_r A$, which is $\eta_r A$ itself.
\end{proof}

We are now ready to show that that $\phi\from B\to A$ is an isomorphism.  We will first show it is surjective.  If $a\in A$, use Lemma \ref{Surjmodpi} to find $b\in B$ with $a=\phi(b_0)+\pi a_1$, $a_1\in A$.  Write $a_1=\phi(b_1)+\pi a_2$, and so on.  Then $b=b_0+\pi b_1+\pi^2 b_2+\dots$ satisfies $\phi(b)=a$.

We now turn to injectivity.  Suppose that $b\in B$ is an element with $\phi(b)=0$.  Let $m\geq 1$.
Since Frobenius is surjective on $B/\pi$ (Lemma \ref{Bsemiperfect}), we may write $b\equiv c^{q^m}\pmod{\pi B}$, with $c\in B$.  Then $\phi(c)^{q^m}\in \pi A$.  By Lemma \ref{AIsIntegrallyClosed}, $\phi(c)\in \eta_mA$.  Since $\phi$ is an isomorphism $B/J\to A/I$ (Lemma \ref{BJAI}), we have $c\in \eta_mB+J$, and therefore $b\in \pi B + J^{[q^m]}$ for all $m\geq 1$.  Since $B$ is $J$-adically complete, this implies that $b$ lies in the closure of $\pi B$, which is $\pi B$ itself.   This shows that $\ker\phi\subset \pi B$.  But then $b/\pi\in \ker\phi$, so that in fact $b\in \pi^2B$.  Inductively, we find $b\in\pi^mB$ for all $m\geq 1$.  Since $B$ is $\pi$-adically separated, $b=0$.  This completes the proof of Theorem \ref{Cartesian}.

\begin{Cor} \label{Ainf} If $K$ has positive characteristic, then
\[ A\isom \overline{\FF}_q\powerseries{X_1^{1/q^\infty},\dots,X_n^{1/q^\infty}}.\]
If $K$ has characteristic 0, then put $A^{\flat}=\varprojlim A/\pi$, where the limit is taken with respect to the $q$th power Frobenius map.  Then
\[ A^{\flat}\isom \overline{\FF}_q\powerseries{X_1^{1/q^\infty},\dots,X_n^{1/q^\infty}} \]
\end{Cor}

\begin{proof}  By Thm. \ref{Cartesian}, we have an isomorphism of complete local $\OO_{\hat{K}_\infty}$-algebras
\begin{equation}
\label{tensor}
\OO_{K_0}\powerseries{X_1^{1/q^\infty},\dots,X_n^{1/q^\infty}}\hat{\otimes}_{\OO_{K_0}\powerseries{X^{1/q^\infty}}} \OO_{\hat{K}_\infty} \to A
\end{equation}
In the tensor product appearing in Eq. \eqref{tensor}, the image of $X\in \OO_{K_0}\powerseries{X^{1/q^\infty}}$ is $\delta(X_1,\dots,X_n)$ in the left factor and $t$ in the right.
First assume that $K$ has positive characteristic.  Then the map $\OO_{K_0}\powerseries{X^{1/q^\infty}}\to \OO_{\hat{K}_\infty}$ is surjective with kernel generated by $\pi-g(X)$ for some fractional power series $g(X)\in\overline{\FF}_q\powerseries{X^{1/q^\infty}}$ without constant term.  Recalling that $\OO_{K_0}=\overline{\FF}_q\powerseries{\pi}$, we have
\begin{eqnarray*}
A &\isom& \OO_{K_0}\powerseries{X_1^{1/q^\infty},\dots,X_n^{1/q^\infty}}/(\pi-g(\delta(X_1,\dots,X_n))) \\
&=& \overline{\FF}_q\powerseries{X_1^{1/q^\infty},\dots,X_n^{1/q^\infty}}.
\end{eqnarray*}

Now assume $K$ has characteristic 0.  We have
\[ A/\pi \isom \overline{\FF}_q\powerseries{X_1^{1/q^\infty},\dots,X_n^{1/q^\infty}}\hat{\otimes}_{\overline{\FF}_q\powerseries{X^{1/q^\infty}}}
\OO_{K_\infty}/\pi. \]
Now take the inverse limit along the $q$th power Frobenius maps.  In doing so, the surjection $\overline{\FF}_q\powerseries{X^{1/q^\infty}} \to \OO_{K_\infty}/\pi$ becomes an isomorphism $\overline{\FF}_q\powerseries{X^{1/q^\infty}}\to \OO_{K_\infty^{\flat}}$.  Thus $A^\flat = \overline{\FF}_q\powerseries{X_1^{1/q^\infty},\dots,X_n^{1/q^\infty}}$ as required.
\end{proof}

\subsection{Group actions and geometrically connected components}
Thm. \ref{Cartesian} extends to the entire formal scheme $\MM_{G_0,\infty}$, so that we get a cartesian diagram
\begin{equation}
\label{Mdiagram}
\xymatrix{
\mathcal{M}_{G_0,\infty} \ar[r] \ar[d] & \mathcal{M}_{\wedge^n G_0,\infty} \ar[d] \\
\tilde{G}^n \ar[r]_{\lambda} & \widetilde{\wedge^n G}
}
\end{equation}
The action of $\GL_n(K)\times D^\times$ on $\mathcal{M}_{G_0,\infty}$ can be described directly in terms of this diagram.  $\GL_n(K)$ acts on the right on $\tilde{G}^n$ via the usual right action of matrices on row vectors.  The determinant morphism $\lambda$ transforms the action of $g_1\in \GL_n(K)$ into $\det g_1\in K^\times$, which preserves $\MM_{\wedge^n G_0,\infty}$.  Thus $g_1$ preserves $\MM_{G_0,\infty}$.  

Recall that $D=\End G_0\otimes_{\OO_K} K$.  An element $g_2\in D^\times$ determines an autmorphism of $\tilde{G}_0$, hence of $\tilde{G}$ by Prop. \ref{crystalline}.  The determinant morphism transforms the action of $g_2$ into $N(g_2)\in K^\times$, where $N\from D\to K$ is the reduced norm map.  Thus $g_2$ preserves $\MM_{G_0,\infty}$ as well.  However, in order to get a right action of $D^\times$ on $\MM_{G_0,\infty}$ which is consistent with the previously described action on the finite-level spaces, we define  
\[ (\X_1,\dots,\X_n)g_2 = (g_2^{-1}\X_1,\dots,g_2^{-1}\X_n), \]
whenever $\X_1,\dots,\X_n$ are sections of $\tilde{G}$.  

Let $\C$ be the completion of an algebraic closure of $K$.  Choose an embedding $\hat{K}_\infty\injects\C$.  This is tantamount to choosing a generator $\mathbf{t}$ for the free rank one $\OO_K$-module $\varprojlim \wedge^nG[\pi^m](\OO_\C)$.  The inclusion $\varprojlim \wedge^nG[\pi^m](\OO_\C)\subset\varprojlim \wedge^nG(\OO_\C)$ allows us to view $\mathbf{t}$ as an element of $\widetilde{\wedge^nG}(\OO_\C)$.  This in turn corresponds (see Rmk. \ref{boldfaceX}) to a topologically nilpotent element $t\in\OO_\C$ together with a compatible system of $q$th power roots $t^{1/q^m}$ for $m\geq 1$.

Let
\[ A_m^{\circ}=A_m\hat{\otimes}_{\OO_{K_m}}\OO_\C, \]
and let $A^{\circ}$ be the completion of $\varinjlim_m A_m^{\circ}$.
Let $\MM^{\circ}_{G_0,\infty}=\Spf A^{\circ}$.  From Thm. \ref{Cartesian} we may identify $\MM^{\circ}_{G_0,\infty}$ with the fiber of $\lambda\from \tilde{G}^n_{\OO_\C}\to \widetilde{\wedge^n G}_{\OO_\C}$ over the single point $\mathbf{t}\in\widetilde{\wedge^n G}(\OO_\C)$.  Thus
\begin{equation}
\label{Acirc}
 A^{\circ}\isom \OO_{\C}\powerseries{X_1^{1/q^\infty},\dots,X_n^{1/q^\infty}}/\left(\delta(X_1,\dots,X_n)^{1/q^m}-t^{1/q^m}\right). \end{equation}
The $\OO_\C$-algebra $A^{\circ}$ admits an action of the group $(\GL_n(K)\times D^\times)^{\det=N}$ consisting of pairs $(g,b)$ with $\det g = N(b)$.

\subsection{The Lubin-Tate space at infinite level as a perfectoid space}
We wish to take the adic generic fiber of the formal scheme $\MM_{\infty}^{(0)}=\MM_{G_0,\infty}^{(0)}$.  The construction of this generic fiber is as follows.  Let $\MM_{\infty}^{(0),\ad}=\Spa(A,A)$ be the set of continuous valuations on $A$, as in \cite{Hub94}.  This is fibered over the two-point space $\Spa(\OO_K,\OO_K)$, and we can let $\MM^{(0),\ad}_{\infty,\eta}$ be the fiber over $\eta=\Spa(K,\OO_K)$.  The trouble with this is that the structure presheaf of $\Spa(A,A)$ isn't necessarily a sheaf, and therefore we don't know {\em a priori} that $\MM^{(0),\ad}_{\infty,\eta}$ is an adic space.  We will resolve this problem by observing that $\MM^{(0),\ad}_{\infty,\eta}$ is a perfectoid space after passing to $\C$.  It is known that the structure presheaf on a perfectoid affinoid is a sheaf (\cite{ScholzePerfectoidSpaces}, Thm. 6.3).

Let $\overline{\eta}=\Spa(\C,\OO_\C)$, and let $\MM^{(0),\ad}_{\infty,\overline{\eta}}$ be the base change to $\C$.

\begin{lemma} $\MM^{(0),\ad}_{\infty,\overline{\eta}}$ can be covered by perfectoid affinoids (and therefore is a perfectoid space).
\end{lemma}

\begin{proof} This is a consequence of the fact that $A_{\OO_\C}=A\hat{\otimes}\OO_\C$ is a reduced complete flat adic $\OO_\C$-algebra admitting a finitely generated ideal of definition $I$ containing $\pi$, such that the Frobenius map is surjective on $A_{\OO_\C}/\pi$.

Indeed, let $f_1,\dots,f_n$ be a set of generators for $I$ (such as the elements $X_i^{(1)}$).  Then any valuation $\abs{\;}$ belonging to $\MM_{\infty,\overline{\eta}}^{(0),\ad}$ must satisfy $\abs{f_i}<1$ for $i=1,\dots,n$.  Since $\abs{\pi}\neq 0$, there exists $r\geq 1$ for which $\abs{f_i^r}\leq \abs{\pi}$.  Let $R_r=A_{\OO_\C}\tatealgebra{f_i^r/\pi}[1/\pi]$, and let $R_r^+\subset R_r$ be the integral closure of $A_{\OO_\C}\tatealgebra{f_i^r/\pi}$ in $R_r$.  The Frobenius map is surjective on $R_r^+/\pi$ (this follows quickly from the corresponding property of $A$), and thus $(R_r, R_r^+)$ is a perfectoid $\OO_{\hat{K}_\infty}$-algebra, cf. \cite{ScholzePerfectoidSpaces}, Def. 6.1.  The valuation $\abs{\;}$ extends uniquely to $R_r$ and satisfies $\abs{R_r^+}\leq 1$, so it belongs to $\Spa(R_r,R_r^+)$.  We have thus shown that $\MM^{(0),\ad}_{\eta}$ is the union of perfectoid affinoids $\Spa(R_r,R_r^+)$.
\end{proof}

\section{Representation-Theoretic preparations}
\label{RepresentationTheoreticPreparations}

\subsection{Non-abelian Lubin-Tate theory}

Let $\M_{m,\overline{\eta}}^{\ad}$ be the adic geometric generic fiber of the formal scheme $\M_m$, and let
\[ H^i_c=\varinjlim_m H^i_c\left(\M_{m,\overline{\eta}}^{\ad}, \overline{\Q}_{\ell}\right),\]
where $\ell$ is a prime distinct from the residue characteristic of $K$.

Then $H^i_c$ admits an action of $\GL_n(K)\times D^\times\times W_K$, in which elements of the form $(\alpha,\alpha,1)$, $\alpha\in K^\times$, act trivially.  (See \cite{StrauchDeformationSpaces}, \S2.2.2 for a detailed discussion of this action.) {\em Non-abelian Lubin-Tate theory} refers to the realization of Langlands functoriality by the $H^i_c$, as predicted by the conjectures made in \cite{CarayolNonabelianLTTheory}.  Carayol's conjectures have been settled (at least for supercuspidal representations) in \cite{HarrisTaylor} as part of the proof of the local Langlands conjectures for $\GL_n$ over a $p$-adic field.  For a complete historical account of non-abelian Lubin-Tate theory, see the introduction to \cite{StrauchDeformationSpaces}.

For the remainder of this discussion we assume that $n=2$.  In that case, a complete description of $H^1_c$ was given in \cite{Carayol:ladicreps2}, 12.4 Proposition.  See also \cite{CarayolNonabelianLTTheory}, \S3.3.  We present it here because it will be indispensable to the proof of our main theorem.

Let $\Pi\mapsto\LLC(\Pi)$ be the bijection between irreducible admissible representations of $\GL_2(K)$ (with complex coefficients) and two-dimensional Frobenius-semisimple Weil-Deligne representations of $K$ afforded by the local Langlands correspondence.  Write $\mathcal{H}(\Pi)=\LLC(\Pi\otimes\abs{\det}^{-1/2})$;  then $\Pi\mapsto\mathcal{H}(\Pi)$ is compatible under automorphisms of the complex field.  Thus $\Pi\mapsto\mathcal{H}(\Pi)$ may be extended unambiguously to representations with coefficients in any algebraically closed field of characteristic zero, e.g. $\overline{\Q}_\ell$.

Let $\chi\from K^\times\to \overline{\Q}_\ell^\times$ be a character, and let $H^1_c[\chi]$ be the subspace of $H^1_c$ on which the center $K^\times$ of $\GL_2(K)$ acts as $\chi$.

\begin{Theorem}
\label{JLinCohomology}
The representation $H^1_c[\chi]$ decomposes as the direct sum:
\[ H^1_c[\chi]=\bigoplus_{\Pi\in A^2(\chi)} \Pi\otimes \JL(\check{\Pi}) \otimes \mathcal{H}(\Pi)',\]
where $\Pi$ varies through the set $A^2(\chi)$ of discrete series representations of $\GL_2(K)$ with central character $\chi$.  The representation $\mathcal{H}(\Pi)'$ is the unique irreducible quotient of $\mathcal{H}(\Pi)$.
\end{Theorem}
Thus if $\Pi$ is supercuspidal then $\mathcal{H}(\Pi)'=\mathcal{H}(\Pi)$.  The non-supercuspidal discrete series representations of $\GL_2(K)$ are exactly the representations $\St\otimes(\psi\circ\det)$, where $\psi$ is a character of $K^\times$ and $\St$ is the Steinberg representation.  We will write $\St\otimes\psi$ as a shorthand for $\St\otimes(\psi\circ\det)$.  Note that the central character of $\St\otimes\psi$ is $\psi^2$.  For $\Pi$ of this form we have $\dim\mathcal{H}(\Pi)'=1$.

\begin{proof}  We use the notation of \cite{CarayolNonabelianLTTheory}, \S3.3. Carayol shows that if $\mathcal{U}^v$ is the space of vanishing cycles in degree 1 attached to the tower of formal schemes $\MM_m$, then the $\chi$-isotypic component of $\mathcal{U}^v$ decomposes as in the theorem.  By using the comparison theorem of \cite{BerkovichVanishingCyclesII}, Prop. 2.4, $\mathcal{U}^v$ is isomorphic to our $H^1_c$.   See \cite{StrauchDeformationSpaces}, proof of Lemma 2.5.1 for details.
\end{proof}

Let $\MM_{m,\overline{\eta}}^{\ad}/\pi^\Z$ denote the quotient of $\MM_{m,\overline{\eta}}^{\ad}$ by the subgroup of $\GL_n(K)$ generated by the central element $\text{diag}(\pi,\pi)$.  Multiplication by this matrix induces an isomorphism between $\MM_{m,\overline{\eta}}^{(j),\ad}$ and $\MM_{m,\overline{\eta}}^{(j+2),\ad}$, so the quotient $\MM_{m,\overline{\eta}}^{\ad}/\pi^\Z$ is isomorphic to two copies of $\MM_{m,\overline{\eta}}^{{(0)},\ad}$.
As a consequence of Thm. \ref{JLinCohomology} we can give a formula for the dimension of the cohomology of $\MM_{m,\overline{\eta}}^{\ad}/\pi^\Z$.

\begin{Cor} \label{dimJL} We have
\begin{eqnarray*}
\dim H^1_c(\MM_{m,\overline{\eta}}^{\ad}/\pi^\Z,\overline{\Q}_\ell)
&=&
\sum_{\chi_\Pi(\pi)=1} 2\dim \Pi^{\Gamma(\pi^m)}\otimes\JL(\Pi)\\
&&+2q^{m-1}(q-1)(\#\mathbf{P}^1(\OO_K/\pi^m)-1), 
\end{eqnarray*}
where $\Pi$ runs over supercuspidal representations of $\GL_2(K)$ whose central character $\chi_\Pi$ is trivial on $\pi$, and $\Gamma(\pi^m)$ is the congruence subgroup $1+\pi^mM_2(\OO_K)$.
\end{Cor}


\begin{proof}
In Thm. \ref{JLinCohomology}, taking $\Gamma(\pi^m)$-invariants and summing over all $\chi$ with $\chi(\Pi)=1$ yields
\[ H^1_c(\MM_{m,\overline{\eta}}^{\ad}/\pi^\Z,\overline{\Q}_\ell)=\bigoplus_{\chi_{\Pi}(\pi)=1} \Pi^{\Gamma(\pi^m)}\otimes\JL(\Pi)\otimes\mathcal{H}(\Pi)'. \]
For supercuspidal representations $\Pi$ we have $\dim\mathcal{H}(\Pi)'=2$.  All other representations are of the form $\Pi=\St\otimes\psi$, where $\psi$ runs over characters of $K^\times$ with $\psi(\pi)=\pm 1$.  For such $\Pi$, we have that $\dim \JL(\Pi)=\dim \mathcal{H}(\St)'=1$.

It remains to compute $\dim(\St\otimes\psi)^{\Gamma(\pi^m)}$.  Certainly this is only nonzero if $\psi(1+\pi^m\OO_K)=1$.  In that case, $\dim (\St\otimes\psi)^{\Gamma(\pi^m)}=\dim\St^{\Gamma(\pi^m)}$.  Now recall that $\St$ can be modeled on the space of all locally constant functions on $\mathbf{P}^1(K)$, modulo constants.  It follows that $\dim\St^{\Gamma(\pi^m)}=\#\mathbf{P}^1(\OO_K/\pi^m)-1$.   We conclude the proof by noting that the number of characters $\psi$ of $K^\times$ satisfying $\psi(\pi)=\pm 1$ and $\psi(1+\pi^m\OO_K)=1$ is $2\#(\OO_K/\pi^m)=2q^{m-1}(q-1)$. 
\end{proof}

\subsection{Chain orders and strata}
\label{ChainOrders}
We now review the theory of types for $\GL_2(K)$ as presented in \cite{BushnellHenniart}.

A {\em lattice chain} is an $K$-stable family of lattices $\Lambda=\set{L_i}$, with each $L_i\subset K\oplus K$ an $\OO_K$-lattice and $L_{i+1}\subset L_i$ for all $i\in\Z$.  There is a unique integer $e(\Lambda)\in\set{1,2}$ for which $\pi L_i=L_{i+e(\Lambda)}$.  Let $\A_\Lambda$ be the stabilizer in $M_2(K)$ of $\Lambda$.  Up to conjugacy by $\GL_2(K)$ we have
\[\A_{\Lambda} =
\begin{cases}
M_2(\OO_K),& e_{\Lambda}=1,\\
\tbt{\OO_K}{\OO_K}{\gp_K}{\OO_K},& e_{\Lambda}=2
\end{cases}
\]

\begin{defn} A {\em chain order} in $M_2(K)$ is an $\OO_K$-order $\A\subset M_2(K)$ which is equal to $\A_{\Lambda}$ for some lattice chain $\Lambda$.  We say $\A$ is unramified or ramified as $e_{\Lambda}$ is 1 or 2, respectively.
\end{defn}

Suppose $\A$ is a chain order in $M_2(K)$;  let $\mathfrak{P}$ be its Jacobson radical. Then $\mathfrak{P}=\pi\A$ if $\A$ is unramified and $\mathfrak{P}=\tbt{\gp_K}{\OO_K}{\gp_K}{\gp_K}$ in the case that $\A=\tbt{\OO_K}{\OO_K}{\gp_F}{\OO_K}$.  We have a filtration of $\A^\times$ by subgroups $U_{\A}^n=1+\gP^n$, $n\geq 1$.

These constructions have obvious analogues in the quaternion algebra $D$:  If $\A=\OO_D$ is the maximal order in $D$, then the maximal two-sided ideal of $\A$ is generated by a prime element $\pi_D$ of $D$.  We let $U_{\A}^n = 1+\pi_D^n\OO_D$.

\subsection{Characters and Bushnell-Kutzko types}
\label{Characters}

In the following discussion, $\psi\from K\to \overline{\Q}_\ell^\times$ is a fixed additive character.   We assume that $\psi$ is of level one, which means that $\psi(\gp_K)$ is trivial but $\psi(\OO_K)$ is not.  (The choice of level of $\psi$ is essentially arbitrary, but it has become customary to use characters of level one.)  Write $\psi_{M_2(K)}$ for the (additive) character of $M_2(K)$ defined by $\psi_{M_2(K)}(x)=\psi(\tr x)$. Similarly define $\psi_D\from D\to \C$ by $\psi_D(x)=\psi(\tr_{D/K}(x))$, where $\tr_{D/K}$ is the reduced trace.

Now let $A$ be either $M_2(K)$ or $D$.  Let $\A\subset A$ be an $\OO_K$-order which equals a chain order (if $A=M_2(K)$) or the maximal order in $D$ (if $A=D$).
Let $n\geq 1$.  We have a character $\psi_\alpha$ of $U^m_{\A}$ defined by
\begin{eqnarray*}
 U^m_{\A}/U_{\A}^{m+1} &\to& \overline{\Q}_\ell^\times \\
 1+x &\mapsto& \psi(\alpha x)
\end{eqnarray*}
If $\pi$ is an admissible irreducible representation of $\GL_2(K)$, one may ask for which $\alpha$ is the character $\psi$ contained in $\pi\vert_{U^m_{\A}}$.  This is the basis for the classification of representations by Bushnell-Kutzko types, c.f.~\cite{BushnellKutzko}.

\begin{defn}
A {\em stratum} in $A$ is a triple of the form $S=(\A,m,\alpha)$, where $m\geq 1$ and $\alpha\in\gP^{-m}_{\A}$.  Two strata $(\A,m,\alpha)$ and $(\A,m,\alpha')$ are equivalent if $\alpha\equiv\alpha'\pmod{\gP^{1-m}}$.
\end{defn}

\begin{defn}  Let $S=(\A,m,\alpha)$ be a stratum.
\begin{enumerate}
\item $S$ is {\em ramified simple} if $L=K(\alpha)$ is a ramified quadratic extension field of $K$, $m$ is odd, and $\alpha\in L$ has valuation exactly $-m$.
\item $S$ is {\em unramified simple} if $L=K(\alpha)$ is an unramified quadratic extension field of $K$, $\alpha\in L$ has valuation exactly $-m$, and the minimal polynomial of $\pi^m\alpha$ over $K$ is irreducible mod $\pi$.
\item $S$ is {\em simple} if it is ramified simple or unramified simple.
\end{enumerate}
\end{defn}

If $S=(\A,m,\alpha)$ is a stratum in $M_2(K)$ (resp., $D$) and $\Pi$ is an admissible representation of $\GL_2(K)$ (resp., smooth representation of $D^\times$), we say that $\Pi$ {\em contains the stratum $S$} if $\pi\vert_{U^m_{\A}}$ contains the character $\psi_\alpha$.

We call $\Pi$ {\em minimal} if its conductor cannot be lowered by twisting by one-dimensional characters of $F^\times$.

From~\cite{BushnellHenniart} we have the following classification of supercuspidal representations of $\GL_2(K)$:
\begin{Theorem}
\label{classA}
A minimal irreducible admissible representation $\Pi$ of $\GL_2(K)$ is supercuspidal if and only if one of the following holds:
\begin{enumerate}
\item $\Pi$ contains the trivial character of $U^1_{M_2(\OO_K)}$ ({\em i.e.} $\Pi$ has ``depth zero").
\item $\Pi$ contains a simple stratum.
\end{enumerate}
\end{Theorem}

The analogous statement for $D$ is:

\begin{Theorem}
\label{classB}
A minimal irreducible representation $\Pi$ of $D^\times$ of dimension greater than 1 satisfies exactly one of the following properties:
\begin{enumerate}
\item $\Pi$ contains the trivial character of $U^1_{\OO_D}$ ({\em i.e.} $\Pi$ has ``depth zero").
\item $\Pi$ contains a simple stratum.
\end{enumerate}
\end{Theorem}

\section{CM points and linking orders}
\label{CMPointsAndLinkingOrders}

\subsection{CM points}
Once again, let $G_0$ be a one-dimensional formal $\OO_K$-module of height 2 over $\overline{k}$.  Let $L/K$ be a separable quadratic extension, which we consider embedded in a fixed complete algebraically closed field $\C/K$.   In this section we pay special attention to those deformations $(G,\iota)$ of $G_0$ which admit endomorphisms by an order in $L$.  These are investigated in \cite{GrossCanonicalLiftings}.

\begin{defn}  Let $\x$ be an $\OO_{\C}$-point of $\MM_{\infty}=\MM_{G_0,\infty}$ corresponding to the triple $(G,\iota,\phi)$, where $(G,\iota)$ is a deformation of $G_0$ and $\phi\from K^2\to V(G)$ is a basis for its rational Tate module.  Say that $\x$ has {\em CM by $L$} if $G$ admits endomorphisms by an order in $L$.
\end{defn}

Suppose $\x$ is a CM point.  Recall that $D$ is the quaternion algebra of endomorphisms of $G_0$ up to isogeny.  Then $x$ induces embeddings $i_1\from L\to M_n(K)$ and $i_2\from L\to D$, characterized by the commutativity of the diagrams
\[
\xymatrix{
K^2 \ar[d]_{i_1(\alpha)} \ar[r]^{\phi} & V(G) \ar[d]^{\alpha} \\
K^2 \ar[r]_{\phi} & V(G)
}
\]
and
\[
\xymatrix{
G_0 \ar[d]_{i_2(\alpha)} \ar[r]^{\iota} & G \ar[d]^{\alpha} \\
G_0 \ar[r]_{\iota} & G
}
\]
for $\alpha\in K$.  At the risk of minor confusion, from this point forward we will usually suppress $i_1$, $i_2$ from the notation and instead think of $L$ as a subfield of $M_n(K)$ and $D$.  Let $\Delta_{\x}\from L\to M_n(K)\times D$ be the diagonal embedding.  The group $GL_n(K)\times D^\times$ acts transitively on the set of $\C$-points of $\MM_{\infty}$ with CM by $L$;  the stabilizer of $x$ is $\Delta_{\x}(L^\times)$.  We note that after replacing $G$ by an isogenous $\OO_K$-module one can assume that $\End G$ is the maximal order $\OO_L$.

Points in $\MM_\infty(\OO_\C)$ which have CM by some $L/K$ will be called {\em CM points}.  These points give rise to $\C$-points of the perfectoid space $\MM_{\infty,\overline{\eta}}^{\ad}$, which we will also call CM points.

\subsection{Linking Orders}
\label{linkingorders}
To a CM point $\x$ and an integer $m\geq 0$ we will associate an $\OO_K$-order $\L_{\x,m}\subset M_2(K)\times D$ which we have called a ``linking order'' in \cite{WeinsteinFourier}.

The CM point $\x$ corresponds to a triple $(G,\iota,\phi)$, where $(G,\iota)$ is a deformation of $G_0$ to $\OO_\C$ such that $\End G=\OO_L$ and $\phi\from K^2\to V(G)$ is an isomorphism.  The isomorphism $\phi$ allows us to identify $M_2(K)$ with the algebra of $K$-linear endomorphisms of $V(G)\isom K^2$.

The integral Tate module $T(G)$ is a lattice in the $V(G)$.   Since $G$ admits endomorphisms by $L$ up to isogeny, $V(G)$ becomes an $L$-vector space of dimension 1, and it makes sense to talk about the family of lattices $\gp_L^iT(G)\subset V(G)$, for $i\in\Z$.  Observe that these form a lattice chain, cf. \S\ref{ChainOrders}.

\begin{defn}  Let $\A_{\x}\subset M_2(K)$ be the chain order corresponding to the lattice chain $\set{\gp_L^iT(G)}_{i\in\Z}$.
\end{defn}

Up to conjugacy by an element of $\GL_2(K)$ we have
\[  \A_{\x}=\begin{cases}
M_2(\OO_K),&L/K\text{ unramified},\\
\tbt{\OO_K}{\OO_K}{\gp_K}{\OO_K},&L/K\text{ ramified}.\end{cases} \]

Since $T(G)$ is an $\OO_L$-module, $\A_\x$ contains $\OO_L$.  It will be helpful to give a basis for $\A_\x$ as an $\OO_L$-module.  The basis of the $\OO_K$-module $T(G)$ corresponding to $\phi$ takes the form $\set{\alpha_1w,\alpha_2w}$, where $w$ generates $T(G)$ as an $\OO_L$-module and $\alpha_1,\alpha_2$ is a basis for $\OO_L/\OO_K$.   Let $\sigma$ be the nontrivial automorphism of $L/K$.  Define $\varpi_1\in \End_{\OO_K}T(G)$ by $\varpi_1(\alpha_iw)=\alpha_i^\sigma w$, $i=1,2$.  Then we have
\[ \mathcal{A}=\OO_L\oplus\OO_L\varpi_1. \]
In fact this is an orthogonal decomposition with respect to the trace pairing on $M_2(K)$.  Note that $\varpi_1^2=1$ in $\GL_2(K)$.

Now we turn to the corresponding structures in the quaternion algebra $D$.  We have an $L$-linear pairing $D\times D\to L$ given by the reduced trace, which induces an orthogonal decomposition $D=L\oplus C$.  Then $C\cap \OO_D$ is a free $\OO_L$-module of rank 1, generated by an element $\varpi_2$, and then
\[ \OO_D=\OO_L\oplus \OO_L\varpi_2.  \]

Since $\tr_{D/K}(\varpi_2)=0$, we have $\varpi^2_2\in\OO_K$.  Note that $\varpi^2_2$ lies in $\OO_K^\times$ if $L/K$ is ramified, and $\varpi^2_2$ is a uniformizer of $\OO_K$ if $L/K$ is unramified.

We are now ready to define the linking orders $\L_{\x,m}$.  Let $\Delta_{\x}\from L\to M_2(K)\times D$ be the diagonal embedding $\alpha\mapsto (\alpha,\alpha)$.

\begin{defn} Let $m\geq 0$ be an integer.  The {\em linking order} of conductor $m$ associated to the CM point $\mathbf{x}$ is
\[ \L_{\mathbf{x},m}=\Delta_{\mathbf{x}}(\mathcal{\OO}_L)+(\gp_L^m\times\gp_L^m)
+\left(\gp_L^{r_1(m)}\varpi_1\times\gp_L^{r_2(m)}\varpi_2\right),\]
where $r_1(m)=\floor{(m+1)/2}$ and $r_2(m)$ is defined by
\[ r_2(m)=\begin{cases}
\floor{m/2},& L/K\text{ unramified}\\
\floor{(m+1)/2},& L/K\text{ ramified}.
\end{cases}\]

We also define a double-sided ideal $\L_{\mathbf{x},m}^{\circ}\subset \L$ by
\[\L^{\circ}_{\mathbf{x},m}=\Delta_{\mathbf{x}}(\gp_L)+(\gp_L^{m+1}\times\gp_L^{m+1})
+\left(\gp_L^{r_1(m+1)}\varpi_1\times\gp_L^{r_2(m+1)}\varpi_2\right).\]
Let $\mathcal{R}_{\mathbf{x},m}=\L_{\mathbf{x},m}/\L^{\circ}_{\mathbf{x},m}$, a finite-dimensional algebra over $\OO_L/\gp_L$.  Finally, let $\mathcal{R}^1_{\mathbf{x},m}$ be the image of $\L_{\x,m}^\times\cap (\GL_2(K)\times D^\times)^{\det=N}$ in $\mathcal{R}^\times_{\mathbf{x},m}$.
\end{defn}

\begin{rmk} The values of $r_1(m)$ and $r_2(m)$ are the least possible such that the expression given for $\L_{\mathbf{x},m}$ is closed under multiplication.  Also note that if $L/K$ is unramified then $\mathcal{L}_{\x,0}$ is conjugate to $M_2(\OO_K)\times\OO_D$.  Finally, note that $\mathcal{L}_{\x,m}^\times$ is normalized by $\Delta_{\x}(L^\times)$.
\end{rmk}


The following Lemma is \cite{WeinsteinFourier}, Prop. 4.3.4.

\begin{lemma}
\label{RxmValues} The ring $\mathcal{R}_{\mathbf{x},m}$ and the subgroup $\mathcal{R}^1_{\mathbf{x},m}\subset\mathcal{R}^\times_{\mathbf{x},m}$ take the following values.
\begin{enumerate}
\item If $L/K$ is unramified, then $\mathcal{R}_{\mathbf{x},0}\isom M_2(\FF_q)\times\FF_{q^2}$, and $\mathcal{R}_{\mathbf{x},0}^1$ is the subgroup of pairs $(g_1,g_2)$ with $\det g_1=g_2^{q+1}$.
\item If $L/K$ is unramified and $m>0$, then $\L_{\mathbf{x},m}/\L_{\mathbf{x},m}^{\circ}$ is isomorphic to the ring of $3\times 3$ matrices of the form
\[
[\alpha,\beta,\gamma]:=
\begin{pmatrix}
\alpha & \beta & \gamma \\
& \alpha^q & \beta^q \\
& & \alpha
\end{pmatrix}
\]
with $\alpha,\beta,\gamma\in k_L=\OO_L/\gp_L\isom\FF_{q^2}$.  The description of this isomorphism depends on the parity of $m$.

If $m\geq 2$ is even, then a typical element of $\L_{x,m}$ is of the form $\Delta(\alpha)+(\pi^m\gamma,0)+(\pi^{m/2}\beta\varpi_1,\pi^{m/2}\delta\varpi_2)$, with $\alpha,\beta,\gamma,\delta\in\OO_L$.  The isomorphism carries the image of this element in $\R_{\x,m}$ onto the matrix $[\overline{\alpha},\overline{\beta},\overline{\gamma}]$.  (The overline indicates reduction modulo $\gp_L$.)

If $m$ is odd, then a typical element of $\L_{x,m}$ is of the form
$\Delta(\alpha)+(0,\pi^m\gamma)+(\pi^{(m+1)/2}\delta\varpi_1,\pi^{(m-1)/2}\beta\varpi_2)$, with $\alpha,\beta,\gamma,\delta\in\OO_L$.  The isomorphism carries the image of this element in $\R_{\x,m}$ onto the matrix $[\overline{\alpha},\overline{\beta},\overline{\gamma}]$.

In either case, the subgroup $\mathcal{R}_{\mathbf{x},m}^1\subset \mathcal{R}_{\x,m}^\times$ corresponds to the group of matrices $[\alpha,\beta,\gamma]$ with $\alpha\gamma^q+\alpha^q\gamma=\beta^{q+1}$.

\item If $L/K$ is ramified and $m$ is odd, then we have an isomorphism $\mathcal{R}_{\x,m}\isom k[e]/e^2$.  A typical element of $\L_{x,m}$ is of the form $\Delta(\alpha)+(\beta\pi_L^m,0)+(\pi^{(m+1)/2}\gamma\varpi_1,\pi^{(m+1)/2}\delta\varepsilon)$, with $\alpha,\beta,\gamma,\delta\in\OO_L$. The isomorphism carries the image of this element in $\R_{\x,m}$ onto $\overline{\alpha}+\overline{\beta}e$.  We have $\mathcal{R}_{\mathbf{x},m}^1=\mathcal{R}_{\mathbf{x},m}^\times.$
\end{enumerate}
\end{lemma}

We will not be needing an explicit presentation for $\R_{\x,m}$ in the case that $L/K$ is ramified and $m$ is even.

\begin{defn} Define groups $\K_{\x,m}$ and $\K^1_{\x,m}$ by
\begin{eqnarray*}
\K_{\x,m}&=&\Delta_{\x}(L^\times)\L^\times_{\x,m}\\
\K_{\x,m}^1&=&\K_{\x,m}\cap (\GL_2(K)\times D^\times)^{\det=N}
\end{eqnarray*}
\end{defn}

We now construct a family of representations of the groups $\K_{\x,m}$, as in Thm. 5.0.3 of \cite{WeinsteinFourier}.

\begin{defn} \label{RhoSDefinition} Let $x$ be a CM point, and Let $S$ be a simple stratum of the form $(\A_x,m,\alpha)$.  We define a certain irreducible representation $\rho_S$ of $\mathcal{K}_{\x,m}$ with coefficients in $\overline{\Q}_\ell$.  The representation $\rho_S$ will have the following properties:
\begin{enumerate}
\item The restriction of $\rho_S$ to $\L_{\x,m}^\times$ factors through $\R_{\x,m}^\times$.
\item The restriction of $\rho_S$ to $\Delta(K^\times)\subset \K_{\x,m}$ is a sum of copies of the trivial representation.
\item The restriction of $\rho_S$ to $U_{A_{\x}}^m\times U_{\OO_D}^m$ is sum of copies of the character $(1+x,1+y)\mapsto\psi_{M_2(K)}(\alpha x)\psi_{D}(\alpha x)^{-1}$.
\end{enumerate}

If $L/K$ is ramified, then we have an isomorphism $\R_{\x,m}^\times\isom (k[e]/e^2)^\times$ under which $U^m_{\A_x}\times U^m_{\OO_D}$ corresponds to $1+ke$.  The stratum $S$ determines a nontrivial character of this group, which we extend to all of $\R_{\x,m}^\times$ by requiring that it be trivial on $k^\times$.  Inflate this to get a character $\rho_S$ of $\L_{\x,m}^\times$.  Finally, extend $\rho_S$ to all of $\K_{\x,m}=\Delta_\x(L^\times)\L_{\x,m}^\times$ by declaring $\rho_S(\Delta_{\x}(\pi_L))=-1$.

If $L/K$ is unramified and $m\geq 1$, then $\R_{\x,m}^\times$ is isomorphic to the group of matrices of the form
\[ [\alpha,\beta,\gamma]=\begin{pmatrix}
\alpha & \beta & \gamma \\
& \alpha^q & \beta^q \\
& & \alpha
\end{pmatrix}, \;\alpha\in k_L^\times,\;\beta,\gamma\in k_L.
\]
Under this isomorphism, the image of $U^m_{\A_x}\times U^m_{\OO_D}$ in $\R_{\x,m}^\times$ corresponds to the subgroup $U=\set{[0,0,\gamma]\vert \gamma\in k_L}$, which lies in the center of $\R_{\x,m}
$.
The most direct way to construct $\rho_S$ uses the $\ell$-adic cohomology of a curve admitting an action by $\R_{\x,m}^\times$.  Recall that $\R_{\x,m}^1$ is isomorphic to the subgroup of $\GL_3(k_L)$ consisting of matrices of the form $[\alpha,\beta,\gamma]$ with $\beta^{q+1}=\alpha^q\gamma+\alpha\gamma^q$.  Observe that this group preserves the affine plane curve $C^1$ defined by the equation $Z_1^q+Z_1=Z_2^{q+1}$, under the action $(Z_1,Z_2)\mapsto (Z_1+\alpha^{-1}\beta^q Z_2+\alpha^{-1}\gamma,\alpha^{q-1} Z_2+\alpha^{-1}\beta)$.  Let $C=\R_{\x,m} \times_{\R_{\x,m}^1} C^1$.  By inflation, we get an action of $\L_{\x,m}^\times$ on $C$.  Extend this action to all of $\K_{\x,m}$ by having $\Delta_{\x}(\pi)$ act trivially.


Let $\rho$ be the action of $\K_{\x,m}$ on $H^1_c(C,\overline{\Q}_\ell)$.  Since $U$ lies in the center of $\R_{\x,m}^{\times}$, and $\psi$ is a character of $U$, it makes sense to define $\rho_S$ as the $\psi$-isotypic component of $\rho$.  It is irreducible of dimension $q$.  For proofs of these claims, see \cite{WeinsteinFourier}, \S5.1.

Finally, we turn to the case $L/K$ unramified and $m=0$.  We have $\R_{\x,0}=M_2(\OO_K)\times k_L$.  Let $\theta$ be a character of $k_L^\times$ which does not factor through the norm map $k_L^\times\to k^\times$.  There exists an irreducible cuspidal representation $\eta_\theta$ of $\GL_2(k)$ corresponding to $\theta$.  The character of this representation takes the value $-(\theta(\alpha)+\theta(\alpha^q))$ on an element $g\in \GL_2(k)$ with distinct eigenvalues $\alpha,\alpha^q\in k_E$ not lying in $k$.  Let $\rho_\theta$ be the character $\eta_\theta\otimes\theta^{-1}$ of $\R_{\x,0}^\times=\GL_2(k)\times k_L^\times$.  Extend $\rho_\theta$ to all of $\K_{x,0}$ by having $\Delta_{\x}(\pi)$ act trivially.
\end{defn}

\begin{Theorem} \label{MainThmOfWei10} Let $\x$ be a point with CM by $L$.
\begin{enumerate}
\item Let $S$ be a simple stratum of the form $(\A_{\x},m,\alpha)$, with $\alpha\in L$.  Let $\rho_S$ be the representation of $\K_{\x,m}$ described in Defn. \ref{RhoSDefinition}.  Let
\[ \Pi_S=\Ind_{\K_{\x,m}}^{\GL_2(K)\times D^\times} \rho_S \]
(compactly supported smooth induction).  For a character $\chi$ of $K^\times$, let $\Pi_S[\chi]$ be the subspace of $\Pi_S$ on which the center of $\GL_2(K)$ acts through $\chi$.  Then $\Pi_S[\chi]$ is the direct sum of representations of $\GL_2(K)\times D^\times$ of the form $\Pi\otimes\JL(\check{\Pi})$, where $\Pi$ is a minimal supercuspidal representation having central character $\chi$ and containing the stratum $S$.  Every such $\Pi$ appears in $\Pi_S[\chi]$.

\item Now suppose $L/K$ is unramified, and let $\theta$ be a character of $k_L^\times$ which does not factor through the norm map $k_L^\times\to k^\times$.  Let
\[ \Pi_\theta=\Ind_{\K_{\x,0}}^{\GL_2(K)\times D^\times} \rho_\theta, \]
and define $\Pi_\theta[\chi]$ as above.  Then $\Pi_\theta[\chi]$ is the direct sum of representations of $\GL_2(K)\times D^\times$ of the form $\Pi\otimes JL(\check{\Pi})$, where $\Pi$ is a minimal supercuspidal representation of depth 0 and central character $\chi$.  Every such $\Pi$ appears in some $\Pi_\theta[\chi]$.
\end{enumerate}
\end{Theorem}

\begin{proof}  This is a restatement of Thm. 6.0.1 of \cite{WeinsteinFourier}.  There, the center is treated a little differently.   In \cite{WeinsteinFourier} one first defines a representation $\rho_S$ of $\L_{\x,m}^\times$, and then (after choosing a central character $\chi$) extends of $\rho_S$ to a representation $\rho_{S,\chi}$ of $(K^\times\times K^\times)\L_{x,m}^\times=(K^\times\times\set{1})\K_{\x,m}^\times$ by having $K^\times\times\set{1}$ act through $\chi$.   Thm. 6.0.1 of \cite{WeinsteinFourier} says that the representation of $\GL_2(K)\times D^\times$ induced from $\rho_{S,\chi}$ is a direct sum of representations $\Pi\otimes\JL{\check{\Pi}}$ as claimed in Thm. \ref{MainThmOfWei10}.  But it is easy to see that the induced representation of $\rho_{S,\chi}$ is the same thing as our $\Pi_S[\chi]$ above.  The argument for depth zero supercuspidals is similar.
\end{proof}

\subsection{Curves over $\overline{\FF}_q$ and the Jacquet-Langlands correspondence}
\begin{defn} \label{DefnOfCxM} We define a smooth affine curve $C_{\mathbf{x},m}$ over $\overline{\FF}_q$.  In each case we will define an action of $\mathcal{K}_{\mathbf{x},m}^1$ on $C_{\mathbf{x},m}$, whose restriction to $\L_{\mathbf{x},m}^\times$ factors through the finite group $\mathcal{R}_{\mathbf{x},m}^\times$.
\begin{enumerate}
\item When $L/K$ is unramified and $m=0$, let $C_{\mathbf{x},0}$ be the affine plane curve with equation $Z_1Z_2^q-Z_1^qZ_2=1$.  An element of $\mathcal{R}_{\x,0}^1$ corresponds to a pair $(g_1,g_2)$ in $\GL_2(\FF_q)\times\FF_{q^2}^\times$ satisfying $\det g_1=g_2^{q+1}$.  Suppose $g_1=\tbt{a}{b}{c}{d}$.  This pair will act on $C_{\mathbf{x},0}$ by sending $(Z_1,Z_2)$ to $g_2^{-1}(aZ_1+cZ_2,bZ_1+dZ_2)$.

\item When $L/K$ is unramified and $m>0$ is even, let $C_{\mathbf{x},m}$ be the affine plane curve with equation $Z_1^q+Z_1=Z_2^{q+1}$.  An element of $\mathcal{R}_{\mathbf{x},m}^1$ corresponds to a matrix $[\alpha,\beta,\gamma]\in\GL_3(\FF_{q^2})$ which satisfies $\alpha \gamma^q+\alpha^q\gamma=\beta^{q+1}$.  This element will act on $C_{\mathbf{x},m}$ by sending $(Z_1,Z_2)$ to $\alpha^{-1}(\alpha Z_1+\beta^q Z_2+\gamma,\alpha^q Z_2+\beta)$.

\item When $L/K$ is unramified and $m>0$ is odd, let $C_{\mathbf{x},m}$ be the affine plane curve with equation $Z_1^q+Z_1=Z_2^{q+1}$.  The action of $\K_{\x,m}$ is defined as follows:  $[\alpha,1,1]$ acts as $(Z_1,Z_2)\mapsto (Z_1,\alpha^{q-1}Z_2)$ (as it did in the $m$ even case), but $[1,\beta,\gamma]^{-1}$ acts as $(Z_1,Z_2)\mapsto (Z_1+\beta^q Z_2+\gamma^q,Z_2+\beta)$.

\item When $L/K$ is ramified and $m$ is odd, let $C_{\mathbf{x},m}$ be the affine plane curve with affine equation $Z_1^q-Z_1=Z_2^2$.  An element of $\mathcal{R}_{\mathbf{x},m}^1$ corresponds to an element $a+be\in \FF_q[e]/e^2$, $a\neq 0$.  This element will act on $C_{\mathbf{x},m}$ by sending $(Z_1,Z_2)$ to $(Z_1+a^{-1}b,Z_2)$.
\end{enumerate}
In each case $C_{\x,m}$ has an action of $\mathcal{R}_{\x,m}^\times$, hence of $\L_{\x,m}^{\times}$.  Extend this to an action of $\mathcal{K}_{\x,m}$ by the following rule:  if $L/K$ is unramified, have $\Delta_{\x}(\pi)$ act trivially, and if $L/K$ is ramified, with uniformizer $\pi_L$, have $\Delta_{\x}(\pi_L)$ send $(Z_1,Z_2)$ to $(Z_1,-Z_2)$.
\end{defn}

The compactly supported cohomology $H^1_c(C_{\x,m},\Q_\ell)$ is a smooth representation of $\mathcal{K}_{\mathbf{x},m}$.

\begin{prop} \label{H1OfC} Suppose $m>0$.  If $S=(\A_{\x},m,\alpha)$ is a simple stratum then $\rho_S\vert_{\K^1_{\x,m}}$ is a direct summand of $H^1_c(C_{\x,m},\overline{\Q}_\ell)$.

Now suppose that $L/K$ is unramified and $m=0$.  If $\theta$ is character of $k_L^\times$ which does not factor through the norm map $k_L^\times\to k^\times$, then $\rho_\theta\vert_{\K^1_{\x,0}}$ is a direct summand of $H^1_c(C_{\x,0},\overline{\Q}_\ell)$.
\end{prop}



\begin{proof}  The proof proceeds by cases.

\begin{enumerate}
\item When $L/K$ is ramified (so that $m$ is odd), $C_{\x,m}$ is the curve $Z_1^q-Z_1=Z_2^2$, and $\mathcal{R}_{\x,m}^\times$ acts through a quotient $\FF_q$ acting through substitutions $(Z_1,Z_2)\mapsto (Z_1+a,Z_2)$, with $a\in\FF_q$.  The claim is reduced to showing that $H^1_c(C_{\x,0},\overline{\Q}_\ell)$ is the direct sum of the nontrivial characters of $\FF_q$.  This is an exercise, see \cite{KatzCrystallineCohomology}, Cor. 2.2, which proves a more general statement about curves of the form $Z_1^q-Z_1=Z_2^N$:  the isotypic component of a nontrivial character in $H^1_c$ has dimension $N-1$.

\item When $L/K$ is unramified and $m$ is even, $C_{\x,m}$ is the very curve $C^1$ used to define the representations $\rho_S$ in the first place:  the $\rho_S$ were defined as summands of $\Ind_{\R_{\x,m}^1}^{\R_{\x,m}} H^1_c(C_{\x,m},\overline{\Q}_\ell)$, inflated to $\L_{\x,m}$ and extended to $\K_{\x,m}$ by having $\Delta_{\x}(\pi)$ act trivially.  From this description it is clear that $\rho_S\vert_{\K_{\x,m}^1}$ is a direct summand of $H^1_c(C_{\x,m},\overline{\Q}_\ell)$.

\item When $L/K$ is unramified and $m$ is odd, $C_{\x,m}$ carries a different action of $\R_{\x,m}$ as the curve $C^1$ used to define the $\rho_S$, but the same representations of $\R_{\x,m}^\times$ appear in the cohomology $H^1_c(C_{\x,0},\overline{\Q}_\ell)$.  To prove this, it is enough to show that the trace of an element of $\R_{\x,m}^{\times}$ acting on the Euler characteristic of $C^1$ is the same in either action, and this is easily done using the Lefschetz fixed-point formula.

\item When $L/K$ is unramified and $m=0$, $C_{\x,0}$ is one connected component of the Deligne-Lusztig curve associated to $\GL_2(k)$, which is $C_{\text{DL}}=\R_{\x,0}^\times\times_{\R_{x,0}^1} C_{\x,0}$.  The observation that $H^1(C_{\text{DL}},\overline{\Q}_\ell)$ is the direct sum of all representations of the form $\eta_{\theta}\otimes\theta^{-1}$, with $\theta\neq\theta^q$, goes back to Drinfeld.  It follows that if $\eta_{\theta}\otimes\theta^{-1}$ is any such representation of $\R_{\x,0}^\times$, then $\eta_{\theta}\otimes\theta^{-1}\vert_{\R_{\x,0}^1}$ is a direct summand of $H^1(C_{\x,0},\overline{\Q}_\ell)$ as a representation of $\R_{\x,0}^1$.  This shows that $\rho_\theta\vert_{K_{\x,0}^1}$ is a direct summand of $H^1(C_{\x,0},\overline{\Q}_\ell)$ as a representation of $\K_{\x,0}^1$.
\end{enumerate}
\end{proof}




\begin{Theorem}  \label{CurvesRealizeJLC}  Let $\x$ be a point with CM by $L$.  Let $m\geq 0$, assumed odd if $L/K$ is ramified.  The representation
\[ V_{\x,m}=\Ind_{\K_{\x,m}^1}^{(\GL_2(K)\times D^\times)/(\pi,1)^{\Z}} H^1_c(C_{\x,m},\overline{\Q}_\ell) \]
(compactly supported smooth induction) contains
\[ \bigoplus_{\Pi}\Pi\oplus\JL(\check{\Pi})^{\oplus 2}, \]
where $\Pi$ runs over supercuspidal representations of $\GL_2(K)$ whose central character is trivial on $\pi$, and which have the property that
\begin{enumerate}
\item if $m\geq 1$, then some twist of $\Pi$ contains a simple stratum of the form $(\A_x,m,\alpha)$.
\item if $m=0$, then some twist of $\Pi$ has depth zero.
\end{enumerate}
\end{Theorem}

\begin{proof}  Let $\Pi$ be a supercuspidal representation of $\GL_2(K)$ whose central character is trivial on $\pi$.  Suppose that some twist of $\Pi$, call it $\Pi'$, contains a simple stratum $S$ of the form $(\A_x,m,\alpha)$.  By Thm. \ref{MainThmOfWei10}, $\Pi'\otimes\JL(\check{\Pi}')\vert_{\K_{\x,m}}$ contains $\rho_S$.  By Prop. \ref{H1OfC}, $\rho_S\vert_{\K^1_{\x,m}}$ is a summand of $H^1_c(C_{\x,m},\overline{\Q}_\ell)$.  This shows that $\Pi\otimes\JL(\check{\Pi})\vert_{\K^1_{\x,m}}=\Pi'\otimes\JL(\check{\Pi'})\vert_{\K^1_{\x,m}}$ contains a summand of $H^1_c(C_{\x,m},\overline{\Q}_\ell)$.  Therefore by Frobenius reciprocity, $\Pi\otimes\JL(\check{\Pi})\vert_{\K^1_{\x,m}}$ is contained in $V_{\x,m}$.  Now we observe that if $\Pi'$ contains $S=(\A_x,m,\alpha)$, then it also contains the distinct stratum $S^\sigma=(\A_x,m,\alpha^{\sigma})$, where $\alpha^{\sigma}\in L=K(\alpha)$ is $K$-conjugate to $\alpha$.  By repeating the argument we find that $\Pi\otimes\JL(\check{\Pi}')$ is contained in $V_{\x,m}$ with multiplicity at least 2.

The argument for depth zero supercuspidals is similar.
\end{proof}

\section{Special affinoids in the Lubin-Tate tower}
\label{specialaffinoids}

\subsection{Special affinoids: an overview}
Let $\Pi$ be a supercuspidal representation of $\GL_2(K)$ with coefficients in $\overline{\Q}_\ell$.  Thm. \ref{JLinCohomology} shows that $\Pi\otimes\JL(\check{\Pi})$ appears in $H^1_c$ of the Lubin-Tate tower.  On the other hand, Thm. \ref{CurvesRealizeJLC} shows that there exists a nonsingular affine curve $C_{\x,m}$ over $\overline{k}$ admitting an action of a subgroup $\mathcal{K}^1_{\x,m}\subset (\GL_2(K)\times D^\times)^{\det=N}$, such that $\Pi\otimes\JL(\check{\Pi})$ is contained in the induced representation of $H^1_c(C_{\x,m},\overline{\Q}_\ell)$.  (Here $\x$ is a CM point and $m\geq 0$ is an integer, both of which depend on $\Pi$.)  The next theorem finds a link between the Lubin-Tate tower and the curve $C_{\x,m}$.  Essentially, the presence of $\Pi\otimes\JL(\check{\Pi})$ in the cohohomology of the Lubin-Tate tower can be traced to the existence of an open affinoid subset of $\MM_{\infty,\overline{\eta}}^{\ad}$ whose reduction is related to $C_{\x,m}$.

First we must define what we mean by the reduction of an affinoid.

\begin{defn} Let $Z=\Spa(R,R^+)$ be an affinoid adic space over $\Spa(\C,\OO_\C)$.  The {\em reduction} of $Z$ is $\overline{Z}=\Spec R^+/\gm_\C$, where $\gm_\C\subset\OO_{\C}$ is the maximal ideal.  It is a reduced affine scheme over $\overline{k}=\OO_\C/\gm_\C$.
\end{defn}

Note that if $Z$ is a perfectoid affinoid, then $\overline{Z}$ is the spectrum of a perfect $\overline{k}$-algebra.

To state the theorem precisely, it is convenient to work with one connected component $\MM_{\infty,\overline{\eta}}^{\circ,\ad}$ of $\MM_{\infty,\overline{\eta}}^{\ad}$.  Recall that $\MM_{\infty,\overline{\eta}}^{\circ,\ad}$ admits an action of $(\GL_n(K)\times D^\times)^{\det=N}$.

\begin{Theorem}  \label{ExistenceOfAffinoid} Assume that the residue characteristic of $K$ is odd.  Let $\x\in \MM_{\infty,\overline{\eta}}^{\circ,\ad}(\OO_\C)$ be a point with CM by a quadratic field $L/K$.  Let $m\geq 0$ be an integer, assumed to be odd if $L/K$ is ramified.  Then there exists an open affinoid subset $\ZZ_{\x,m}\subset \MM_{\infty,\overline{\eta}}^{\circ,\ad}$ with the following propoerties:
\begin{enumerate}
\item $\ZZ_{\x,m}$ is stabilized by the action of $\K_{\x,m}^1$.
\item There exists a nonconstant morphism $\overline{\ZZ}_{\x,m}\to C_{\x,m}$ which is equivariant for the action of $\K_{\x,m}^1$.
\item For $x$ fixed, the $\ZZ_{\x,m}$ form a decreasing sequence of open neighborhoods of $\x$, and $\cap_m \ZZ_{\x,m}=\set{\x}$.
\item For $g\in (\GL_2(K)\times D^\times)^{\det=N}$ we have $\ZZ_{\x,m}^g=\ZZ_{\x^g,m}$.
\end{enumerate}
\end{Theorem}

\begin{rmk}  If $X=\Spec R$ is a reduced affine scheme over $\overline{k}$, we write $X^{\perf}$ for the scheme $\Spec R^{\perf}$, where $R^{\perf}=\varprojlim_{x\mapsto x^p} R$ is the perfect closure of $R$.  Since $\overline{\ZZ}_{\x,m}$ is the spectrum of a perfect ring, $\overline{\ZZ}_{\x,m}\to C_{\x,m}$ factors through a morphism $\overline{\ZZ}_{\x,m}\to C_{\x,m}^{\perf}$.  In fact, $\overline{\ZZ}_{\x,m}\to C_{\x,m}^{\perf}$ is an isomorphism, but we can only prove this {\em a posteriori}.  We will also define an affinoid $\ZZ_{\x,m}$ when $L/K$ is ramified and $m$ is even, but we will not need to analyze it as intensely.  We will show that $\ZZ_{\x,m}$ is the inverse limit of curves whose completion has genus 0, also {\em a posteriori}.
\end{rmk}

The proof of Thm. \ref{ExistenceOfAffinoid} is a case-by-case calculation which we will undertake in the following sections.  This calculation fits the following pattern.  Recall that $\MM_{\infty,\overline{\eta}}^{\ad}$ is isomorphic to a subspace of $\tilde{G}^{\ad}_{\overline{\eta}}\times\tilde{G}^{\ad}_{\overline{\eta}}$.  We will first construct an open affinoid $\YY_{\x,m}\subset \tilde{G}^{\ad}_{\overline{\eta}}\times\tilde{G}^{\ad}_{\overline{\eta}}$, and then define $\ZZ_{\x,m}=\YY_{\x,m}\cap\MM^{\circ,\ad}_{\overline{\eta}}$.  We will show that the group $\K_{\x,m}$ stabilizes $\YY_{\x,m}$, and therefore that $\K_{\x,m}^1$ stabilizes $\ZZ_{\x,m}$.  For its part, $\YY_{\x,m}$ may be described as a ``rectangle" centered around $\x$:  it is given by a pair of bounds on two linear forms which vanish on $\x$.  As $m\to\infty$, the bounds decrease to 0, so that $\cap_m \YY_{\x,m}=\set{\x}$.

The affinoids $\YY_{\x,m}$ will be defined in \S\ref{DefinitionOfAffinoids}.  In the meantime we need to develop some language for the geometry of the formal vector space $\tilde{G}$.

\subsection{Some convenient coordinate systems}
In our calculations it will be helpful to fix coordinates on the formal $\OO_K$-modules $G$ and $\wedge^2 G$.  By a coordinate on $G$ mean an isomorphism $G\isom \Spf\OO_K\powerseries{T}$, or equivalently a collection of power series (a ``law'') defining a formal $\OO_K$-module isomorphic to $G$.

For $G$ we take the formal $\OO_K$-module law whose logarithm is
\[ \log_G(T) = T+\frac{T^{q^2}}{\pi}+\frac{T^{q^4}}{\pi^2}+\dots. \]
This is the series obtained by setting $v_1=0$, $v_2=1$ and $v_n=0$ for $n\geq 3$ in Hazewinkel's functional equation, cf. \S\ref{ModuliOfFormalModules}.

The formal $\OO_K$-module $\wedge^2 G$ is the (unique up to isomorphism) formal $\OO_K$-module over $\OO_{K_0}$ of dimension 1 and height 1.  Let us fix $\wedge^2 G_0=\Spf\OO_K\powerseries{T}$ by specifying its logarithm as
\[ \log_{\wedge^2 G}(T) = T-\frac{T^q}{\pi}+\frac{T^{q^2}}{\pi^2}-\frac{T^{q^3}}{\pi^3} +\dots.\]

\begin{lemma}
\label{AdditionLawG}
\begin{enumerate}
\item $[\pi]_{G}(T)\equiv T^{q^2}$ modulo $(\pi,T^{q^2+1})$.
\item $[\pi]_{\wedge^2 G}(T)\equiv -T^q$ modulo $(\pi,T^{q+1})$.
\item $X_1+_{G} X_2\equiv X_1+X_2$ modulo terms of total degree $q^2$.
\item $X_1+_{\wedge^2 G}X_2\equiv X_1+X_2$ modulo terms of total degree $q$.
\end{enumerate}
\end{lemma}

\begin{proof}  The first two claims are special cases of the congruences in Eq. \eqref{picongs}.  For the third, note that $\log_G(X_1+_G X_2)=\log_G (X_1)+\log_G (X_2)$ and use the fact that $\log_G(T)\equiv T$ modulo terms of degree $q^2$.  The fourth claim is similar.
\end{proof}

Now we turn to formal vector spaces.  Our choice of coordinate $G\isom\Spf \OO_{K_0}\powerseries{T}$ gives rise to a coordinate $\tilde{G}\isom \Spf\OO_{K_0}\powerseries{T^{1/q^\infty}}$, along the lines of Rmk. \ref{boldfaceX}.  Suppose $R\in\Alg_{\OO_{K_0}}$ and $\X\in \tilde{G}(R)$.  Suppose $\X=(X_0,X_1,\dots)$ with $X_i\in \Nil(R)$ satisfying $[\pi]_G(X_{i+1})=X_i$, $i\geq 0$.  Let
\[  X=\lim_{i\to\infty} X_i^{q^i}. \]
We call $X$ the {\em coordinate} of $\X$, and we write
\[ X=[\X] \]
Note that $X$ comes equipped with a privileged root $X^{1/q^i}$ for every $i\geq 1$.

Recall that $D=\End G_0\otimes K$ acts on $\tilde{G}$.  It will be helpful to give an approximation for this action in terms of coordinates.  Let $R$ be a Banach $\C$-algebra, with norm $\abs{\;}$, and let $R^+\subset R$ be the $\OO_{\C}$-subalgebra of elements $f$ with $\abs{f}\leq 1$.  Assume that $R^+$ is bounded (and thus complete for the $\pi$-adic topology).  We say that $X_1,X_2\in R$ are {\em equal up to smaller terms} if $\abs{X_1-X_2}<\abs{X_1}=\abs{X_2}$.  Finally, let $\pi_D\in D$ be the Frobenius element $T\mapsto T^q$, which is a uniformizer for $D$.

\begin{lemma}
\label{ApproximateActionOnTildeG}
\begin{enumerate}
\item Let $\X,\Y\in\tilde{G}(R^+)$.  Then $[\X+\Y]=[\X]+[\Y]$ up to smaller terms.
 Let $g\in D^\times$.  Suppose that $g=u\pi_D^m$, with $u\in\OO_D^\times$, $m\in\Z$.  Let $\overline{u}$ be the image of $u$ in $\OO_D/\pi_D\injects\overline{k}$.  Then $[g\X]=\overline{u}[\X]^{q^{m}}$ up to smaller terms.   In particular $[\pi\X]=X^{q^2}$ plus smaller terms.
\item Similarly, suppose $\X,\Y\in \widetilde{\wedge^2 G}(R)$.  Then $[\X+\Y]=[\X]+[\Y]$ up to smaller terms.  Let $g\in K^\times$, with $g=u\pi^m$, $u\in\OO_K^\times$, $m\in\Z$.  Then $[g\X]=(-1)^m\overline{u}[\X]^{q^m}$ up to smaller terms.
\end{enumerate}
\end{lemma}
\begin{proof} These statements follow easily from the corresponding statements about the formal group laws $G$ and $\wedge^2 G$.  For instance, if $\X,\Y\in\tilde{G}(R^+)$ with $X=[\X]$, $Y=[\Y]$, then
\[ [\X+\Y]=\lim_{i\to \infty} (X^{1/q^i}+_G Y^{1/q^i})^{q^i}. \]
Thus for some $i$ large enough, $[\X+\Y]=(X^{1/q^i}+_G Y^{1/q^i})^{q^i}$ up to smaller terms.  Since $+_G=+$ up to quadratic terms, and because exponentiation by $q^i$ commutes with addition up to smaller terms, we get $[\X+\Y]=X+Y$ up to smaller terms.
\end{proof}

\subsection{An approximation for the determinant morphism}
Recall from \S\ref{detmodules} that we have a determinant map $\lambda\from \tilde{G}\times\tilde{G}\to\widetilde{\wedge^2 G}$.  This map corresponds to a continuous homomorphism from $\OO_{K_0}\powerseries{T^{1/q^\infty}}$ into $\OO_{K_0}\powerseries{X_1^{1/q^\infty},X_2^{1/q^\infty}}$.  Let $\delta(X_1,X_2)$ be the image of $T$ under this homomorphism, and let $\delta(X_1,X_2)^{1/q^m}$ be the image of $T^{1/q^m}$, for $m\geq 1$.  Let $\delta_0(X_1,X_2)$ be the image of $\delta(X_1,X_2)$ in $\overline{k}\powerseries{X_1^{1/q^\infty},X_2^{1/q^\infty}}$.

\begin{prop} \label{deltaestimate}   Possibly after replacing $\delta_0$ with $[\alpha]_{\widetilde{\wedge^2 G}}(\delta_0)$ for some $\alpha\in K^\times$, the congruence
\[ \delta_0(X_1,X_2)\equiv X_1X_2^q-X_1^qX_2 \]
holds modulo terms of total degree $q^2$.
\end{prop}

\begin{proof}   Let $d$ be the least degree of any term appearing in $\delta_0(X_1,X_2)$, and let $F$ be the homogeneous part of $d$ of degree $d$.  We have that $\delta$ is $\OO_K$-alternating with respect to the operations $+_{G_0}$ and $[g]_{G_0}$ ($g\in\OO_K$), and similarly for $\wedge^2G_0$.  These operations are simply addition and scalar multiplication modulo quadratic terms.  Thus $F$ is a $k$-bilinear alternating form, which is to say it is of the form
\[ F = \sum_{(a_1,a_2)} c_{a_1,a_2} X_1^{q^{a_1}}X_2^{q^{a_2}} \]
where $(a_1,a_2)$ runs over pairs of integers with $q^{a_1}+q^{a_2}=d$, and $c_{a_1,a_2}\in \overline{k}$ satisfies $c_{a_2,a_1}=-c_{a_1,a_2}$.  After replacing $\delta_0$ with $[\pi^m]_{\widetilde{\wedge^2 G_0}}(\delta_0)$ for some $m$, we may assume that $F$ contains a nonzero term of the form $c_{0,a_2}X_1X_2^{q^{a_2}}$, with $a_2\geq 1$.

Since
\[ \delta_0([\pi]_{G_0}(X_1),X_2)) = [\pi]_{\wedge^2 G_0}(\delta_0(X_1,X_2)), \]
we find (using Lemma \ref{AdditionLawG}) that $-\delta_0^q$ contains the term $c_{0,a_2}X_1^{q^2}X_2^{q^{a_2}}$, which shows that $\delta_0$ contains the term $-c_{0,a_2}^{1/q}X_1^qX_2^{q^{a_2-1}}$.  By definition of $d$ we have $d\leq q+q^{a_2-1}$, but on the other hand $d=1+q^{a_2}$, which shows that $a_2=1$ and $d=q+1$.

The only integral solutions to $q^{a_1}+q^{a_2}=q+1$ are $(1,0)$ and $(0,1)$.  Thus
\[ F=c_{0,1}(X_1X_2^q-X_1^qX_2). \]
By the above observation, $\delta_0$ contains the term $-c_{0,1}^{1/q}X_1^qX_2$, which shows that $c_{0,1}=c_{0,1}^q$, thus $c_{0,1}\in k\isom\FF_q$.  After replacing $\delta_0$ with some multiple $[\alpha]_{\wedge^2 G_0}(\delta_0)$, we can assume that $c_{0,1}=1$, so that $F=X_1X_2^q-X_1^qX_2$.

Now consider the difference $E=\delta-F$.  By Lemma \ref{AdditionLawG}, the addition and scalar multiplication laws in $G_0$ and $\wedge G_0$ equal ordinary addition and scalar multiplication up to degree $q^2$.  This shows that $E$ is also $\OO_K$-bilinear and alternating modulo degree $q^2$.  Suppose for the sake of contradiction that the leading homogeneous part of $E$, call it $F'$, has degree $d'<q^2$.   Then $F'$ is $k$-bilinear and alternating.  The foregoing argument shows that $[\alpha]_{\wedge^2 G_0}(F')=X_1X_2^q-X_1^qX_2$ for some $\alpha\in K^\times$.  If $\alpha\in \OO_K$, this contradicts $d'>d=q+1$, and if $\alpha\not\in\OO_K$, this contradicts $d<q^2$.  Thus the degree of $E$ is at least $q^2$.
\end{proof}

We need to translate Prop. \ref{deltaestimate} into an approximation for $\delta$.  We have shown that $\delta \equiv F+E$ modulo $\pi$, where $F(X_1,X_2)=X_1X_2^q-X_1^qX_2$ and $E(X_1,X_2)\in\OO_{K_0}\powerseries{X_1^{1/q^\infty},Y_1^{1/q^\infty}}$ has degree $\geq q^2$.  Then
\begin{equation}
\label{DeltaLimit}
\delta=\lim_{m\to\infty} \left(F(X_1^{1/q^m},X_2^{1/q^m})+E(X_1^{1/q^m},X_2^{1/q^m})\right)^{q^m}
\end{equation}

\begin{lemma} \label{DeltaEstimateForYi}
Let $R$ be a Banach $\C$-algebra with multiplicative norm $\abs{\;}$, and let $R^+\subset R$ be $\OO_{\C}$-subalgebra of elements $f$ with $\abs{f}\leq 1$.  Assume that $R^+$ is bounded (and thus complete for the $\pi$-adic topology).  Let $\Y_1,\Y_2\in \tilde{G}(R^+)$, and let $Y_1,Y_2\in R^+$ be the topologically nilpotent elements which correspond to $\Y_1,\Y_2$.  Finally, let $\lambda\in R^+$ be the topologically nilpotent element corresponding to the determinant $\lambda(Y_1,Y_2)\in \widetilde{\wedge^2 G}(R^+)$.
\begin{enumerate}
\item Suppose that $\abs{Y_2}^q\leq \abs{Y_1}<\abs{Y_2}$.  Then $\lambda=Y_1Y_2^q$ plus smaller terms.
\item Suppose that $\abs{Y_1}=\abs{Y_2}$.  Then $\lambda=Y_1Y_2^q-Y_1^qY_2$ plus smaller terms.
\item In general, there exists a unique $m\in\Z$ such that one of the following inequalities holds:
\begin{enumerate}
\item $\abs{Y_2}^q\leq \abs{Y_1}^{q^{-2m}}<\abs{Y_2}$
\item $\abs{Y_1}^q\leq \abs{Y_2}^{q^{-2m}}<\abs{Y_1}$
\item $\abs{Y_1}=\abs{Y_2}^{q^{-2m}}$.
\end{enumerate}
Then respectively we have
\begin{enumerate}
\item $\lambda=(-1)^mY_1^{q^m}Y_2^{q^{1-m}}$
\item $\lambda=(-1)^{m+1}Y_1^{q^{1-m}}Y_2^{q^m}$
\item $\lambda=(-1)^m(Y_1^{q^m}Y_2^{q^{1-m}}-Y_1^{q^{m+1}}Y_2^{q^{-m}})$
\end{enumerate}
plus smaller terms.
\end{enumerate}
\end{lemma}

\begin{proof}   To prove the lemma we will show that $F(Y_1,Y_2)-R(Y_1,Y_2)$ equals the desired approximation plus strictly smaller terms.  The same arguments will apply to $F(Y_1^{1/q^m},Y_2^{1/q^m})-R(X_1^{1/q^m},X_2^{1/q^m}))$.  Then Eq. \eqref{DeltaLimit} will show that $\delta(X,Y)$ equals the desired approximation plus strictly smaller terms.

We have $\lambda=\delta(Y_1,Y_2)$.  If $\abs{Y_1}=\abs{Y_2}$, then $\abs{R(Y_1,Y_2)}<\abs{F(Y_1,Y_2)}$ and we get claim (2).  If $\abs{Y_2}^q\leq\abs{Y_1}<\abs{Y_2}$, then $F(Y_1,Y_2)=Y_1Y_2^q$ plus smaller terms.  Since $R(X_1,0)=R(0,X_2)=0$ (this follows from the same properties of $\delta_0$ and $F$) we observe that every term of $R(Y_1,Y_2)$ contains both $Y_1$ and $Y_2$, and therefore (since $R$ has degree $\geq q^2$) we have a strict inequality $\abs{R(Y_1,Y_2)}< \abs{Y_2}^{q^2}$.  This is bounded by $\abs{Y_2}^{q^2}=\abs{Y_2}^q\abs{Y_2}^{q^2-q}\leq \abs{Y_1Y_2^{q^2-q}}\leq \abs{Y_1Y_2^q}=\abs{F(Y_1,Y_2)}$, which establishes claim (1).

In claim (3), the existence and uniqueness of $m$ is easy to see.  In the first case, apply claim (1) to the pair $\Y_1,\pi^{-m}\Y_2$ and note that $\lambda(\Y_1,\Y_2)=\pi^m\lambda(\pi^{-m}\Y_1,\Y_2)$.  For the second case, apply claim (1) to the pair $\Y_2,\pi^{-m}\Y_1$ (and recall that $\lambda$ is alternating).  For the third case, apply claim (2) to $\Y_1,\pi^{-m}\Y_2$.
\end{proof}

\subsection{Definition of the affinoids $\YY_{\x,m}$}
\label{DefinitionOfAffinoids}
Let $B$ be the coordinate ring of the affine formal scheme $\tilde{G}_{\OO_\C}\times\tilde{G}_{\OO_\C}$, so that $\tilde{G}_{\OO_\C}\times\tilde{G}_{\OO_\C}=\Spf B$.  We have two distinguished elements $\X_1,\X_2\in\tilde{G}(B)$, corresponding to the two projections $\tilde{G}\times\tilde{G}\to \tilde{G}$.  Let $X_i=[\X_i]$ be their coordinates;  then we have
\[ B\isom \OO_\C\powerseries{X_1^{1/q^\infty},X_2^{1/q^\infty}}.\]

Recall that $\GL_2(K)\times D^\times$ acts on the right of $\tilde{G}\times\tilde{G}$.  By our conventions, an element $g=(g_1,g_2)$ acts by the rule
\[ (g(X_1),g(X_2))=(g_2^{-1}X_1,g_2^{-1}X_2)g_1. \]

Let $(\x_1,\x_2)\in \tilde{G}(\OO_\C)\times\tilde{G}(\OO_\C)$ denote the image of $\x$ under the morphism $\MM_{\infty}(\OO_\C)\to \tilde{G}(\OO_\C)\times\tilde{G}(\OO_\C)$.
Then there exists a basis $\alpha_1,\alpha_2$ for $L/K$ and an element $\x_0\in \tilde{G}(\OO_\C)$ for which $\x_i=\alpha_i\x_0$.  Let $A\in \GL_2(L)$ be the matrix
\[ A=\tbt{\alpha_1}{\alpha_2}{\alpha_1^{\sigma}}{\alpha_2^{\sigma}}, \]
where $\sigma$ denotes the nontrivial automorphism of $L/K$.  Then we have $(\x_0,0)A=(\x_1,\x_2).$

Recall the elements $\varpi_1\in\A_\x$ and $\varpi_2\in \OO_D$ from \S\ref{linkingorders}
\begin{lemma}
\label{PropertiesOfA}  
\begin{enumerate}
\item Let $\alpha\in L^\times$, and let $g$ be the image of $\alpha$ in $\GL_2(K)$.  Then $AgA^{-1}=\tbt{\alpha}{}{}{\alpha^\sigma}$ (equality in $\GL_2(L)$).
\item We have $A\varpi_1 A^{-1}=\tbt{}{1}{1}{}$ (equality in $\GL_2(L)$).
\item We have $A\varpi_2 A^{-1}=\tbt{}{1}{1}{}\varpi_2$ (equality in $\GL_2(D)$, with $\varpi_2$ considered as a scalar matrix).
\end{enumerate}
\end{lemma}

\begin{proof}  By the definition of the embedding of $L^\times$ into $\GL_2(K)$, we have $(\x_1,\x_2)g=(\alpha\x_1,\alpha\x_2)$, so that $(\alpha_1,\alpha_2)g=(\alpha\alpha_1,\alpha\alpha_2)$.  Applying $\sigma$, we see that $(\alpha_1^\sigma,\alpha_2^{\sigma})g=(\alpha^\sigma\alpha_1^\sigma,\alpha^\sigma\alpha_2^\sigma)$, and so
\[ Ag=\tbt{\alpha_1}{\alpha_2}{\alpha_1^\sigma}{\alpha_2^{\sigma}}g=\tbt{\alpha\alpha_1}{\alpha\alpha_2}{\alpha^\sigma\alpha_1^\sigma}{\alpha^\sigma}{\alpha_2^\sigma}
=\tbt{\alpha}{}{}{\alpha^\sigma}A, \]
proving (1).

Similarly, (2) and (3) follow from $(\alpha_1,\alpha_2)\varpi_1=(\alpha_1^\sigma,\alpha_2^\sigma)$ and $(\alpha_1,\alpha_2)\varpi_2=\varpi_2(\alpha_1^{\sigma},\alpha_2^{\sigma})$.
\end{proof}

Define elements $\Y_1,\Y_2\in\tilde{G}(B)$ through an affine change of variables
\begin{equation}
\label{DefinitionOfYi}
(\X_1,\X_2) = (\x_1,\x_2) + (Y_1,Y_2)A.
\end{equation}

Let $Y_i=[\Y_i]$, $i=1,2$.  Also let $x_0=[\x_0]$.

\begin{defn} The affinoid $\YY_{\x,m}$ is defined by the inequalities
\[ \abs{Y_i}\leq \abs{x_0}^{s_i},\; i=1,2,\]
where $s_1$ and $s_2$ are defined by the following table.
\begin{center}
\begin{tabular}{l|ll}
& $s_1$ & $s_2$ \\
\hline
$L/K$ unramified & $q^{2m}$ & $q^m$ \\
$L/K$ ramified, $m$ even & $q^m$ & $q^{m/2}$\\
$L/K$ ramified, $m$ odd & $q^{m}$ & $\frac{q+1}{2}q^{(m+1)/2}$
\end{tabular}
\end{center}
\end{defn}
Thus if we let $Z_i=Y_i/x_0^{s_i}$, then $\YY_{\x,m}=\Spa(R,R^+)$, where
\[ R^+=\OO_\C\tatealgebra{Z_1^{1/q^\infty},Z_2^{1/q^\infty}}. \]  The reduction of $\YY_{\x,m}$ is then  $\overline{\YY}_{\x,m}=\Spec \overline{k}[Z_1^{1/q^\infty},Z_2^{1/q^\infty}]=\mathbf{A}^{2,\perf}_{\overline{k}}$.

We intend to prove Thm. \ref{ExistenceOfAffinoid} for the affinoid $\ZZ_{\x,m}=\YY_{\x,m}\cap \MM_{\infty,\overline{\eta}}^{\circ,\ad}$, under the assumption that $m$ is odd if $L/K$ is ramified.  The calculations (presented in \S\ref{CaseOfRamified}-\S\ref{CaseOfLevelZero}) follow the same general pattern.  First we verify that $\YY_{\x,m}$ is stabilized by the group $\K_{\x,m}$, and we compute the action of $\K_{\x,m}$ on $\overline{\YY}$ in terms of the variables $Z_1,Z_2$.  In each case the formulas are identical to the formulas given in Defn. \ref{DefnOfCxM}.

We then prove the existence of the claimed map $\overline{\ZZ}_{\x,m}\to C_{\x,m}^{\perf}$.  Recall the determinant morphism $\tilde{G}\times\tilde{G}\to \tilde{\wedge^2G}$ from \S\ref{detmodules}.  After choosing coordinates on the formal groups $G$ and $\wedge^2 G$, this morphism corresponds to an element $\delta(X_1,X_2)\in B$ admitting arbitrary $q$th power roots.  Let $t=\delta(x_1,x_2)\in \OO_\C$.  The key calculation is an approximation for $\delta(X_1,X_2)$ in terms of the variables $Z_1,Z_2$ inside the ring $R^+$:
\begin{equation}
\label{DeterminantCongruence} \delta(X_1,X_2)\equiv t+t^{r}f(Z_1,Z_2)+\text{ smaller terms}
\end{equation}
where $r\geq 1$ and $f(Z_1,Z_2)$ a polynomial which (up to replacing the $Z_i$ with $q$th powers) equals the polynomial defining the curve $C_{\x,m}$ from Thm. \ref{CurvesRealizeJLC}.

By Thm. \ref{Cartesian}, $\MM_{\infty,\eta}^{\circ,\ad}$ is the fiber of the determinant map $\tilde{G}^{\ad}_{\overline{\eta}}\times \tilde{G}^{\ad}_{\overline{\eta}}\to \widetilde{\wedge G}^{\ad}_{\overline{\eta}}$ over $\t$.  Thus $\delta(X_1,X_2)=t$ on $\M^{\circ,\ad}_{\infty,\eta}$.  Combining this fact with Eq. \eqref{DeterminantCongruence} shows that $f(Z_1,Z_2)=0$ holds in the coordinate ring of the reduction $\overline{\ZZ}$.  Since $\ZZ$ is the spectrum of a perfect ring, there must exist a map $\ZZ\to C_{\x,m}^{\perf}$ as claimed.

It will be helpful to record the action of element of $\GL_2(K)\times D^\times$ on the elements $\Y_1,\Y_2\in \tilde{G}(B)$.

\begin{lemma}\label{GroupActionOnYi}
Suppose $g\in \GL_2(K)$.  Let $g(\Y_1),g(\Y_2)$ denote the images of $\Y_1,\Y_2\in \tilde{G}(B)$ under the automorphism $\tilde{G}(B)\to\tilde{G}(B)$ coming from the automorphism of $B$ induced by $g$.  Then
\[ (g(\Y_1),g(\Y_2))=(\x_0,0)AgA^{-1}-(\x_0,0)+(\Y_1,\Y_2)AgA^{-1}. \]

Now suppose $g\in D^\times$, and define $g(\Y_1),g(\Y_2)$ similarly.  Then
\[ (g(\Y_1),g(\Y_2))=(\x_0,0)Ag^{-1}A^{-1}-(\x_0,0)+(\Y_1,Y_2)Ag^{-1}A^{-1}. \]
\end{lemma}

\begin{proof}  Using Eqs. \eqref{DefinitionOfYi} and \eqref{GroupActionOnYi} we have
\[ (\Y_1,\Y_2)=(\X_1,\X_2)A^{-1}-(\x_1,\x_2)A^{-1}=(\X_1,\X_2)A^{-1}-(\x_0,0). \]
We have
\begin{eqnarray*}
(g(\Y_1),g(\Y_2)) &=& (g(\X_1),g(\X_2))A^{-1}-(\x_0,0) \\
&=& (\X_1,\X_2)gA^{-1}-(\x_0,0) \\
&=& ((\x_1,\x_2)+(Y_1,Y_2)A)gA^{-1}-(\x_0,0) \\
&=& ((\x_0,0)A+(Y_1,Y_2)A)gA^{-1}-(\x_0,0) \\
&=& (\x_0,0)AgA^{-1} -(\x_0,0)+(Y_1,Y_2)AgA^{-1}.
\end{eqnarray*}
The case of $g\in D^\times$ is done the same way; the sign $g^{-1}$ appears because of our convention concerning the right action of $D^\times$ on $\tilde{G}\times\tilde{G}$.
\end{proof}

\subsection{Case: $L/K$ ramified, $m$ odd}
\label{CaseOfRamified} Let $\pi_L$ be a uniformizer of $L$.  Then $\pi_L$ is also a uniformizer of $D$.   We may assume $\pi_L^2\in \OO_K$, so that $\pi_L^\sigma=-\pi_L$.  After replacing $\x$ with a $\GL_2(K)$-translate we may also assume that the basis $\alpha_1,\alpha_2$ for $L/K$ is $1,\pi_L$.  Recall that $\pi_D\in D$ is the Frobenius endomorphism of $G_0$;  let us write $\pi_L=u\pi_D$ for some $u\in\OO_D^\times$.

We have the orthogonal decompositions
\begin{eqnarray*}
\A_\x&=&\OO_L\oplus \OO_L\varpi_1\\
\OO_D&=&\OO_L\oplus \OO_L\varpi_2
\end{eqnarray*}
where $\varpi_1^2=1$ and $\varpi_2^2\in\OO_K^\times$.  The linking order is
\[ \LL_{\x,m}=\Delta_{\x}(\OO_L)+(\gp_L^m\times\gp_L^m)+\left(\gp_L^{\frac{m+1}{2}}\varpi_1\times\gp_L^{\frac{m+1}{2}}
\varpi_2\right). \]
The affinoid $\YY_{\x,m}=\Spa(R,R^+)$ is defined by the conditions
\begin{eqnarray*}
\abs{\Y_1}&\leq& \abs{\x_0}^{q^m}\\
\abs{\Y_2}&\leq& \abs{\x_0}^{(q+1)q^{(m+1)/2}},
\end{eqnarray*}
where $\Y_1,\Y_2$ are defined as in Eq. \ref{DefinitionOfYi}.  We write $Y_i=[\Y_i]$ for the coordinate of $\Y_i$ in $R^+$, and define $Z_1,Z_2\in R^+$ by
\begin{eqnarray*}
Y_1&=&x_0^{q^m}Z_1^{q^{(1-m)/2}}\\
Y_2&=&x_0^{\frac{q+1}{2}q^{(m+1)/2}} Z_2.
\end{eqnarray*}
Then $R=\C\tatealgebra{Z_1^{1/q^\infty},Z_2^{1/q^\infty}}$.  Then $R$ is a Banach $\C$-algebra for the sup norm, which we write as $\abs{\;}$.  Note that $\abs{\;}$ is multiplicative (this is the perfectoid version of Gauss' lemma).

We wish to check that $\K_{\x,m}$ preserves $\YY_{\x,m}$, and to calculate the induced action of $\K_{\x,m}$ on the reduction $\overline{\YY}_{\x,m}=\mathbf{A}^{2,\perf}_{\overline{k}}$.  The group $\K_{\x,m}$ is generated by three types of elements:

\begin{enumerate}
\item $\Delta_{\x}(\alpha)$, for $\alpha\in L^\times$.
\item Elements of $\GL_2(K)$ of the form $1+\pi_L^m\beta+\pi_L^{(m+1)/2}\gamma\varpi_1$, with $\beta,\gamma\in\OO_L$.
\item Elements of $D^\times$ of the form $1+\pi_L^{(m+1)/2}\beta\varpi_2$, with $\beta\in\OO_L$.
\end{enumerate}

Let $\alpha\in L^\times$, let $g=\Delta_{\x}(\alpha)$, and let $g_1$ be the image of $\alpha$ in $\GL_2(K)$.   We have $\x^g=\x$.  By Lemmas \ref{GroupActionOnYi} and \ref{PropertiesOfA}, the action of $g$ on the variables $Y_1,Y_2$ is given by
\begin{eqnarray*}
(g(\alpha)(\Y_1),g(\alpha)(\Y_2))
&=&(\alpha^{-1}\Y_1,\alpha^{-1}\Y_2)Ag_1A^{-1}\\\\
&=&(\alpha^{-1}\Y_1,\alpha^{-1}\Y_2)\tbt{\alpha}{}{}{\alpha^{\sigma}}\\
&=&(\Y_1,\alpha^{\sigma}/\alpha \Y_2).
\end{eqnarray*}
Since $\alpha^{\sigma}/\alpha$ is always a unit in $\OO_L^\times$, it is clear that $g$ preserves $\YY_{\x,m}$.   If $\alpha\in\OO_K^\times$, then $\Delta_{\x}$ acts trivially on $\overline{\YY}_{\x,m}$, whereas $\Delta_{\x}(\pi_L)$ acts as $(Z_1,Z_2)\mapsto (Z_1,-Z_2)$.

Now suppose $g=1+\pi_L^m\beta+\pi_L^{(m+1)/2}\gamma\varpi_1\in\GL_2(K)$, where $\beta$ and $\gamma$ are to be interpreted as lying in the image of $\OO_L$ in $M_2(L)$.  By Lemma \ref{PropertiesOfA} we have $AgA^{-1}=1+\pi^m\beta+\pi^{(m+1)/2}\tbt{0}{\gamma}{\gamma^\sigma}{0}$, where now $\beta$ and $\gamma$ are to be interpreted as scalars in $L^\times\subset \GL_2(L)$.  Therefore by Lemma \ref{GroupActionOnYi} we have
\begin{eqnarray*}
g(\Y_1)&=&\Y_1+\pi^m\beta\Y_1+\pi_L^{(m+1)/2}\gamma^\sigma\Y_2+\pi^m\beta \x_0 \\
g(\Y_2)&=&\Y_2+\pi^m\beta\Y_2+\pi_L^{(m+1)/2}\gamma\Y_1+\pi^{(m+1)/2}\gamma \x_0.
\end{eqnarray*}

Taking coordinates and using Lemma \ref{ApproximateActionOnTildeG}, we find that the following equations hold modulo smaller terms:
\begin{eqnarray*}
g(Y_1)&=& Y_1+ \beta x_0^{q^{m}}\\
g(Y_2)&=& Y_2
\end{eqnarray*}

In terms of the coordinates $Z_1,Z_2$ on the reduction $\overline{\YY}_{\x,m}$, we have
\begin{eqnarray*}
g(Z_1)&=& Z_1+\overline{\beta}\\
g(Z_2)&=& Z_2
\end{eqnarray*}
The calculation for elements of the form $g=1+\pi^{(m+1)/2}\gamma\varpi_2$ is similar, with the result that such elements act trivially on $\overline{\YY}_{\x,m}$.

To complete the proof of Thm. \ref{ExistenceOfAffinoid} in this case, we need an approximation for the determinant $\lambda(\X_1,\X_2)$ as an element of the Banach algebra $R$.  The elements $\Y_1,\Y_2\in \tilde{G}(B)$ are defined by
\begin{eqnarray*}
\X_1&=&\x_0+\Y_1+\Y_2 \\
\X_2&=&\pi_L\x_0+\pi_L\Y_1-\pi_L\Y_2.
\end{eqnarray*}

Thus
\begin{eqnarray*}
\lambda(\X_1,\X_2)
&=&\lambda(\x_0,\pi_L\x_0)\\
&+&\lambda(\x_0,\pi_L\Y_1)+\lambda(\Y_1,\pi_L\x_0)\\
&-&\lambda(\Y_2,\pi_L\Y_2)\\
&+&\lambda(Y_2,\pi_L\x_0,\Y_2)-\lambda(\x_0,\pi_L\Y_2)\\
&-&\lambda(\Y_1,\pi_L\Y_2)+\lambda(\pi_L\Y_1,\Y_2)+\lambda(\Y_2,\pi_L \Y_2)
\end{eqnarray*}

We analyze each of these lines in turn to find an estimate for $[\lambda(X_1,X_2)]$ in the ring $R^+$.  The first line is $\lambda(\x_0,\pi_L\x_0)=\lambda(\x_1,\x_2)=\t$.  We have $[\t]=t$, $[\x_0]=x_0$, $[\pi_L\x_0]=ux_0^q$ up to smaller terms.  By Lemma \ref{DeltaEstimateForYi} we get that $ux_0^{2q}=t$ plus smaller terms.

The other lines can be treated using Lemma \ref{DeltaEstimateForYi}.  For instance, the first term on the second line is $\lambda(\x_0,\pi_L\Y_1)$.  We have $\abs{\pi_L\Y_1}=\abs{\Y_1}^q=\abs{\x_0}^{q^{m+1}}$.  Since $m+1$ is even, we are in the third case of claim (3) of Lemma \ref{DeltaEstimateForYi}, and so up to smaller terms we have
\begin{eqnarray*}
[\lambda(\x_0,\pi_L\Y_1)]
&=&(-1)^{m/2}\left([\x_0]^{q^{(m+1)/2}}[\pi_L\Y_1]^{q^{(3-m)/2}}-[\x_0]^{q^{(m+3)/2}}[\pi_L\Y_1]^{q^{(1-m)/2}}\right)  \\
&=&(-1)^{m/2}\left(x_0^{q^{(m+1)/2}}(ux_0^{q^{m+1}}Z_1^{q^{(1-m)/2}})^{q^{(3-m)/2}} - x_0^{q^{(m+3)/2}}(ux_0^{q^{m+1}}Z_1^{q^{(1-m)/2}})^{q^{(1-m)/2}}\right) \\
&=&(-1)^{m/2}ux_0^{(q+1)q^{(m+1)/2}}(Z_1^q-Z)\\
&=&(-1)^{m/2}t^{q^{(m+1)/2}}(Z_1^q-Z_1).
\end{eqnarray*}
The second term of the second line equals the first.

For the third line, we have $\abs{\pi_L\Y_1}=\abs{\Y_1}^q$, and we are in the situation of claim (1) of Lemma \ref{DeltaEstimateForYi}.  Up to smaller terms we have
\begin{eqnarray*}
\lambda(\Y_2,\pi_L\Y_2)
&=& -[\Y_2]^q[\pi_L\Y_2] \\
&=& -(x_0^{\frac{q+1}{2}q^{(m-1)/2}}Z_2^{1/q})^q(ux_0^{\frac{q+1}{2}q^{(m-1)/2}}Z_2) \\
&=& -ux_0^{(q+1)q^{(m+1)/2}}Z_2^2 \\
&=& -t^{q^{(m+1)/2}} Z_2^2
\end{eqnarray*}

The fourth line is zero on the nose, and the contribution of the fifth line is smaller than the contributions of the previous lines.  We get (up to a benign change of the variables $Z_1$, $Z_2$)
\[ \delta(Y_1,Y_2)=t+t^{q^{(m+1)/2}}(Z_1^q-Z_1-Z_2^2) \]
up to smaller terms.  Thus we have proved Prop. \ref{ExistenceOfAffinoid} in this case.

\subsection{Case: $L/K$ unramified, $m$ odd}
Once again we have the orthogonal decompositions
\begin{eqnarray*}
\A_{\x}&=&\OO_L\oplus\OO_L\varpi_1\\
\OO_D &=& \OO_L\oplus \OO_L\varpi_2,\\
\end{eqnarray*}
but this time $\varpi_2$ is a uniformizer in $D$.  The linking order is
\[ \LL_{\x,m} = \Delta_{\x}(\OO_L)+(\gp_L^m\times\gp_L^m)
+\left(\gp_L^{\frac{m+1}{2}}\varpi_1 \times \gp_L^{\frac{m-1}{2}}\varpi_2\right). \]
Recall that $\K_{\x,m}=\Delta_{\x}(L^\times)\L_{\x,m}^\times$.

The elements $\Y_1,\Y_2\in \tilde{G}(B)$ are defined by
\begin{eqnarray}
\label{XsInTermsOfYs}
\X_1 &=& \alpha_1\x_0+\alpha_1 \Y_1 +\alpha_1^{\sigma} \Y_2 \\
\X_2 &=& \alpha_2\x_0+\alpha_2 \Y_1 +\alpha_2^{\sigma} \Y_2,
\end{eqnarray}
and $\YY_{\x,m}$ is defined by the conditions
\begin{eqnarray*}
\abs{\Y_1}&\leq& \abs{x_0}^{q^{2m}}\\
\abs{\Y_2}&\leq& \abs{x_0}^{q^m}
\end{eqnarray*}

Let $Y_i=[\Y_i]$, and define variables $Z_1,Z_2$ by
\begin{eqnarray*}
Y_1&=&x_0^{q^{2m}}Z_1^{q^{m}}\\
Y_2&=&x_0^{q^m}Z_2
\end{eqnarray*}

We wish to check that $\K_{\x,m}$ preserves $\YY_{\x,m}$, and to calculate the induced action of $\K_{\x,m}$ on the reduction $\overline{\YY}_{\x,m}=\mathbf{A}^{2,\perf}_{\overline{k}}$.  The group $\K_{\x,m}$ is generated by three types of elements:
\begin{enumerate}
\item $\Delta_{\x}(\alpha)$, for $\alpha\in L^\times$,
\item Elements of $g\in D^\times$ with $g^{-1}=1+\pi^m\gamma+\pi^{(m-1)/2}\beta\varpi$, with $\beta,\gamma\in\OO_L$.
\item Elements of $g\in \GL_2(K)$ of the form $1+\pi^{(m+1)/2}\beta\sigma$, with $\beta\in \OO_L$.
\end{enumerate}

Let $\alpha\in L^\times$, and let $g=\Delta_{\x}(\alpha)$.  As in the previous case we have $g(\Y_1)=\Y_1$ and $g(\Y_2)=\alpha^{\sigma}/\alpha \Y_2$.  Since $\alpha^{\sigma}/\alpha$ is always a unit in $\OO_L^\times$, it is clear that $g$ preserves $\YY_{\x,m}$.   We have that $\Delta_{x}(\pi)$ acts trivially, and if $\alpha\in\OO_L^\times$, then the action of $g=\Delta_{\x}(\alpha)$ on the reduction $\overline{\YY}_{\x,m}$ is given in terms of the coordinates $Z_1,Z_2$ by $g(Z_1)=Z_1$ and $g(Z_2)=\overline{\alpha}^{q-1}Z_2$.

Now suppose $g^{-1}=1+\pi^m\gamma+\pi^{(m-1)/2}\beta\varpi_2$, with $\beta,\gamma$ lying in the image of $\OO_L$ in $M_2(\OO_K)$.   By Lemma \ref{PropertiesOfA} we have $Ag^{-1}A^{-1}=1+\pi^m\gamma+\pi^{(m-1)/2}\beta\tbt{0}{1}{1}{0}\varpi_2$, where now $\beta,\gamma$ are to be interpreted as scalars in $L\subset D\subset M_2(D)$.  Therefore by Eq. \eqref{GroupActionOnYi} we have
\begin{eqnarray*}
g(\Y_1)&=& \Y_1+ \pi^m\gamma \Y_1+\pi^{(m-1)/2}\beta\varpi_2\Y_2 +\pi^m\gamma \x_0 \\
g(\Y_2)&=& \Y_2+\pi^m\gamma \Y_2+\pi^{(m-1)/2}\beta\varpi_2 \Y_1 +\pi^{(m-1)/2}\beta \varpi_2\x_0.
\end{eqnarray*}
Taking coordinates and using Lemma \ref{ApproximateActionOnTildeG}, we find
\begin{eqnarray*}
g(Y_1)&=& Y_1+\beta Y_2^{q^m} + \gamma x_0^{q^{2m}}\\
g(Y_2)&=& Y_2+\beta x_0^{q^m}
\end{eqnarray*}
plus smaller terms.

In terms of the coordinates $Z_1,Z_2$ on the reduction $\overline{\YY}_{\x,m}$, we have
\begin{eqnarray*}
g(Z_1)&=& Z_1+\overline{\beta}^qZ_2+\overline{\gamma}^q\\
g(Z_2)&=& Z_2+\overline{\beta}
\end{eqnarray*}
(Here we have used $\beta^{q^m}=\beta^q$, because $m$ is odd.)  Note the accord with Defn. \ref{DefnOfCxM}. The calculation is similar for elements of $\GL_2(K)$ of the form $g=1+\pi^{(m+1)/2}\varpi_1$, with the result that such elements act trivially on $\overline{\YY}_{\x,m}$.

To complete the proof of Thm. \ref{ExistenceOfAffinoid} in this case, we need an approximation for the determinant $\lambda(\X_1,\X_2)$ as an element of the Banach algebra $R$.  From Eq. \eqref{XsInTermsOfYs} we find
\begin{eqnarray*}
\lambda(\X_1,\X_2)&=&\lambda(\alpha_1\x_0,\alpha_2\x_0)\\
&+&\lambda(\alpha_1\x_0,\alpha_2\Y_1)+\lambda(\alpha_1\Y_1,\alpha_2\x_0)\\
&+&\lambda(\alpha_1^\sigma\Y_2,\alpha_2^\sigma\Y_2)\\
&+&\lambda(\alpha_1\x_0,\alpha_2^{\sigma}\Y_2)+\lambda(\alpha_1^\sigma\Y_2,\alpha_2\x_0)\\
&+&\lambda(\alpha_1\Y_1,\alpha_2\Y_1)+\lambda(\alpha_1\Y_1,\alpha_2^{\sigma}  \Y_2)+\lambda(\alpha_2\Y_1,\alpha_1^{\sigma}\Y_2).
\end{eqnarray*}
We analyze each of these lines in turn.  The first line is $\lambda(\alpha_1\x_0,\alpha_2\x_0)=\t$, which by Prop. \ref{deltaestimate} implies that $(\alpha_1\alpha_2^q-\alpha_1^q\alpha_2)x_0^{q+1}=t$ plus smaller terms.

The other lines can be treated using Lemma \ref{DeltaEstimateForYi}.  For instance, the first term of the second line is $\lambda(\alpha_1\x_0,\alpha_2\Y_1)$.  We have $\abs{\alpha_2\Y_1}=\abs{\alpha_1\x_0}^{q^{2m}}$, so we are in the third case of part (3) of Lemma \ref{DeltaEstimateForYi}.  Thus up to smaller terms we have
\begin{eqnarray*} [\lambda(\alpha_1\x_0,\alpha_2\Y_1)]
&=&-([\alpha_1\x_0]^{q^m}[\alpha_2\Y_1]^{q^{1-m}}-[\alpha_1\x_0]^{q^{m+1}}[\alpha_2\Y_1]^{q^{-m}}] \\
&=&-((\alpha_1x_0)^{q^m}(\alpha_2x_0^{q^{2m}}Z_1^{q^{m}})^{q^{1-m}}-(\alpha_1x_0)^{q^{m+1}}(\alpha_2x_0^{q^{2m}}Z_1^{q^m})^{q^{-m}}\\
&=&-x_0^{(q+1)q^{m}}(\alpha_1^q\alpha^2 Z_1^q - \alpha_1\alpha_2^q Z_1) \\
\end{eqnarray*}
Similarly, the contribution of the second term of the second line is
\[
[\lambda(\alpha_1\Y_1,\alpha_2\x_0)]=x_0^{(q+1)q^{m}}(\alpha_1\alpha_2^q Z_1 - \alpha_1^q\alpha_2 Z_1^q),\]
so that the total contribution of the second line is
\[ (\alpha_1\alpha_2^q-\alpha_1^q\alpha_2)x_0^{(q+1)q^{m}}(Z_1^q+Z_1)=t^{q^{m}}(Z_1^q+Z_1). \]
The contribution of the third line is $t^{q^m}Z_2^{q+1}$ plus smaller terms.  The fourth line is zero on the nose, and the fifth line's contribution is smaller than any of the others.  We find that
\[ \delta(X_1,X_2)\equiv t+t^{q^m}(Z_1^q+Z_1-Z_2^{q+1})\]
modulo smaller terms.  This proves Prop. \ref{ExistenceOfAffinoid} in this case.

\subsection{Case:  $L/K$ unramified, $m\geq 2$ even}   This time the linking order is
\[ \L_{\x,m}=\Delta_{\x,m}(\OO_L)+(\gp_L^m\times\gp_L^m)\times \left(\gp_L^{m/2}\varpi_1\times\gp_L^{m/2}\varpi_2\right).\]
The affinoid $\YY_{\x,m}$ is defined the same way as in the previous section.  We must show that $\K_{\x,m}$ preserves $\YY_{\x,m}$ and compute its action on the reduction $\overline{\YY}_{\x,m}$.  The calculation is very similar to what occurred in the case of $m$ odd, except that in a sense the actions of $\GL_2(K)$ and $D^\times$ are reversed.

The group $\K_{\x,m}$ is generated by three types of elements:
\begin{enumerate}
\item $\Delta_{\x}(\alpha)$, for $\alpha\in L^\times$,
\item Elements of $\GL_2(K)$ of the form $1+\pi^m\gamma+\pi^{m/2}\beta\varpi_1$, with $\beta,\gamma\in \OO_L$.
\item Elements of $D^\times$ of the form $1+\pi^{m/2}\beta\varpi_2$, with $\beta\in\OO_L$.
\end{enumerate}
The action of $\Delta_{\x}(\alpha)$ works the same as in the case of $m$ odd.

Suppose $g=1+\pi^m\gamma+\pi^{m/2}\beta\varpi_1$, with $\beta,\gamma\in \OO_L$.  We have
\[ AgA^{-1}=1+\pi^m\gamma+\pi^{m/2}\tbt{0}{\beta}{\beta^\sigma}{0}. \]
By Lemma \ref{GroupActionOnYi} we get
\begin{eqnarray*}
g(\Y_1)&=&\Y_1+\pi^m\gamma(\x+\Y_1)+\pi^{m/2}\beta^{\sigma}\Y_2\\
g(\Y_2)&=&\Y_2+\pi^m\gamma\Y_2+\pi^{m/2}\beta\Y_1+\pi^{m/2}\beta(\x+\Y_1)
\end{eqnarray*}

Taking coordinates, we find
\begin{eqnarray*}
g(Y_1)&=&Y_1+\beta^{\sigma} Y_2^{q^m} + \gamma x_0^{q^{2m}}\\
g(Y_2)&=&Y_2+\beta x_0^{q^m}
\end{eqnarray*}
modulo smaller terms.  Recall that $Y_1=x_0^{q^{2m}}Z_1^{q^m}$, $Y_2=x_0^{q^m}Z_2$.  In terms of the coordinates $Z_1,Z_2$ on the reduction $\overline{\YY}_{\x,m}$, we have
\begin{eqnarray*}
g(Z_1)&=& Z_1+\overline{\beta}^qZ_2+\overline{\gamma}\\
g(Z_2)&=& Z_2+\overline{\beta}.
\end{eqnarray*}

The approximation for $\delta(X_1,X_2)$ in terms of $Z_1$ and $Z_2$ proceeds exactly as in the case of $m$ odd, with the result that
\[ \delta(X_1,X_2)=t+t^{q^m}(Z_1^q-Z_1-Z_2^{q+1}) \]
up to smaller terms.  This establishes Thm. \ref{ExistenceOfAffinoid} in this case.

\subsection{Case:  $L/K$ unramified, $m=0$}
\label{CaseOfLevelZero}
In this case the affinoid $\ZZ_{\x,0}$ we construct is related to the semistable model of the Lubin-Tate space of level 1 described by Yoshida in \cite{Yos}.

The linking order is $\L_{\x,0}=M_2(\OO_K)\times \OO_D$.  In this case $\YY_{\x,0}=\Spa(R,R^+)$ is the affinoid described by the conditions
\[ \abs{X_i}\leq \abs{x_0},\]
which is clearly preserved by $\K_{\x,0}=\Delta_{\x}(L^\times)\L_{\x,0}^\times$.

Let us write $X_i=x_0 Z_i$, and observe the effect of $\K_{\x,0}$ on the coordinates $Z_i$.  Suppose $g=(g_1,g_2)\in \L_{\x,0}^\times=\GL_2(\OO_K)\times \OO_D^\times$, with $g_1=\tbt{a}{b}{c}{d}$.  Then
\[ (g(\X_1),g(\X_2))= (g_2^{-1}\X_1,g_2^{-1}\X_2)\tbt{a}{b}{c}{d}.\]

In terms of the coordinates $Z_1,Z_2$, this means
\begin{eqnarray*}
g(Z_1) &\equiv & \overline{g}_2^{-1}(aZ_1+cZ_2) \\
g(Z_2) &\equiv & \overline{g}_2^{-1}(bZ_1+dZ_2).
\end{eqnarray*}

By Prop. \ref{deltaestimate} we have
\[ \delta(X_1,X_2)\equiv t(Z_1Z_2^q-Z_1^qZ_2)\]
up to smaller terms, thus establishing Thm. \ref{ExistenceOfAffinoid} in this case.

\section{Semistable coverings for the Lubin-Tate tower of curves}
\label{semistablecoverings}

\subsection{Generalities on semistable coverings of wide open curves}
The following notions are taken from~\cite{Coleman:StableMapsOfCurves}.  We will assume in the following that the field of scalars is $\C$.

\begin{defn}  A {\em wide open} (curve) is an adic space isomorphic to $C\backslash D$, where $C$ is the adic space attached to a smooth complete curve and $D\subset C$ is a finite disjoint union of closed discs.
Each connected component of $C$ is required to contain at least one disc from $D$.

If $W$ is a wide open, an {\em underlying affinoid} $Z\subset W$ is an open affinoid subset for which $W\backslash Z$ is a finite disjoint union of annuli $U_i$.  It is required that no annulus $U_i$ be contained in any open affinoid subset of $W$.

An {\em end} of $W$ is an element of the inverse limit of the set of connected components of $W\backslash Z$, where $Z$ ranges over open affinoid subsets of $W$.

Finally, $W$ is {\em basic} if it has an underlying affinoid $Z$ whose reduction $\overline{Z}$ is a semistable curve over $\overline{\FF}_q$.  (Recall that if $Z=\Spa(R,R^+)$, then $\overline{Z}=\Spec R^+\otimes_{\OO_\C}\overline{\FF}_q$.)
\end{defn}

For an affinoid $X$, there is a reduction map $\redu\from X\to\overline{X}$.  The following is a special case of Thm. 2.29 of \cite{ColemanMcMurdy}.
\begin{Theorem} \label{wideopencriterion}  If $X$ is a smooth one-dimensional affinoid, and $x$ is a closed point of $\overline{X}$, then the residue region $\redu^{-1}(x)$ is a wide open.
\end{Theorem}

In particular, $\M_{m,\eta}^{\circ,\ad}$ is a wide open, because it is the residue region $\redu^{-1}(x)$ of a supersingular point $x$ of the special fiber of an appropriate Shimura curve or Dinfeld modular curve.  

We adapt the definition of semistable covering in~\cite{Coleman:StableMapsOfCurves}, \S 2, which only applies to coverings of proper curves.  Our intention is to construct semistable coverings of the spaces $\M_{m,\eta}^{\circ,\ad}$.  Therefore we define:

\begin{defn}
\label{semistablecoveringW}
Let $W$ be a wide open curve.  A {\em semistable covering} of $W$ is a covering $\mathcal{D}$ of $W$ by connected wide opens satisfying the following axioms:
\begin{enumerate}
\item If $U,V$ are distinct wide opens in $\mathcal{D}$, then $U\cap V$ is a disjoint union of finitely many open annuli.
\item No three wide opens in $\mathcal{D}$ intersect simultaneously.
\item For each $U\in \mathcal{D}$, if
\[ Z_U = U\backslash\left(\bigcup_{U\neq V\in \mathcal{D}} V \right),  \]
then $Z_U$ is a non-empty affinoid whose reduction is nonsingular.
\end{enumerate}
In particular $U$ is a basic wide open and $Z_U$ is an underlying affinoid of $U$.
\end{defn}

Suppose $\mathcal{D}$ is a semistable covering of a wide open $W$.  For each $U\in\mathcal{D}$, let $\OO_W^+(U)$ be the ring of analytic functions on $U$ of norm $\leq 1$, and let $X_U=\Spf \OO_W^+(U)$.  Similarly if $U,V\in\mathcal{D}$ are overlapping wide opens, similarly define $X_{U\cap V}=\Spf \OO_W^+(U\cap V)$.  Let $\mathscr{W}$ denote the formal scheme over $\OO_{\C}$ obtained by gluing the $X_U$ together along the maps
\[ X_{U,V}\to X_U\coprod X_V. \]
Then $\mathscr{W}$ has generic fiber $W$.  The special fiber $\mathscr{W}_s$ of $\mathscr{W}$ is a scheme whose geometrically connected components are exactly the nonsingular projective curves $\overline{Z}_U^{\cl}$ with affine model $\overline{Z}_U$;  the curves $\overline{Z}_U^{\cl}$ and $\overline{Z}_V^{\cl}$ intersect exactly when $U$ and $V$ do.

\begin{exmp}
Note that $\mathcal{D}$ will in general not be finite. Suppose $W$ is the adic open unit disc over $\C$.
We construct a semistable covering $\mathcal{D}=\set{U_n}$ of $W$ indexed by integers $n\geq 0$.  First let $U_0=\set{\abs{z}<\abs{\pi^{1/2}}}$, and for $n\geq 1$ let $U_n=\set{\abs{\pi}^{1/n}<\abs{z}<\abs{\pi}^{1/(n+2)}}$.  Then $Z_{U_0}$ is the closed disc $\set{\abs{z}\leq \abs{\pi}}$ and for $n\geq 1$, $Z_{U_n}$ is the ``circle" $\set{\abs{z}=1/(n+1)}$.  The resulting formal scheme $\mathscr{W}$ has special fiber which is an infinite union of rational components;  the dual graph $\Gamma$ is a ray.
\end{exmp}

Let $C$ be the adic space attached to a smooth complete curve, let $D\subset C$ be a disjoint union of closed discs, and let $W=C\backslash D$.  A semistable covering of $W$ yields a semistable covering of $C$ (in the sense of~\cite{Coleman:StableMapsOfCurves}) by the following procedure.  Let $\mathcal{D}$ be a semistable covering of $W$ corresponding to the formal model $\mathscr{W}$.  Let $\Gamma$ be the dual graph attached to the special fiber of $\mathscr{W}$.  There are bijections among the following three finite sets:
\begin{enumerate}
\item ends of $W$,
\item ends of $\Gamma$, and
\item connected components of $D$.
\end{enumerate}
Suppose $v_1,v_2,\dots$ is a ray in $\Gamma$ corresponding to the wide opens $U_1,U_2,\dots\subset W$.  Then there exists $N>0$ such that for all $i\geq N$, $U_i$ is an open annulus.  If $D_0\subset D$ is the connected component corresponding to the ray $v_1,v_2,\dots$, then (possibly after enlarging $N$) $D_0\cup \bigcup_{i\geq N} U_i$ is an open disc, which intersects $U_{N-1}$ in an open annulus.  Repeating this process for all ends of $\Gamma$ yields a semistable covering $\mathcal{D}_0$ of $C$ by finitely many wide opens.  Let $\Gamma_0$ be the dual graph corresponding to $\mathcal{D}_0$.

In \cite{ColemanMcMurdy}, \S2.3, the genus $g(W)$ of a wide open curve $W$ is defined.  It is shown (Prop. 2.32) that in the above context that the genus of $W$ equals the genus (in the usual sense) of the smooth complete curve whose rigidification is $C$.  It is also shown (in the remark preceding Prop. 2.32) that if $U$ is a basic wide open whose underlying affinoid $Z_U$ has nonsingular reduction, then $g(U)=g(\overline{Z}_U^{\cl})$, where $\overline{Z}_U$ is the reduction of $Z_U$ and $\overline{Z}_U^{\cl}$ is the unique nonsingular projective curve containing $\overline{Z}_U$.  In Prop. 2.34 we find the formula
\begin{equation}
\label{prop234}
g(C)=\sum_{U\in\mathcal{D}_0} g(U)+\dim H^1(\Gamma_0,\Q).
\end{equation}

\begin{prop}\label{H1W}  Let $W$ be a wide open curve, and let $\mathcal{D}$ be a semistable covering of $W$.  Let $\Gamma$ be the dual graph of the special fiber of the corresponding semistable model $\mathscr{W}$ of $W$, so that the irreducible components $\overline{Z}_v^{\cl}$ of the special fiber of $\mathscr{W}$ are indexed by the vertices of $\Gamma$.  Then
\[
\dim H^1_c(W,\Q_\ell)=\sum_{v\in\Gamma} \dim H^1_c(\overline{Z}_v,\Q_\ell)+ \dim H^1_c(\Gamma,\Q_\ell)
\]
\end{prop}
\begin{proof}
Part of the long exact sequence in compactly supported cohomology for the pair $(C,D)$ reads
\begin{equation}
\label{H0CD}
 0\to H^0(C,\Q_\ell)\to H^0(D,\Q_\ell)\to H^1_c(W,\Q_\ell)\to H^1(C,\Q_\ell)\to 0.
\end{equation}
Note that $H^0_c(D,\Q_\ell)\isom \Q_\ell[\Ends(\Gamma)]$ is the space of $\Q_\ell$-valued functions on the set of ends of $\Gamma$.  On the other hand, $\Gamma_0$ is up to homotopy the graph obtained by deleting the ends from $\Gamma$, so that
\begin{equation}
\label{H1G}
\dim H^1_c(\Gamma,\Q_\ell)=\dim H^1(\Gamma_0,\Q_\ell)+\#\Ends(\Gamma)-\dim H^0(\Gamma,\Q_\ell).
\end{equation}
The result now follows from Eqs. \eqref{prop234}, \eqref{H0CD}, and \eqref{H1G}.
\end{proof}

\subsection{The fundamental domain}
\label{funddomainsection}
Let $\mathcal{F}\subset\MM^{\circ,\ad}_{\infty,\overline{\eta}}$ be the open affinoid subset defined by the conditions
\[ \abs{X_1}\geq \abs{X_{2}} \geq \abs{X_1}^q \]
\begin{prop}
\label{funddomain}
The translates of $\mathcal{F}$ under $(\GL_n(K)\times D^\times)^{\det=N}$ cover $\MM^{\circ,\ad}_{\infty,\overline{\eta}}$.
\end{prop}

\begin{proof} Let $\abs{\;}$ be a valuation on $A^{\circ}_{\OO_{\C}}$ with $\abs{\pi}\neq 0$.  In $A^{\circ}_{\OO_{\C}}$ we have the equation $\delta(X_1,X_2)=t$.  Since $\delta(X_1,X_2)$ has no constant term, and is alternating in its variables, we have $\abs{X_i}\neq 0$ for $i=1,2$.  On the other hand each $X_i$ is topologically nilpotent, so $\abs{X_i}<1$.

Certainly there exists $m\in \Z$ with $\abs{X_1}\geq \abs{X_2}^{q^m}>\abs{X_1}^q$.  Let $\pi_D\in \OO_D=\End G_0$ be the Frobenius element, so that $\pi_D$ acts on the variables $X_1,X_2$ by the rule $X_i\mapsto X_i^q$ plus smaller terms.  For the purposes of the lemma we may assume that $N(\pi_D)=\pi$.  If $m=2k$ is even, then the pair $\left(\tbt{1}{}{}{\pi^{-k}},\pi_D^{-k}\right)$ translates $\abs{\;}$ into $\mathcal{F}$.   If $m=2k+1$ is odd, then the pair $\left(\tbt{0}{1}{-\pi^{-k}}{0},\pi_D^{-k}\right)$ translates $\abs{\;}$ into $\mathcal{F}$.
\end{proof}

\subsection{A covering of the Lubin-Tate perfectoid space}
\label{graphT}
We now define a graph which is in a sense the dual graph for our semistable model of the Lubin-Tate tower.  We consider pairs $(\x,m)$, where $\x$ is a CM point in $\MM_{\infty,\overline{\eta}}^{\ad}$ and $m\geq 0$.

\begin{defn} Pairs $(\x,m)$ and $(\y,n)$ are {\em equivalent} if $m=n$ and there exists $g\in \mathcal{K}_{\x,m}^1$ for which $\y=\x^g$.  Call such a pair ramified or unramified as the CM field of $\x$ is ramified or unramified.   Also, we call $(\x,m)$ {\em imprimitive} if it is ramified and $m$ is even.  Otherwise, $(\x,m)$ is {\em primitive}.
\end{defn}
Define a graph $\T$ as follows:  the vertices are equivalence classes of pairs $(\x,m)$, and vertices $(\x,n)$ and $(\y,m)$ will be adjacent if (up to exchanging the pairs) one of the following conditions holds:
\begin{enumerate}
\item $m=n=0$, $\x$ is unramified, $\y$ is ramified, and $\A_\y\subseteq\A_\x$.
\item $(\y,m)$ is equivalent to $(\x,n+1)$.
\end{enumerate}
Then $\T$ admits a  action of $\GL_2(K)\times D^\times$.  To each vertex $v=(\x,m)$ of $\T$ we have an associated affinoid $\ZZ_{\x,m}$ which is open in the connected component of $\MM_{\infty,\overline{\eta}}^{\ad}$ containing $\x$, but not in $\MM_{\infty,\overline{\eta}}^{\ad}$ itself.  For the moment we work with one connected component:  Let $\T^\circ$ be the subgraph of $\T$ on the vertices $(x,m)$ with $x\in\MM_{\infty,\overline{\eta}}^{\circ,\ad}$, so that $T^\circ$ admits an action of $(\GL_2(K)\times D^\times)^{\det=N}$.

For each vertex $(\x,m)$ of $\T^{\circ}$, let
\[ \redu_{\x,m}\from \ZZ_{\x,m}\to\overline{\ZZ}_{\x,m} \]
be the reduction map, and let $S_{\x,m}\subset\overline{\ZZ}_{\x,m}$ be the set of images of CM points.

We now define an open cover $\set{W_v}$ of $\MM_{\overline{\eta}}^{\circ,\ad,\nonCM}$ indexed by vertices of $\T^\circ$.
If $v=(\x,0)$ for $\x$ unramified, assume that $\A_\x=M_2(\OO_K)$ and put
\[ W_v=\set{\abs{X_1}\geq \abs{X_2}>\abs{X_1}^q}\backslash\bigcup_{y\in S_{\x,0}} \redu_{\x,0}^{-1}(y), \]
If $v=(\x,0)$ for $\x$ ramified, assume that $\A_\x$ is the standard Iwahori algebra, and put
\[ W_v = \set{\abs{X_1}> \abs{X_2}\geq \abs{X_1}^q}\backslash\bigcup_{y\in S_{\x,0}}\redu_{\x,0}^{-1}(y).\]
If $v=(\x,m)$ for $m>0$ we set
\[ W_v = \redu^{-1}_{\x,m-1}(\overline{\x}) \backslash \bigcup_{y\in S_{\x,m}} \redu^{-1}_{\x,m}(y) \]
The assignment $v\mapsto \W_v$ can be extended to all vertices $v\in\T^{\circ}$ in such a way that $\W_v^g=\W_{v^g}$ for all $g\in (\GL_2(K)\times D^\times)^{\det=N}$.

\begin{prop}
The $W_v$ cover $\MM_{\overline{\eta}}^{\circ,\ad,\nonCM}$ (this being the complement in $\MM_{\infty,\overline{\eta}}^{\circ,\ad,\nonCM}$ of the set of CM points).
\end{prop}

\begin{proof} By Prop. \ref{funddomain} it suffices to show that the $W_v$ cover the set of non-CM points in $\mathcal{F}=\set{\abs{X_1}\geq \abs{X_2}\geq \abs{X_1}^q}$.  It is clear from the definitions of the $W_v$ that any point in $\mathcal{F}$ not lying in one of the $W_v$ must lie in $\cap_{m\geq 1}\ZZ_{\x,m}$ for some CM point $x$.  But $\cap_{m\geq 1}\ZZ_{\x,m}=\set{\x}$.
\end{proof}

Let
\[ \ZZ_v=W_v\backslash\bigcup_w W_w \]
where $w$ runs over vertices adjacent to $v$.  Then if $v=(\x,m)$:
\[ \ZZ_v = \ZZ_{\x,m} \backslash \bigcup_{y\in S_{(\x,m)}} \redu_{\x,m}^{-1}(y) \]
is the complement in $\ZZ_{\x,m}$ of finitely many residue regions.

\subsection{A semistable covering of $\MM_{m,\overline{\eta}}^{\ad}$.}
\label{SemistableCoveringAtFiniteLevel}

In this paragraph we translate our results about $\MM_{\infty,\overline{\eta}}^{\ad}$ into results about the Lubin-Tate spaces of finite level.  Recall the tower of complete local rings $A_m$, with $A$ defined as the completion of $\varinjlim A_m$.  Passing to adic spaces, we have a morphism from $\MM_{\infty,\eta}^{(0),\ad}$ to the projective system $\varprojlim \MM_{m,\eta}^{(0),\ad}$.

\begin{lemma}  \label{affinoidtoaffinoid} For each $m$, the morphism $\MM_{\eta}^{(0),\ad}\to\MM_{m,\eta}^{(0),\ad}$ is surjective, and carries open affinoids onto open affinoids.
\end{lemma}

\begin{proof}  The maps between the local rings $A_r$ are finite.  This implies that each continuous valuation on $A_m$ can be extended to $A_{m+r}$ for all $r\geq 0$, hence to $A$.  This shows that $\MM_{\eta}^{\ad}\to\MM_{\eta,m}^{\ad}$ is surjective.

Now suppose $\ZZ=\Spa(R,R^+)$ be an open affinoid in $\MM_{\infty,\eta}^{\ad}$.  Then $R^+=A\tatealgebra{f_1/\varpi,\dots,f_n/\varpi}$ for elements $f_1,\dots,f_n\in A$ generating an ideal of definition of $A$ and an element $\varpi\in\OO_K$ of positive valuation.   Since $\varinjlim A_r$ is dense in $A$, we may assume that the elements $f_i$ live in $A_N$ for some sufficiently large $N\geq m$.  Then the image of $\ZZ$ in $\MM_{N,\overline{\eta}}^{\ad}$ is an affinoid $\ZZ_N=\Spa(R_N,R_N^+)$, with $R_N=A_N\tatealgebra{f_1/\varpi,\dots,f_n/\varpi}$.  Since $\MM_{N,\overline{\eta}}^{\ad}\to\MM_{m,\overline{\eta}}^{\ad}$ is an \'etale map of adic curves, the image of $\ZZ_N$ in $\MM_{m,\eta}^{\ad}$ is again affinoid.
\end{proof}

Of course, Lemma \ref{affinoidtoaffinoid} holds for the tower of geometrically connected components $\MM_{m,\overline{\eta}}^{\circ,\ad}$ as well.

Fix $m\geq 0$.  For each vertex $v$ of $\T$, let $W_v^{(m)}$ be the image of $W_v$ in $\MM_{m,\overline{\eta}}$.   Similarly define $\ZZ_v^{(m)}$ as the image of $\ZZ_v$.  By repeatedly applying Lemmas \ref{wideopencriterion} and \ref{affinoidtoaffinoid}, we deduce that $W_v$ is a wide open.

Let $\Gamma^1(\pi^m)=\Gamma(\pi^m)\cap\SL_2(K)$.  We have a map $\overline{\ZZ}_v\to \overline{\ZZ}_v^{(m)}$.  Since $\overline{\ZZ}_v$ is the spectrum of a perfect ring, this map extends to a map $\overline{\ZZ}_v\to\overline{\ZZ}_v^{(m),\perf}$.

\begin{lemma} \label{quotient} Assume $m\geq 1$.  The map $\overline{\ZZ}_v\to\overline{\ZZ}_v^{(m),\perf}$ is a quotient by $\mathcal{K}_v^1\cap \Gamma^1(\pi^m)$.  That is, the coordinate ring of $\overline{\ZZ}_v^{(m),\perf}$ is the ring of $\mathcal{K}_v^1\cap \Gamma^1(\pi^m)$-invariants in the coordinate ring of $\overline{\ZZ}_v$.
\end{lemma}

\begin{proof}  Let $H=\mathcal{K}_v^1\cap \Gamma^1(\pi^m)$.  Let $S^{(m)}$ (resp. $S$) be the integral coordinate ring of $\ZZ_v^{(m)}$ (resp. of $\ZZ_v$), and let $\overline{S}^{(m)}$ (resp. $\overline{S}$) be its reduction.  It suffices to show that the map $\overline{S}^{(m),\perf}\to \overline{S}^{H}$ is surjective.  Let $\overline{f}\in \overline{S}$ be invariant by $H$, and let $f\in S$ be any lift.  We may assume that $f$ is invariant by $\Gamma^1(\pi^M)$ for some sufficiently large $M$, for the set of such elements (as $M$ varies) is dense in $\mathscr{S}$.  Let $H'=\mathcal{K}_v^1\cap \Gamma^1(\pi^M)$.  Let $g$ be the product of translates of $f$ by a set of coset representatives for $H/H'$, so that $g$ is $H$-invariant and therefore belongs to $S^{(m)}$.  Since $m\geq 1$, $H$ is a $p$-group, so that $[H:H']=p^n$ for some $n$.  Since $\overline{f}$ is $H$-invariant, we have $\overline{g}=\overline{f}^{p^n}$.  Thus $\overline{f}$ is the image of an element of $\overline{S}^{(m),\perf}$, namely $\overline{g}^{1/p^n}$.
\end{proof}

We now use the wide opens $W_v$ to construct an open covering of $\MM_{m,\overline{\eta}}^{\ad}$.  For a CM point $\x\in\MM^{\ad}_{\infty,\overline{\eta}}$, let $\x^{(m)}$ be the image of $\x$ in $\MM^{\ad}_{m,\overline{\eta}}$.  Let $U_\x$ be a sufficiently small affinoid neighborhood of $\x^{(m)}$, so that $U_\x$ is a disc.  Then $U_\x$ contains $\ZZ_{(\x,m)}$ for $m$ sufficiently large, say $m\geq N_{\x}+1$.  Now let $\T^{(m)}$ be the graph described by the following procedure:  Start with $\T/\Gamma(\pi^m)$, but remove $(\x,m)$ whenever $m> N_{\x}$.  Then we have a covering of $\MM_{m,\overline{\eta}}^{\ad}$ by wide opens $V_v$ indexed by the vertices of $\T^{(m)}$, where we have put
\[ V_{(\x,m)}=
\begin{cases}
W_{(\x,m)},& m<N_x \\
W_{(\x,m)}\cup U_x,& m=N_x.
\end{cases} \]

Recall that $\MM_{m,\overline{\eta}}^{\ad}/\pi^\Z$ is the quotient of $\MM_{m,\overline{\eta}}^{\ad}$ by the subgroup of $\GL_2(K)$ generated by the scalar $\pi$.   We get a get a covering of $\MM_{m,\overline{\eta}}^{\ad}/\pi^\Z$ by wide opens $V_v$ indexed by vertices of the quotient graph $\T^{(m)}/\pi^\Z$.  We will show this is a semistable covering.  For this we will have to show that the affinoids $\ZZ_v^{(m)}\subset V_v^{(m)}$ have good reduction.   A priori, we only know that the reduction $\overline{\ZZ}_v^{(m)}$ is an integral scheme over $\Spec \overline{k}$ of dimension 1.  Write $\overline{\ZZ}_v^{(m),\cl}$ for the smooth projective curve associated to the function field of $\overline{\ZZ}_v^{(m)}$.

\begin{prop} \label{ZvJL} We have
\[ \sum_{v\in \T^{(m)}/\pi^\Z}\dim H^1(\overline{\ZZ}_v^{(m),\cl},\Q_\ell) \geq
2\sum_{\Pi} \dim\Pi^{\Gamma(\pi^m)}\dim \JL(\check{\Pi}), \]
where the sum ranges over supercuspidal representations of $\GL_2(K)$ whose central character is trivial on $\pi$.
\end{prop}

\begin{proof}
For each vertex $v=(x,m)\in \T$, the perfection $\overline{\ZZ}_v^{(m),\perf}$ is the quotient of $\overline{\ZZ}_v$ by $\mathcal{K}^1_v\cap\Gamma^1(\pi^m)$, by Lemma \ref{quotient}.  On the other hand, by Thm. \ref{ExistenceOfAffinoid} there exists a nonconstant $\K_v^1$-equivariant morphism $\overline{\ZZ}_v\to C_v$, where $C_v$ is the curve of Defn. \ref{DefnOfCxM}.  From this we can conclude there exists a nonconstant morphism $\overline{\ZZ}_v^{(m)}\to C_v/(\K^1_v\cap \Gamma^1(\pi^m))$.  This morphism extends to a morphism of smooth projective curves $\overline{\ZZ}_v^{(m),\cl}\to C_v^{\cl}/(\K^1_v\cap\Gamma^1(\pi^m))$, so that
\[ \dim H^1(\overline{\ZZ}_v^{(m),\cl},\Q_\ell) \geq \dim H^1_c(C_v^{\cl},\Q_\ell)^{\K^1_v\cap\Gamma^1(\pi^m)}. \]
We now sum this inequality over all $v\in\T^{(m)}/\pi^\Z$.  Suppose $R\subset\T^{\circ}$ is a set of representatives for the quotient $\T^{\circ}/(\GL_2(K)\times D^\times)$.  Then every vertex in $\T^{(m)}/\pi^\Z$ is the translate of some uniquely determined $v\in R$ by an element $g\in \GL_2(K)\times D^\times$ which is well-defined up to left multiplication by $\K_v^1$ (the stabilizer of $v$) and up to right multiplication by $\Gamma^1(\pi^m)\pi^\Z$ (which fixes $\T^{(m)}$ pointwise).  We get
\begin{eqnarray*}
\sum_{v\in \T^{(m)}/\pi^\Z}\dim H^1(\overline{\ZZ}_v^{(m),\cl},\Q_\ell)
&\geq&
\sum_{v\in R}\;\;\sum_{g\in \K_v^1\backslash(\GL_2(K)\times D^\times)/\Gamma^1(\pi^m)\pi^\Z} \dim H^1_c(C_{v^g},\Q_\ell) \\
&=&\sum_{v\in R} \dim \left(\Ind_{\K_v^1}^{\GL_2(K)/\pi^\Z\times D^\times} H^1_c(C_v,\Q_\ell)\right)^{\Gamma^1(\pi^m)}
\end{eqnarray*}
by Mackey's theorem.  The result now follows from Prop. \ref{CurvesRealizeJLC}.
\end{proof}

\begin{Theorem} $\set{V_v}_{v\in \T^{(m)}}$ constitutes a semistable covering of $\MM^{\ad}_{m,\overline{\eta}}$.
\end{Theorem}

\begin{proof} It suffices to show that $\set{V_v}_{v\in \T^{(m)}/\pi^\Z}$ constitutes a semistable covering of $\MM^{\ad}_{m,\overline{\eta}}/\pi^\Z$.  Let us abbreviate $T=\T^{(m)}/\pi^\Z$ and $Z_v=\ZZ_v^{(m)}$ for $v\in T$.  The wide open curve $\MM^{\ad}_{m,\overline{\eta}}/\pi^\Z$ admits {\em some} semistable covering, so suppose there is a graph $T'$ and a collection of wide opens $V_v'$ satisfying the criteria in Defn. \ref{semistablecoveringW}, with underlying affinoids $Z_v$.  After refining the covering, we may assume that $T'$ contains $T$ as a subgraph, that $Z_v'\subset Z_v$ for all vertices $v\in T$, and that (for all $v\in T$) $\ZZ_v\subset\ZZ_v'$ is open, so that $\overline{\ZZ}_v^{\cl}=(\overline{\ZZ}'_v)^{\cl}.$

By Prop. \ref{H1W} we have
\[ \dim H^1_c(\MM^{\ad}_{m,\overline{\eta}}/\pi^\Z,\Q_\ell)=\sum_{v\in T'} \dim H^1_c(\overline{Z}_v',\Q_\ell) + \dim H^1_c(T',\Q_\ell). \]
On the other hand Cor. \ref{dimJL} gives
\[ \dim H^1_c(\MM^{\ad}_{\overline{\eta}}/\pi^\Z,\Q_\ell)=2\sum_\Pi \dim \Pi^{\Gamma(\pi^m)}\dim\JL(\check{\Pi})+2q^{m-1}(q-1)\dim \St^{\Gamma(\pi^m)} \]
We have $\dim\St^{\Gamma(\Pi^m)}=\#\mathbf{P}^1(\OO_K/\pi^m)-1$.  The dimension of $H^1_c(T',\Q_\ell)$ is at least $2q^{m-1}(q-1)(\#\mathbf{P}^1(\OO_K/\pi^m)-1)$, because $T$ has $2q^{m-1}(q-1)=\#K^\times/\pi^{2\Z}(1+\pi^m\OO_K)$ connected components, and each component has ends in correspondence with $\mathbf{P}^1(\OO_K/\pi^m)$.  This has the following consequences:
\begin{enumerate}
\item The inequality in Prop. \ref{ZvJL} is an equality,
\item For all primitive $v\in T$, the morphism $\overline{Z}_v^{(m)}\to C_v/(\K^1_v\cap \Gamma^1(\pi^m))$ induces an isomorphism on the level of $H^1_c$,
\item For all imprimitive $v\in T$, and all $v\in T'\backslash T$, $H^1_c(\overline{Z}_v',\Q_\ell)=0$, so that $\left(\overline{Z}_v'\right)^{\cl}=\mathbf{P}^1$, and
\item $H^1_c(T',\Q_\ell)=H^1_c(T,\Q_\ell)$, so that $T'$ is cycle-free and has no ends other than those of $T$.
\end{enumerate}
These imply that $\set{V_v}_{v\in T}$ was a semistable covering to begin with, because otherwise its semistable refinement would have introduced new curves of positive genus in the special fiber, or else monodromy in the dual graph.
\end{proof}




We can now complete the proof of Thm. \ref{MainTheoremInDepth}.  First, we observe that the wide opens $W_v$ constitute a semistable covering of $\MM_{m,\overline{\eta}}^{\ad,\nonCM}$, from which we get a compatible family of semistable models $\hat{\MM}^{\nonCM}_m$.  Each irreducible component of $\hat{\MM}^{\nonCM}_m$ is the completion of $\overline{\ZZ}_v^{(m)}$ for some $v\in \T/\Gamma(\pi^m)$.  The components corresponding to imprimitive vertices have genus 0.  Let $v$ be any primitive vertex of $\T$.  Then there exists $m$ large enough so that $\Gamma(\pi^m)\cap \K_{\x,v}^1$ acts trivially on $C_v$.  We have the morphism $\overline{\ZZ}_v^{(m)}\to C_v$;  since it induces an isomorphism on $H^1_c$ (see point (2) in the above proof), and since $C_v$ has positive genus, this morphism must be purely inseparable.  The same is true for $\overline{\ZZ}_v^{(m')}\to C_v$ for all $m'>m$;  since this map factors through $\overline{\ZZ}_v^{(m')}\to\overline{\ZZ}_v^{(m)}$, the latter must be purely inseparable as well.

\section{Stable reduction of modular curves:  Figures}
\label{figures}
In Figures 1-3, we draw the graph $\T^\circ$ constructed in \S\ref{graphT}.   Each vertex $v$ is labeled with its corresponding curve appearing on the list of four curves in Thm.~\ref{mainthm}.

We sketch a procedure for calculating the dual graph corresponding to the special fiber of a stable model of one geometrically connected component of the classical modular curve $X_m=X(\Gamma(p^m)\cap\Gamma_1(N))$, where $N\geq 5$.  First one must calculate the quotient $\T^\circ/\Gamma^1(p^m)$, where $\Gamma^1(p^m)=(1+p^mM_2(\Z_p))\cap\SL_2(\Q_p)$.  The image of a vertex $v$ in the quotient is labeled with the nonsingular projective curve constructed by quotienting $C_v$ by $\Gamma^1(p^m)\cap \mathcal{K}^1_v$.  For almost every $v$, the quotient is rational.  The quotient graph $\T^\circ/\Gamma^1(p^m)$ has finitely many ends, and each end is (once one goes far enough) a ray consisting only of rational components.  Erase all rational components lying on an end which corresponds to a CM point.  The remaining ends correspond to the boundary of $\MM^{\circ,\ad}_{m,\overline{\eta}}$;  these are in bijection with $\mathbf{P}^1(\Z/p^m\Z)$.  For each $b\in \mathbf{P}^1(\Z/p^m\Z)$, erase all rational components lying on the end corresponding to $b$, and let $v_b$ be the unique non-rational vertex which is adjacent to one of the vertices just erased.  Call the resulting graph $\T_m$.

Let $\Ig(p^m)$ denote the nonsingular projective model of the Igusa curve parameterizing elliptic curves over $\overline{\FF}_p$ together with Igusa $p^m$ structures and a point of order $N$.  Draw $\mathbf{P}^1(\Z/p^m\Z)$ many vertices $w_b$, and label each with $\Ig(p^m)$.  For each $b\in\mathbf{P}^1(\Z/p^m\Z)$, and each supersingular point of $X_1(N)(\overline{\FF}_p)$, attach a copy of $\T_m$ to $w_b$ in such a way that the vertex $w_b$ is adjacent to each $v_b$.  Finally, blow down any superfluous rational components.  The result is a finite graph representing the special fiber of a stable model of one component of $X_m$.


\begin{figure}
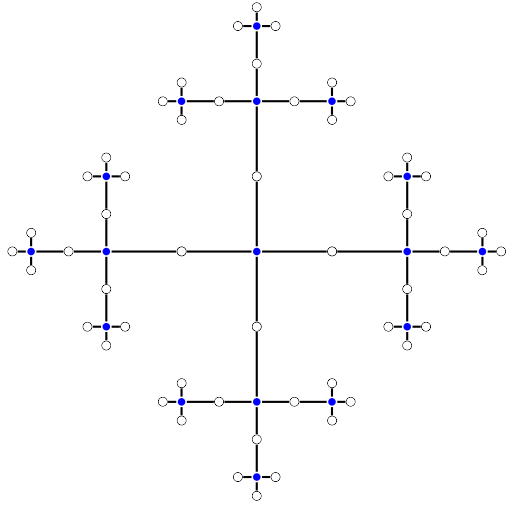

\begin{center}
\depthzerotree
\end{center}
\caption{The ``depth zero" subgraph of $\T^\circ$, consisting of the vertices $v=(\x,0)$.   The blue vertices are unramified.  Each represents a copy of the nonsingular projective curve with affine model $xy^q-x^qy=1$.  The stabilizer of any particular blue vertex in $\SL_2(K)$ is conjugate to $\SL_2(\OO_K)$.  The white vertices are imprimitive.  Each represents a rational component.  The stabilizer of any white vertex in $\SL_2(K)$ is an Iwahori subgroup.}
\end{figure}

\begin{figure}
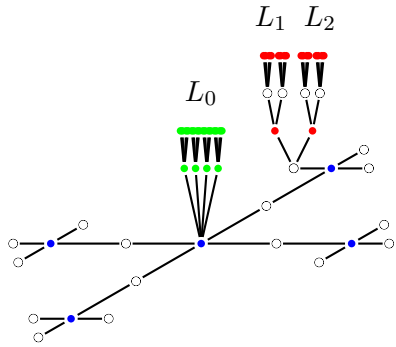

\begin{center}
\BTtreewithsprouts
\end{center}
\caption{Here the depth zero subgraph of $\Gamma$ is shown with several wild vertices to reveal structure.  The green vertices are unramified;  each represents a copy of the nonsingular (disconnected) projective curve with affine model $y^q+y=x^{q+1}$.  The red vertices are ramified;  each represents a copy of the nonsingular projective curve with affine model $y^q-y=x^2$.  The wild vertices $(\x,m)$ labeled with an $L_i$ are those for which $\x$ has CM by $L_i$.  Here $L_0,L_1,L_2$ are the three quadratic extensions of $L$, with $L_0/K$ unramified.  The white vertices are imprimitive;  each represents a rational component.}
\end{figure}

\begin{figure}
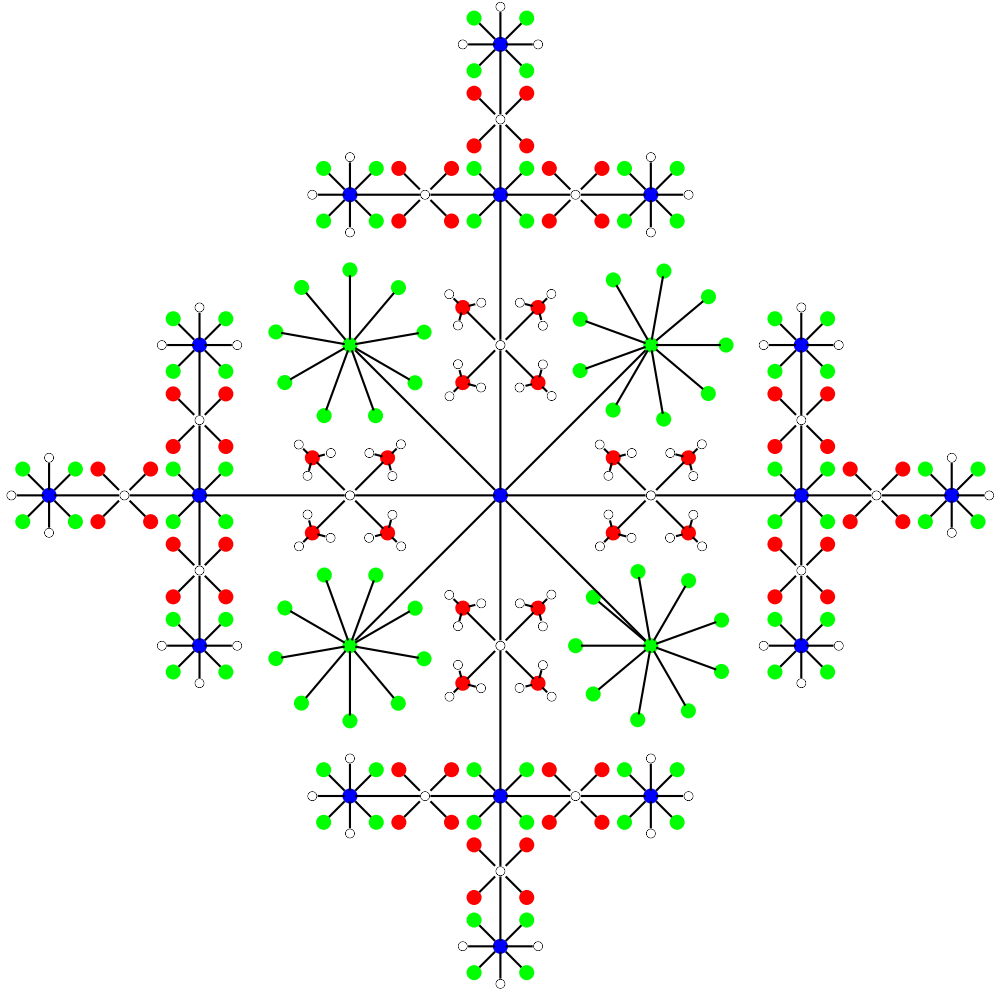

\pspicture(-7,-7)(7,7)
\psset{unit=1, nodesep=2pt}
\fullpicture
\endpspicture
\caption{Dual graph of the special fiber of our semistable model of the tower of Lubin-Tate curves:  complete picture.}
\end{figure}

\bibliographystyle{amsalpha}
\bibliography{bibfile}

\end{document}